\documentclass[11pt]{article}
\usepackage{amssymb}
\usepackage{amsmath}
\usepackage{times}
\usepackage[all]{xy}

\title{Injective tests of low complexity in the plane\indent}
\author{Dominique LECOMTE and Rafael ZAMORA$^1$}
\date{\today}

\def\ufootnote#1{\let\savedthfn\thefootnote\let\thefootnote\relax
\footnote{#1}\let\thefootnote\savedthfn\addtocounter{footnote}{-1}}

\newcommand{\Ana}{{\it\Sigma}^{1}_{1}}

\newcommand{\Ca}{{\it\Pi}^{1}_{1}}

\newcommand{\Borel}{{\it\Delta}^{1}_{1}}
\newcommand{\ana}{{\bf\Sigma}^{1}_{1}}

\newcommand{\Borone}{{\it\Delta}^{0}_{1}}
\newcommand{\boraone}{{\bf\Sigma}^{0}_{1}}
\newcommand{\boratwo}{{\bf\Sigma}^{0}_{2}}

\newcommand{\boraxi}{{\bf\Sigma}^{0}_{\xi}}

\newcommand{\bortwo}{{\bf\Delta}^{0}_{2}}

\newcommand{\bormone}{{\bf\Pi}^{0}_{1}}

\newcommand{\bormtwo}{{\bf\Pi}^{0}_{2}}

\newcommand{\bormlxi}{{\bf\Pi}^{0}_{<\xi}}

\newcommand{\bormxi}{{\bf\Pi}^{0}_{\xi}}

\newcommand{\borxi}{{\bf\Delta}^{0}_{\xi}}

\newtheorem{thm} {Theorem} [section]
\newtheorem{defi} [thm] {Definition}
\newtheorem{cor} [thm] {Corollary}
\newtheorem{lem} [thm] {Lemma}
\newtheorem{prop} [thm] {Proposition}

\begin{document}

\maketitle

\centerline{$\bullet$ Universit\' e Paris 6, Institut de Math\'ematiques de Jussieu, Projet Analyse Fonctionnelle}

\centerline{Couloir 16-26, 4\`eme \'etage, Case 247, 4, place Jussieu, 75 252 Paris Cedex 05, France}

\centerline{dominique.lecomte@upmc.fr}\bigskip

\centerline{Universit\'e de Picardie, I.U.T. de l'Oise, site de Creil,}

\centerline{13, all\'ee de la fa\"\i encerie, 60 107 Creil, France}\bigskip

\centerline{$\bullet^1$ Universit\' e Paris 6, Institut de Math\'ematiques de Jussieu, Projet Analyse Fonctionnelle}

\centerline{Couloir 15-16, 5\`eme \'etage, Case 247, 4, place Jussieu, 75 252 Paris Cedex 05, France}

\centerline{rafael.zamora@imj-prg.fr}\bigskip\bigskip\bigskip\bigskip\bigskip

\ufootnote{{\it 2010 Mathematics Subject Classification.}~Primary: 03E15, Secondary: 26A21, 54H05}

\ufootnote{{\it Keywords and phrases.}~acyclic, Borel, class, dichotomy, difference, homomorphism, injective, locally countable, oriented graph, reduction, Wadge}

\noindent {\bf Abstract.} We study injective versions of the characterization of sets potentially in a Wadge class of Borel sets, for the first Borel and Lavrentieff classes. We also study the case of oriented graphs in terms of continuous homomorphisms, injective or not.

\vfill\eject
 
\section{$\!\!\!\!\!\!$ Introduction}\indent
 
 The reader should see [K] for the standard descriptive set theoretic notation used in this paper. This work is a contribution to the study of analytic subsets of the plane. We are looking for results of the following form: either a situation is simple, or it is more complicated than a situation in a collection of known complicated situations. The notion of complexity we consider is the following, and defined in [Lo3].
 
\begin{defi} (Louveau) Let $X$, $Y$ be Polish spaces, $B$ be a Borel subset of $X\!\times\! Y$, and ${\bf\Gamma}$ be a class of Borel sets closed under continuous pre-images. We say that $B$ is \bf potentially in\it\ $\bf\Gamma$ $\big($denoted ${B\!\in\!\mbox{pot}(\bf{\Gamma})\big)}$ if there are finer Polish topologies $\sigma$ and $\tau$ on $X$ and $Y$, respectively, such that $B$, viewed as a subset of the product $(X,\sigma )\!\times\! (Y,\tau )$, is in $\bf\Gamma$.\end{defi}

 The quasi-order $\leq_B$ of Borel reducibility was intensively considered in the study of analytic equivalence relations during the last decades. The notion of potential complexity is a natural invariant for $\leq_B$: if $E\leq_BF$ and $F\!\in\!\mbox{pot}(\bf{\Gamma})$, then 
$E\!\in\!\mbox{pot}(\bf{\Gamma})$ too. However, as shown in [L1]-[L6] and [L8], $\leq_B$ is not the right notion of comparison to study potential complexity, in the general context, because of cycle problems. A good notion of comparison is as follows. Let $X,Y,X',Y'$ be topological spaces and $A,B\!\subseteq\! X\!\times\! Y$, $A',B'\!\subseteq\! X'\!\times\! Y'$. We write\bigskip

\leftline{$(X,Y,A,B)\leq (X',Y',A',B')\Leftrightarrow$}\smallskip

\rightline{$\exists f\! :\! X\!\rightarrow\! X'~~
\exists g\! :\! Y\!\rightarrow\! Y'\mbox{ continuous with }A\!\subseteq\! (f\!\times\! g)^{-1}(A')
\mbox{ and }B\!\subseteq\! (f\!\times\! g)^{-1}(B').$}\bigskip

 Our motivating result is the following (see [L8]).

\begin{defi} We say that a class $\bf\Gamma$ of subsets of zero-dimensional Polish spaces is a  \bf Wadge class of Borel sets\it\ if there is a Borel subset $\bf A$ of $\omega^\omega$ such that for any zero-dimensional Polish space $X$, and for any $A\!\subseteq\! X$, $A$ is in $\bf\Gamma$ if and only if there is $f\! :\! X\!\rightarrow\!\omega^\omega$ continuous such that 
$A\! =\! f^{-1}({\bf A})$. In this case, we say that $\bf A$ is $\bf\Gamma$-complete.\end{defi}

 If $\bf\Gamma$ is a class of sets, then $\check {\bf\Gamma}\! :=\!\{\neg A\mid A\!\in\! {\bf\Gamma}\}$ is the \bf dual class\rm\ of $\bf\Gamma$, and $\bf\Gamma$ is \bf self-dual\rm\ if 
${\bf\Gamma}\! =\!\check {\bf\Gamma}$. We set 
$\Delta ({\bf\Gamma})\! :=\! {\bf\Gamma}\cap\check {\bf\Gamma}$.

\begin{thm} \label{motivating} (Lecomte) Let $\bf\Gamma$ be a Wadge class of Borel sets, or the class $\borxi$ for some countable ordinal $\xi\!\geq\! 1$. Then there are concrete disjoint Borel relations $\mathbb{S}_0$, $\mathbb{S}_1$ on $2^\omega$ such that, for any Polish spaces $X,Y$, and for any disjoint analytic subsets $A,B$ of $X\!\times\! Y$, exactly one of the following holds:\smallskip

(a) the set $A$ is separable from $B$ by a $\mbox{pot}({\bf\Gamma})$ set,\smallskip

(b) $(2^\omega ,2^\omega ,\mathbb{S}_0,\mathbb{S}_1)\leq (X,Y,A,B)$.\end{thm}

 It is natural to ask whether we can have $f$ and $g$ injective if (b) holds. Debs proved that this is the case if $\bf\Gamma$ is a non self-dual Borel class of rank at least three (i.e., a class $\boraxi$ or $\bormxi$ with $\xi\!\geq\! 3$). As mentioned in [L8], there is also an injectivity result for the non self-dual Wadge classes of Borel sets of level at least three. Some results in [L4] and [L8] show that we cannot have $f$ and $g$ injective if (b) holds and $\bf\Gamma$ is a non self-dual Borel class of rank one or two, or the class of clopen sets, because of cycle problems again.\bigskip

 The work of Kechris, Solecki and Todor\v cevi\'c indicates a way to try to solve this problem. Let us recall one of their results in this direction. All the relations considered in this paper will be binary.
 
\begin{defi} Let $X$ be a set, and $A$ be a relation on $X$.\smallskip

(a) $\Delta (X)\! :=\!\{ (x,y)\!\in\! X^2\mid x\! =\! y\}$ is the \bf diagonal\it\ of $X$.\smallskip

(b) We say that $A$ is \bf irreflexive\it\ if $A$ does not meet $\Delta (X)$.\smallskip

(c) $A^{-1}\! :=\!\{ (x,y)\!\in\! X^2\mid (y,x)\!\in\! A\}$, and $s(A)\! :=\! A\cup A^{-1}$ is the 
\bf symmetrization\it\ of $A$.\smallskip

(d) We say that $A$ is \bf symmetric\it\ if $A\! =\! A^{-1}$.\smallskip

(e) We say that $A$ is a \bf graph\it\ if $A$ is irreflexive and symmetric.\smallskip

(f) We say that $A$ is \bf acyclic\it\ if there is no injective sequence $(x_i)_{i\leq n}$ of points of 
$X$ with $n\!\geq\! 2$, $(x_i,x_{i+1})\!\in\! A$ for each $i\! <\! n$, and $(x_n,x_0)\!\in\! A$.\smallskip

(g)  We say that $A$ is \bf locally countable\it\ if $A$ has countable horizontal and vertical sections (this also makes sense in a rectangular product $X\!\times\! Y$).\end{defi}
 
\noindent\bf Notation.\rm\ Let $(s_{n})_{n\in\omega}$ be a sequence of finite binary sequences with the following properties:\medskip

(a) $(s_n)_{n\in\omega}$ is \bf dense\rm\ in $2^{<\omega}$. This means that for each 
$s\!\in\! 2^{<\omega}$, there is $n\!\in\!\omega$ such that $s_n$ extends $s$ (denoted 
$s\!\subseteq\! s_n$).\smallskip

(b) $|s_n|\! =\! n$.\medskip

\noindent We put $\mathbb{G}_0\! :=\!\{ (s_n0\gamma ,s_n1\gamma )\mid n\!\in\!{\omega}\wedge
\gamma\!\in\! 2^{\omega}\}$. The following result is proved in [K-S-T].
 
\begin{thm} \label{KSTinjgr} (Kechris, Solecki, Todor\v cevi\'c) Let $X$ be a Polish space, and $A$ be an analytic graph on $X$. We assume that $A$ is acyclic or locally countable. Then exactly one of the following holds:\smallskip  

(a) there is $c\! :\! X\!\rightarrow\!\omega$ Borel such that 
$A\!\subseteq\! (c\!\times\! c)^{-1}\big(\neg\Delta (\omega )\big)$,\smallskip  

(b) there is $f\! :\! 2^\omega\!\rightarrow\! X$ injective continuous such that 
$s(\mathbb{G}_0)\!\subseteq\! (f\!\times\! f)^{-1}(A)$.\end{thm} 
 
 This seems to indicate that there is a hope to get $f$ and $g$ injective in Theorem 
\ref{motivating}.(b) for the first classes of the hierarchy if we assume acyclicity or local countability. This is the main purpose of this paper, and leads to the following notation. Let $X,Y,X',Y'$ be topological spaces and $A,B\!\subseteq\! X\!\times\! Y$, $A',B'\!\subseteq\! X'\!\times\! Y'$. We write\bigskip

\leftline{$(X,Y,A,B)\sqsubseteq (X',Y',A',B')\Leftrightarrow$}\smallskip

\rightline{$\exists f\! :\! X\!\rightarrow\! X'~~
\exists g\! :\! Y\!\rightarrow\! Y'\mbox{ injective continuous with }A\!\subseteq\! (f\!\times\! g)^{-1}(A')
\mbox{ and }B\!\subseteq\! (f\!\times\! g)^{-1}(B').$}\bigskip

\noindent We want to study the Borel and Wadge classes of the locally countable Borel relations: the Borel classes of rank one or two, the Lavrentieff classes built with the open sets (the classes of differences of open sets), their dual classes and their ambiguous classes. We will also study the Lavrentieff classes built with the $F_\sigma$ sets and their dual classes.

\begin{defi} Let $\eta\! <\!\omega_1$. If $(O_\theta )_{\theta <\eta}$ is an increasing sequence of subsets of a set $X$, then 
$$D\big( (O_\theta )_{\theta <\eta}\big)\! :=\!\big\{ x\!\in\! X\mid\exists\theta\! <\!\eta\ \ 
\mbox{parity}(\theta )\!\not=\!\mbox{parity}(\eta )\mbox{ and }
x\!\in\! O_\theta\!\setminus\!\big(\bigcup_{\theta'<\theta}\ O_{\theta'}\big)\big\}.$$
Now $D_\eta (\boraxi )(X)\! :=\!\big\{ D\big( (O_\theta )_{\theta <\eta}\big)\mid\forall\theta\! <\!\eta\ \  
O_\theta\!\in\!\boraxi (X)\big\}$, for each $1\!\leq\!\xi\! <\!\omega_1$. The classes 
$D_\eta (\boraxi )$, $\check D_\eta (\boraxi )$ and $\Delta\big( D_\eta (\boraxi )\big)$ form the 
\bf difference hierarchy\rm .\end{defi}

 Some recent work of the first author shows that having $f$ and $g$ injective in Theorem \ref{motivating}.(b) can be used to get results of reduction on the whole product, under some acyclicity and also topological assumptions. Some of the results in the present paper will be used by the first author in a future article on this topic. This work is also motivated by the work of Louveau on oriented graphs in [Lo4].

\begin{defi} Let $X$ be a set, and $A$ be a relation on $X$.\smallskip

(a) We say that $A$ is \bf antisymmetric\it\ if $A\cap A^{-1}\!\subseteq\!\Delta (X)$.\smallskip

(b) We say that $A$ is an \bf oriented graph\it\ if $A$ is irreflexive and antisymmetric.\end{defi}

 It follows from results of Wadge and Martin that inclusion well-orders 
$$\{ {\bf\Gamma}\cup\check {\bf\Gamma}\mid {\bf\Gamma}\mbox{ Wadge class of Borel sets}\}\mbox{,}$$ 
giving rise to an ordinal assignment $w({\bf\Gamma})$. If $G$ is an analytic oriented graph, then we can define $w(G)$ as the least $w({\bf\Gamma})$ such that $G$ is separable from $G^{-1}$ by a $\mbox{pot}({\bf\Gamma})$ set $C$. It is well defined by the separation theorem. Moreover, it is useless in the definition to distinguish between dual classes, for if $C$ separates $G$ from 
$G^{-1}$, then so does $\neg C^{-1}$, which is potentially in $\check {\bf\Gamma}$. The main property of this assignment is that $w(G)\!\leq\! w(H)$ if there is a Borel homomorphism from $G$ into $H$. Louveau also considers a rough approximation of $w(G)$, which is the least countable ordinal $\xi$ for which $G$ is separable from $G^{-1}$ by a $\mbox{pot}(\borxi )$ set. He proves the following.

\begin{thm} \label{Loog} (Louveau) Let $\xi\!\in\!\{ 1,2\}$. Then there is a concrete analytic oriented graph $\mathbb{G}_\xi$ on $2^\omega$ such that, for any Polish space $X$, and for any analytic oriented graph $G$ on $X$,  exactly one of the following holds:\smallskip  

(a) the set $G$ is separable from $G^{-1}$ by a $\mbox{pot}(\borxi )$ set,\smallskip

(b) there is $f\! :\! 2^\omega\!\rightarrow\! X$ continuous such that 
$\mathbb{G}_\xi\!\subseteq\! (f\!\times\! f)^{-1}(G)$.\end{thm} 

 Our main results are the following.\bigskip
 
\noindent $\bullet$ We generalize Theorem \ref{Loog} to all the $\borxi$'s, and all the Wadge classes of Borel sets. 

\begin{thm} \label{Zog} Let $\bf\Gamma$ be a Wadge class of Borel sets, or the class 
$\borxi$ for some countable ordinal $\xi\!\geq\! 1$. Then there is a concrete Borel oriented graph 
$\mathbb{G}_{\bf\Gamma}$ on $2^\omega$ such that, for any Polish space $X$, and for any analytic oriented graph $G$ on $X$,  exactly one of the following holds:\smallskip  

(a) the set $G$ is separable from $G^{-1}$ by a $\mbox{pot}({\bf\Gamma})$ set,\smallskip

(b) there is $f\! :\! 2^\omega\!\rightarrow\! X$ continuous such that 
$\mathbb{G}_{\bf\Gamma}\!\subseteq\! (f\!\times\! f)^{-1}(G)$.\end{thm} 

 We also investigate the injective version of this, for the first classes of the hierarchies again.\bigskip
 
\noindent $\bullet$ In the sequel, it will be very convenient to say that a relation $A$ on a set $X$ is \bf s-acyclic\rm\ if $s(A)$ is acyclic.
 
\begin{thm} \label{main} Let ${\bf\Gamma}\!\in\!\{ D_\eta (\boraone ),\check D_\eta (\boraone ),
D_n(\boratwo ),\check D_n(\boratwo )\mid 1\!\leq\!\eta\! <\!\omega_1,1\!\leq\! n\! <\!\omega\}\cup
\{\bortwo\}$. Then there are concrete disjoint Borel relations $\mathbb{S}_0$, $\mathbb{S}_1$ on $2^\omega$ such that, for any Polish space $X$, and for any disjoint analytic relations $A,B$ on $X$ with s-acyclic union, exactly one of the following holds:\smallskip

(a) the set $A$ is separable from $B$ by a $\mbox{pot}({\bf\Gamma})$ set,\smallskip

(b) $(2^\omega ,2^\omega ,\mathbb{S}_0,\mathbb{S}_1)\sqsubseteq (X,Y,A,B)$.\end{thm}

 In fact, we prove a number of extensions of this result. It also holds\bigskip
 
\noindent - for $\eta\! =\! 0$ if we replace $2^\omega$ with $1$,\bigskip
 
\noindent - with $f\! =\! g$ if 
${\bf\Gamma}\!\notin\!\{ D_\eta (\boraone ),\check D_\eta (\boraone )\mid\eta\! <\!\omega_1\}$; if 
${\bf\Gamma}\!\in\!\{ D_\eta (\boraone ),\check D_\eta (\boraone )\mid\eta\! <\!\omega_1\}$, then there is an antichain basis with two elements for the square reduction (it is rather unusual to have an antichain basis but no minimum object in this kind of dichotomy),\bigskip
 
\noindent - if we assume that $A\cup B$ is locally countable instead of s-acyclic when 
${\bf\Gamma}\!\subseteq\!\bormtwo$ (this also holds in rectangular products $X\!\times\! Y$),\bigskip

\noindent - if we only assume that $A$ is s-acyclic or locally countable when 
${\bf\Gamma}\! =\!\bormtwo$.\bigskip

 The situation is more complicated for the ambiguous classes.

\begin{thm} \label{mainD} Let ${\bf\Gamma}\!\in\!
\big\{\Delta\big( D_\eta (\boraone )\big)\mid1\!\leq\!\eta\! <\!\omega_1\big\}$. Then there is a  concrete finite antichain $\cal A$, made of tuples 
$(2^\omega,2^\omega ,\mathbb{S}_0,\mathbb{S}_1)$ where $\mathbb{S}_0$, $\mathbb{S}_1$ are disjoint Borel relations $\mathbb{S}_0$, $\mathbb{S}_1$ on $2^\omega$, such that, for any Polish space $X$, and for any disjoint analytic relations $A,B$ on $X$ whose union is contained in a potentially closed s-acyclic relation $R$, exactly one of the following holds:\smallskip

(a) the set $A$ is separable from $B$ by a $\mbox{pot}({\bf\Gamma})$ set,\smallskip

(b) there is $(2^\omega ,2^\omega ,\mathbb{A},\mathbb{B})\!\in\! {\cal A}$ with 
$(2^\omega ,2^\omega ,\mathbb{A},\mathbb{B})\sqsubseteq (X,Y,A,B)$.\end{thm}

 Here again, we can say more. This also holds\bigskip
 
\noindent - if we assume that $R$ is locally countable instead of s-acyclic (this also holds in rectangular products $X\!\times\! Y$),\bigskip
 
\noindent - in all those cases, $\cal A$ has size three if $\eta$ is a successor ordinal, and size one if $\eta$ is a limit ordinal (it is quite remarkable that the situation depends on the fact that $\eta$ is limit or not, it confirms the difference observed in the description of Wadge classes of Borel sets in terms of operations on sets present in [Lo1]),\bigskip
 
\noindent - with $f\! =\! g$, but in order to ensure this $\cal A$ must have size six if $\eta$ is a successor ordinal, and size two if $\eta$ is a limit ordinal.\bigskip
 
\noindent $\bullet$ We characterize when part (b) in the injective reduction property holds.
 
\begin{thm} Let ${\bf\Gamma}\!\in\!\{ D_\eta (\boraone ),\check D_\eta (\boraone ),
D_n(\boratwo ),\check D_n(\boratwo )\mid 1\!\leq\!\eta\! <\!\omega_1,1\!\leq\! n\! <\!\omega\}\cup
\{\bortwo\}$. Then there are concrete disjoint Borel relations $\mathbb{S}_0$, $\mathbb{S}_1$ on $2^\omega$ such that, for any Polish space $X$, and for any disjoint analytic relations $A,B$ on $X$, the following are equivalent:\smallskip

\noindent (1) there is an s-acyclic relation $R\!\in\!\ana$ such that $A\cap R$ is not separable from $B\cap R$ by a $\mbox{pot}({\bf\Gamma})$ set,\smallskip

\noindent (2) $(2^\omega ,2^\omega ,\mathbb{S}_0,\mathbb{S}_1)\sqsubseteq (X,Y,A,B)$.
\end{thm}

 The same kind of extensions as before hold (except that we cannot assume local countability instead of s-acyclicity for the classes of rank two).
 
\begin{thm} Let ${\bf\Gamma}\!\in\!
\big\{\Delta\big( D_\eta (\boraone )\big)\mid1\!\leq\!\eta\! <\!\omega_1\big\}$. Then there is a  concrete finite antichain $\cal A$, made of tuples 
$(2^\omega,2^\omega ,\mathbb{S}_0,\mathbb{S}_1)$ where $\mathbb{S}_0$, $\mathbb{S}_1$ are disjoint Borel relations $\mathbb{S}_0$, $\mathbb{S}_1$ on $2^\omega$, such that, for any Polish space $X$, and for any disjoint analytic relations $A,B$ on $X$, the following are equivalent:\smallskip

\noindent (1) there is a potentially closed s-acyclic relation $R\!\in\!\ana$ such that $A\cap R$ is not separable from $B\cap R$ by a $\mbox{pot}({\bf\Gamma})$ set,\smallskip

\noindent (2) there is $(2^\omega ,2^\omega ,\mathbb{A},\mathbb{B})\!\in\! {\cal A}$ with 
$(2^\omega ,2^\omega ,\mathbb{A},\mathbb{B})\sqsubseteq (X,Y,A,B)$.\end{thm} 

 Here again, the same kind of extensions as before hold.\bigskip
 
\noindent $\bullet$ The injective versions of Theorem \ref{Zog} mentioned earlier are as follows.
 
\begin{thm} \label{Zog2} Let ${\bf\Gamma}\!\in\!\{ D_\eta (\boraone ),\check D_\eta (\boraone ),
D_n(\boratwo ),\check D_n(\boratwo )\mid 1\!\leq\!\eta\! <\!\omega_1,1\!\leq\! n\! <\!\omega\}\cup
\{\bortwo\}$. Then there is a concrete Borel oriented graph $\mathbb{G}_{\bf\Gamma}$ on 
$2^\omega$ such that, for any Polish space $X$, and for any analytic s-acyclic oriented graph $G$ on $X$,  exactly one of the following holds:\smallskip  

(a) the set $G$ is separable from $G^{-1}$ by a $\mbox{pot}({\bf\Gamma})$ set,\smallskip

(b) there is $f\! :\! 2^\omega\!\rightarrow\! X$ injective continuous such that 
$\mathbb{G}_{\bf\Gamma}\!\subseteq\! (f\!\times\! f)^{-1}(G)$.\end{thm} 

 This result also holds if we assume that $G$ is locally countable instead of s-acyclic when 
${\bf\Gamma}\!\subseteq\!\bormtwo$.

\begin{thm} \label{Zog3} Let ${\bf\Gamma}\!\in\!
\big\{\Delta\big( D_\eta (\boraone )\big)\mid1\!\leq\!\eta\! <\!\omega_1\big\}$. Then there is a  concrete finite antichain $\cal A$, made of Borel oriented graphs on $2^\omega$, such that, for any Polish space $X$, and for any analytic oriented graph $G$ on $X$ contained in a potentially closed s-acyclic relation, exactly one of the following holds:\smallskip  

(a) the set $G$ is separable from $G^{-1}$ by a $\mbox{pot}({\bf\Gamma})$ set,\smallskip

(b) we can find $\mathbb{G}_{\bf\Gamma}\!\in\! {\cal A}$ and $f\! :\! 2^\omega\!\rightarrow\! X$ injective continuous such that $\mathbb{G}_{\bf\Gamma}\!\subseteq\! (f\!\times\! f)^{-1}(G)$.
\end{thm} 

 The same kind of extensions as before hold, except that $\cal A$ has size three if $\eta$ is a successor ordinal, and size two if $\eta$ is a limit ordinal.\bigskip
 
\noindent $\bullet$ At the end of the paper, we study the limits of our results and give negative results.

\section{$\!\!\!\!\!\!$ Generalities}

$\underline{\mbox{\bf The acyclic and the locally countable cases}}$\bigskip

 In [K-S-T], Section 6, the authors introduce the notion of an almost acyclic analytic graph, in order to prove an injective version of the $\mathbb{G}_0$-dichotomy for acyclic or locally countable analytic graphs. We now give a similar definition, in order to prove injective versions of Theorem \ref{motivating} for the first classes of the hierarchies. This definition is sufficient to cover all our cases, even if it is not always optimal.
 
\begin{defi} Let $X$ be a Polish space, and $A$ be a relation on $X$. We say that $A$ is \bf quasi-acyclic\it\ if there is a sequence $(C_n)_{n\in\omega}$ of 
$\mbox{pot}(\bormone )$ relations on $X$ with disjoint union $A$ such that, for any $s(A)$-path $(z_i)_{i\leq 2}$ with $z_0\!\not=\! z_2$, and for any $n_1,...,n_k\!\in\!\omega$, $C'_{n_i}\!\in\!\{ C^{}_{n_i},C^{-1}_{n_i}\}$ ($1\!\leq\! i\!\leq\! k$), $x_1,y_1,...,x_k,y_k$ in $X\!\setminus\!\{ z_i\mid i\!\leq\! 2\}$, if 
$(z_0,x_1),(z_2,y_1)\!\in\! C'_{n_1}$, $(x_1,x_2),(y_1,y_2)\!\in\! C'_{n_2}$, ..., $(x_{k-1},x_k),(y_{k-1},y_k)\!\in\! C'_{n_k}$ all hold, then $x_k\!\not=\! y_k$.\end{defi}

\begin{lem} \label{suffqa} Let $X$ be a Polish space, and $A$ be a Borel relation on $X$. We assume that $A$ is either s-acyclic and $\mbox{pot}(\boratwo )$, or locally countable. Then $A$ is quasi-acyclic.\end{lem}

\noindent\bf Proof.\rm ~Assume first that $A$ is s-acyclic and $\mbox{pot}(\boratwo )$. Then we can write $A\! =\!\bigcup_{n\in\omega}~C_n$, where $(C_n)_{n\in\omega}$ is a disjoint sequence of potentially closed relations on $X$. The acyclicity of $s(A)$ shows that $A$ is quasi-acyclic.\bigskip

 Assume now that $A$ is locally countable. By 18.10 in [K], $A$ can be written as $\bigcup_{q\in\omega}~G_q$, where $G_q$ is the Borel graph of a partial function $f_q$, and we may assume that the $G_q$'s are pairwise disjoint. By 18.12 in [K], the projections of the $G_q$'s are Borel. By Lemma 2.4.(a) in [L2], there is, for each $q$, a countable partition $(D^q_p)_{p\in\omega}$ of the domain of $f_q$ into Borel sets on which $f_q$ is injective. So the $C_n$'s are the 
$\mbox{Gr}({f_q}_{\vert D^q_p})$'s.\hfill{$\square$}\bigskip

\noindent $\underline{\mbox{\bf Topologies}}$\bigskip

 Let $Z$ be a recursively presented Polish space (see [M] for the basic notions of effective theory).\bigskip
 
\noindent (1) The topology ${\it\Delta}_Z$ on $Z$ is generated by $\Borel (Z)$. This topology is Polish (see (iii) $\Rightarrow$ (i) in the proof of Theorem 3.4 in [Lo3]). The topology $\tau_1$ on $Z^2$ is ${\it\Delta}_Z^2$. If $2\!\leq\!\xi\! <\!\omega_1^{\mbox{CK}}$, then the topology $\tau_\xi$ on $Z^2$ is generated by $\Ana\cap\bormlxi (\tau_1)$.\bigskip

\noindent (2) The \bf Gandy\!\! -\!\! Harrington\ topology\rm\ on $Z$ is generated by $\Ana (Z)$ and denoted ${\it\Sigma}_Z$. Recall the following facts about ${\it\Sigma}_Z$ (see [L7]).\smallskip

(a) ${\it\Sigma}_Z$ is finer than the initial topology of $Z$.\smallskip

(b) We set $\Omega_Z :=\{ z\!\in\! Z\mid\omega_1^z\! =\!\omega_1^{\mbox{CK}}\}$. Then 
$\Omega_Z$ is $\Ana (Z)$ and dense in $(Z,{\it\Sigma}_Z)$.\smallskip

(c) $W\cap\Omega_Z$ is a clopen subset of $(\Omega_Z ,{\it\Sigma}_Z )$ for each 
$W\!\in\!\Ana (Z)$.\smallskip

(d) $(\Omega_Z ,{\it\Sigma}_Z)$ is a zero-dimensional Polish space.

\section{$\!\!\!\!\!\!$ The classes $D_\eta (\boraone )$ and $\check D_\eta (\boraone )$}

$\underline{\mbox{\bf Examples}}$\bigskip

 In Theorem \ref{motivating}, either $\mathbb{S}_0$ or $\mathbb{S}_1$ is not locally countable if 
$\bf\Gamma$ is not self-dual. If ${\bf\Gamma}\!\subseteq\!\bortwo$, we can find disjoint analytic locally countable relations $A,B$ on $2^\omega$ such that $A$ is not separable from $B$ by a 
$\mbox{pot}({\bf\Gamma})$ set, as we will see. This shows that, in order to get partial reductions with injectivity, we have to use examples different from those in [L8], so that we prove the following.\bigskip

\noindent\bf Notation.\rm ~We introduce examples in the style of $\mathbb{G}_0$ in order to prove a dichotomy for the classes $D_\eta (\boraone )$, where $\eta\!\geq\! 1$ is a countable ordinal.\bigskip

\noindent $\bullet$ If $t\!\in\! 2^{<\omega}$, then 
$N_t\! :=\!\{\alpha\!\in\! 2^\omega\mid t\!\subseteq\!\alpha\}$ is the usual basic clopen set.\bigskip

\noindent $\bullet$ As in Section 2 in [L2] we inductively define 
$\varphi_\eta\! :\!\omega^{<\omega}\!\rightarrow\!\{ -1\}\cup (\eta\! +\! 1)$ by 
$\varphi_\eta (\emptyset )\! =\!\eta$ and 
$$\varphi_\eta (sn) = \left\{\!\!\!\!\!\!
\begin{array}{ll} 
& -1\mbox{ if }\varphi_\eta (s)\!\leq\! 0\mbox{,}\cr\cr 
& \theta\mbox{ if }\varphi_\eta (s)\! =\!\theta\! +\! 1\mbox{,}\cr\cr
& \mbox{ an odd ordinal such that the sequence }\big(\varphi_\eta (sn)\big)_{n\in\omega}
\mbox{ is cofinal in }\varphi_\eta (s)\cr 
& \mbox{and strictly increasing if }\varphi_\eta (s)\! >\! 0\mbox{ is limit.}
\end{array}\right.$$
If no confusion is possible, then we will write $\varphi$ instead of $\varphi_\eta$. We set 
$T_\eta\! :=\!\{ s\!\in\!\omega^{<\omega}\mid\varphi_\eta (s)\!\not=\! -1\}$, which is a wellfounded tree.\bigskip

\noindent $\bullet$ Let $(p_q)_{q\in\omega}$ be the sequence of prime numbers, and 
$<.>_\eta\ :\! T_\eta\!\rightarrow\!\omega$ be the following bijection. We define 
$I\! : \! T_\eta\!\rightarrow\!\omega$ by $I(\emptyset )\! :=\! 0$ and 
$I(s)\! :=\! p_0^{s(0)+1}...p_{\vert s\vert -1}^{s(\vert s\vert -1)+1}$ if $s\!\not=\!\emptyset$. As $I$ is injective, there is an increasing bijection $J\! :\! I[T_\eta ]\!\rightarrow\!\omega$. We set 
$<.>_\eta\ :=\! J\circ I$. Note that $<sq>_\eta\! -\! <s>_\eta\ \geq\! q\! +\! 1$ if $sq\!\in\! T_\eta$. Indeed, $I(s0)$, ..., $I\big( s(q\! -\! 1)\big)$ are strictly between $I(s)$ and $I(sq)$.\bigskip

\noindent $\bullet$ Let $\psi\! :\!\omega\!\rightarrow\! 2^{<\omega}$ be the map defined by 
$\emptyset ,\emptyset ,0,0,1,1,0^2,0^2,01,01,10,10,1^2,1^2,...$, so that $\vert\psi (q)\vert\!\leq\! q$ and $\psi [\{ 2n\mid n\!\in\!\omega\} ],\psi [\{ 2n\! +\! 1\mid n\!\in\!\omega\} ]\! =\! 2^{<\omega}$.\bigskip

\noindent $\bullet$ For each $s\in T_\eta$, we define $(t^0_s,t^1_s)\in (2\times 2)^{<\omega}$ by $t^\varepsilon_\emptyset=\emptyset$, and 
$t_{sq}^\varepsilon=t^\varepsilon_s\psi (q)0^{<sq>_\eta -<s>_\eta -|\psi (q)|-1}\varepsilon$. Note that this is well defined, $|t^\varepsilon_s|\! =<s>_\eta$ and $\mbox{Card}\big(\{ l<\; <s>_\eta\mid t^0_s(l)\!\not=\! t^1_s(l)\}\big)\! =\!\vert s\vert$ for each $s\!\in\! T_\eta$.\bigskip

\noindent $\bullet$ We set ${\cal T}^\eta\! :=\!
\big\{\big( t^0_sw,t^1_sw\big)\mid s\!\in\! T_\eta\wedge  w\!\in\! 2^{<\omega}\big\}$. The following properties are satisfied.\bigskip

\noindent - ${\cal T}^\eta$ is a tree on $2\!\times\! 2$, and 
$\lceil {\cal T}^\eta\rceil\!\subseteq\!\mathbb{E}_0\! :=\!
\{ (\alpha ,\beta )\!\in\! 2^\omega\!\times\! 2^\omega\mid
\exists m\!\in\!\omega ~~\forall n\! >\! m~~\alpha (n)\! =\!\beta (n)\}$ is locally countable.\bigskip

\noindent - If $(s,t)\!\in\! {\cal T}^\eta$ and $s(l)\!\not=\! t(l)$, then $s(l)\! <\! t(l)$.\bigskip

\noindent - For each $l\!\in\!\omega$, there is exactly one sequence 
$(u,v)\!\in\! {\cal T}^\eta\cap (2^{l+1}\!\times\! 2^{l+1})$ such that $u(l)\!\not=\! v(l)$ since 
$t^0_{sq}(<sq>_\eta -\! 1)\!\not=\! t^1_{sq}(<sq>_\eta -\! 1)$ (in fact, $(u,v)$ is of the form 
$(t^0_s,t^1_s)$ for some $s$). In particular, 
$s\big( {\cal T}^\eta\cap (2^{l+1}\!\times\! 2^{l+1})\big)\!\setminus\!\Delta (2^{l+1})$ is a connected acyclic graph on $2^{l+1}$, inductively.\bigskip

\noindent $\bullet$ We set, for $\varepsilon\!\in\! 2$, 
$$\mathbb{N}^\eta_\varepsilon\! :=\!\Big\{ (t^0_s\gamma ,t^1_s\gamma)\mid 
s\!\in\! T_\eta\ \wedge\ \mbox{parity}(\vert s\vert )\! =\!\varepsilon\ \wedge\ 
\gamma\!\in\! 2^\omega\Big\} .$$
If $s\!\in\! T_\eta$, then $f_s\! :\! N_{t^0_s}\!\rightarrow\! N_{t^1_s}$ is the partial homeomophism with clopen domain and range defined by $f_s(t^0_s\gamma )\! :=\! t^1_s\gamma$, so that 
$\mathbb{N}^\eta_\varepsilon\! =\!\bigcup_{s\in T_\eta ,\mbox{parity}(\vert s\vert )=\varepsilon}~\mbox{Gr}(f_s)$. We set $C_s\! :=\!\bigcup_{q\in\omega}~\mbox{Gr}(f_{sq})$ when it makes sense (i.e., when $\varphi_\eta (s)\!\geq\! 1$). For $\eta\! =\! 0$, we set $\mathbb{N}^\eta_0\! :=\! 1^2$ and 
$\mathbb{N}^\eta_1\! :=\!\emptyset$ (in $1^2$).

\begin{lem} \label{S^2} Let $\eta$ be a countable ordinal, and $C$ be a nonempty clopen subset of $2^\omega$.\smallskip

(a) If $\varphi_\eta (s)\!\geq\! 1$ and $G$ is a dense $G_\delta$ subset of $2^\omega$, then 
$\overline{C_s}\cap (C\cap G)^2\!\subseteq\!\overline{C_s\cap (C\cap G)^2}$.\smallskip

(b) $\mathbb{N}^\eta_0\cap C^2$ is not separable from $\mathbb{N}^\eta_1\cap C^2$ by a 
$\mbox{pot}\big( D_\eta (\boraone )\big)$ set.\end{lem}

\noindent\bf Proof.\rm ~(a) It is enough to prove that if $q\!\in\!\omega$, then 
$\mbox{Gr}(f_{sq})\cap C^2\!\subseteq\!\overline{\mbox{Gr}(f_{sq})\cap (C\cap G)^2}$. This comes from the proof of Lemma 3.5 in [L1], but we recall it for self-containedness. Let $U,V$ be open subsets of $C$ such that $\mbox{Gr}(f_{sq})\cap (U\!\times\! V)\!\not=\!\emptyset$. Then 
$N_{t^1_{sq}}\cap V\cap G$ is a dense $G_\delta$ subset of $N_{t^1_{sq}}\cap V$, so that 
$f_{sq}^{-1}(V\cap G)$ is a dense $G_\delta$ subset of $f_{sq}^{-1}(V)$. Thus 
$G\cap f_{sq}^{-1}(V)$ and $G\cap f_{sq}^{-1}(V\cap G)$ are dense $G_\delta$ subsets of 
$f_{sq}^{-1}(V)$. This gives $\alpha$ in this last set and $U\cap f_{sq}^{-1}(V)$. Therefore 
$\big(\alpha ,f_{sq}(\alpha )\big)$ is in $\mbox{Gr}(f_{sq})\cap (C\cap G)^2\cap (U\!\times\! V)$.\bigskip

\noindent (b) We may assume that $\eta\!\geq\! 1$. We argue by contradiction, which gives $P\!\in\!\mbox{pot}\big( D_\eta (\boraone )\big)$, and a dense $G_\delta$ subset of 
$2^\omega$ such that $P\cap G^2\!\in\! D_\eta (\boraone )(G^2)$. So let 
$(O_\theta )_{\theta <\eta}$ be a sequence of open relations on $2^\omega$ such that 
$P\cap G^2\! =\!\big(\bigcup_{\theta <\eta ,\mbox{parity}(\theta )\not=\mbox{parity}(\eta )}~
O_\theta\!\setminus\! (\bigcup_{\theta'<\theta}~O_{\theta'})\big)\cap G^2$.\bigskip

\noindent $\bullet$ Let us show that if $\theta\!\leq\!\eta$, $s\!\in\! T_\eta$ and 
$\varphi (s)\! =\!\theta$, then $\mbox{Gr}(f_s)\cap (C\cap G)^2\!\subseteq\!\neg O_\theta$ if 
$\theta\! <\!\eta$, and $\mbox{Gr}(f_s)\cap (C\cap G)^2$ is disjoint from 
$\bigcup_{\theta'<\theta}~O_{\theta'}$ if $\theta\! =\!\eta$. The objects $s\! =\!\emptyset$ and 
$\theta\! =\!\eta$ will give the contradiction.\bigskip

\noindent $\bullet$ We argue by induction on $\theta$. Note that if $s\!\in\! T_\eta$, $\vert s\vert$ is even if and only if $\varphi (s)$ has the same parity as $\eta$. If $\theta =0$, then $\vert s\vert$ has the same parity as $\eta$, thus $\mbox{Gr}(f_s)\cap (C\cap G)^2\!\subseteq\!
\mathbb{N}^\eta_{\mbox{parity}(\eta )}\cap G^2\!\subseteq\!\neg O_0$.\bigskip

\noindent $\bullet$ Assume that the result has been proved for $\theta'\! <\!\theta$. If $\theta$ is the successor of $\theta'$, then the induction assumption implies that 
$\mbox{Gr}(f_{sq})\cap (C\cap G)^2\!\subseteq\!\neg O_{\theta'}$ for each $q$. So 
$C_s\cap (C\cap G)^2\!\subseteq\!\neg O_{\theta'}$ and 
$\overline{C_s\cap (C\cap G)^2}\!\subseteq\!\neg O_{\theta'}$. By (a), we get 
$\overline{C_s}\cap (C\cap G)^2\!\subseteq\!\overline{C_s\cap (C\cap G)^2}$, which gives the desired inclusion if $\theta\! =\!\eta$ since $\mbox{Gr}(f_s)\!\subseteq\!\overline{C_s}$.\bigskip

 If $\theta\! <\!\eta$ and $\vert s\vert $ is even, then $\varphi (s)$ has the same parity as $\eta$ and the parity of $\theta'$ is opposite to that of $\eta$. Note that 
$\mbox{Gr}(f_s)\cap (C\cap G)^2\!\subseteq\!\mathbb{N}^\eta_0\cap G^2\!\subseteq\!\bigcup_{\theta''<\eta ,\mbox{parity}(\theta'')\not=\mbox{parity}(\eta )}~O_{\theta''}\!\setminus\! (\bigcup_{\theta'''<\theta''}~O_{\theta'''})\!\subseteq\!\neg O_\theta$.\bigskip

 If $\vert s\vert $ is odd, then the parity of $\varphi (s)$ is opposite to that of $\eta$ and $\theta'$ has the same parity as $\eta$. But if $s\!\in\! T_\eta$ has odd length, then 
$$\mbox{Gr}(f_s)\cap (C\cap G)^2\!\subseteq\!\mathbb{N}^\eta_1\cap G^2\!\subseteq\! 
G^2\!\setminus\! (\bigcup_{\theta''<\eta}~O_{\theta''})\cup
\bigcup_{\theta''<\eta ,\mbox{parity}(\theta'')=\mbox{parity}(\eta )}~O_{\theta''}\!\setminus\! 
(\bigcup_{\theta'''<\theta''}~O_{\theta'''}).$$ 
This gives the result.\bigskip

\noindent $\bullet$ If $\theta$ is limit, then $\big(\varphi (sn)\big)_{n\in\omega}$ is cofinal in 
$\varphi (s)$, and $\mbox{Gr}(f_{sn})\cap (C\cap G)^2\!\subseteq\!\neg O_{\varphi (sn)}$ by the induction assumption. If $\theta_0\! <\!\varphi (s)$, then there is $n(\theta_0)$ such that 
$\varphi (sn)\! >\!\theta_0$ if $n\!\geq\! n(\theta_0)$. Thus 
$\mbox{Gr}(f_{sn})\cap (C\cap G)^2\!\subseteq\!\neg O_{\theta_0}$ as soon as 
$n\!\geq\! n(\theta_0)$. But   
$$\mbox{Gr}(f_s)\cap (C\cap G)^2\!\subseteq\! (C\cap G)^2\cap\overline{C_s}\!\setminus\! C_s\! =\!
\overline{C_s\cap (C\cap G)^2}\!\setminus\! C_s\!\subseteq\!
\overline{\bigcup_{n\geq n(\theta_0)}~\mbox{Gr}(f_{sn})\cap (C\cap G)^2}\!\subseteq\!\neg 
O_{\theta_0}.$$ 
Thus $\mbox{Gr}(f_s)\cap (C\cap G)^2\!\subseteq\!\neg (\bigcup_{\theta'<\theta}~O_{\theta'})$.\bigskip

 If $\theta<\eta$, as $\vert s\vert$ has the same parity as $\eta$, we get 
$\mbox{Gr}(f_s)\cap (C\cap G)^2\!\subseteq\!\mathbb{N}^\eta_{\mbox{parity}(\eta )}\cap G^2$, so that $\mbox{Gr}(f_s)\cap (C\cap G)^2\!\subseteq\!\neg O_\theta$.\hfill{$\square$}\bigskip

\noindent $\underline{\mbox{\bf A topological characterization}}$\bigskip

\noindent\bf Notation.\rm\ Let $1\!\leq\!\xi\! <\!\omega_1^{\mbox{CK}}$. Theorem 4.1 in [L6] shows that if $A_0,A_1$ are disjoint $\Ana$ relations on 
$\omega^\omega$, then $A_0$ is separable from $A_1$ by a $\mbox{pot}(\boraxi )$ set exactly when $A_0\cap\overline{A_1}^{\tau_\xi}\! =\!\emptyset$. We now define the versions of 
$A_0\cap\overline{A_1}^{\tau_\xi}$ for the classes $D_\eta (\boraxi )$. So let $\varepsilon\!\in\! 2$ and $\eta\! <\!\omega_1^{\mbox{CK}}$. We define 
$\bigcap_{\theta <0}\ F^\varepsilon_{\theta ,\xi}\! :=\! (\omega^\omega )^2$, and, inductively,  
$$F^\varepsilon_{\eta ,\xi}\! :=\!\overline{A_{\vert\mbox{parity}(\eta)-\varepsilon\vert}\cap
\bigcap_{\theta <\eta}\ F^\varepsilon_{\theta ,\xi}}^{\tau_\xi} .$$ 
We will sometimes denote by $F^\varepsilon_{\eta ,\xi}(A_0,A_1)$ the sets $F^\varepsilon_{\eta ,\xi}$ previously defined. By induction, we can check that 
$F^\varepsilon_{\eta ,\xi}(A_1,A_0)\! =\! F^{1-\varepsilon}_{\eta ,\xi}(A_0,A_1)$.\bigskip

  Fix a bijection $l\!\mapsto\!\big( (l)_0,(l)_1\big)$ from $\omega$ onto $\omega^2$, with inverse map 
$(m,p)\!\mapsto <m,p>$. We define, for $u\!\in\!\omega^{\leq\omega}$ and $n\!\in\!\omega$, 
$(u)_n\!\in\!\omega^{\leq\omega}$ by $(u)_n(p)\! :=\! u(<n,p>)$ if $<n,p><\!\vert u\vert$.
 
\begin{thm} \label{kernel1} Let $1\!\leq\!\xi\! <\!\omega_1^{\mbox{CK}}$, $\eta\! =\!\lambda\! +\! 2k\! +\!\varepsilon\! <\!\omega_1^{\mbox{CK}}$ with $\lambda$ limit, $k\!\in\!\omega$ and $\varepsilon\!\in\! 2$, and $A_0$, $A_1$ be disjoint $\Ana$ relations on $\omega^\omega$. Then the following are equivalent:\smallskip

\noindent (1) the set $A_0$ is not separable from $A_1$ by a 
$\mbox{pot}\big( D_\eta (\boraxi )\big)$ set,\smallskip

\noindent (2) the $\Ana$ set $F^\varepsilon_{\eta ,\xi}$ is not empty.\end{thm}
 
\noindent\bf Proof.\rm ~This result is essentially proved in [L8]. However, the formula for 
$F^\varepsilon_{\eta ,\xi}$ is more concrete here, since the more general and abstract case of Wadge classes is considered in [L8]. So we give some details.\bigskip

\noindent $\bullet$ In [Lo-SR], the following class of sets is introduced. Let 
$1\!\leq\!\xi\! <\!\omega_1$ and $\bf\Gamma$, ${\bf\Gamma}'$ be two classes of sets. Then 
$A\!\in\! S_\xi ({\bf\Gamma},{\bf\Gamma}')\ \Leftrightarrow\ 
A\! =\!\bigcup_{p\geq 1}\ (A_p\cap C_p)\cup\left( B\!\setminus\!\bigcup_{p\geq 1}\ C_p\right)$,  
where $A_p\!\in\! {\bf\Gamma}$, $B\!\in\! {\bf\Gamma}'$, and $(C_p)_{p\geq 1}$ is a sequence of pairwise disjoint $\boraxi$ sets. The authors prove the following:
$$\boraxi\! =\! S_\xi (\check\{\emptyset\} ,\{\emptyset\} )\mbox{,}$$
$$D_{\theta +1}(\boraxi )\! =\! S_\xi (\check D_\theta (\boraxi ),\boraxi )
\mbox{ if }\theta\! <\!\omega_1\mbox{,}$$
$$D_\lambda (\boraxi )\! =\! S_\xi (\bigcup_{p\geq 1}~D_{\theta_p}(\boraxi ),\{\emptyset\} )
\mbox{ if }\lambda\! =\!\mbox{sup}_{p\geq 1}~\theta_p\mbox{ is limit.}$$
They also code the non self-dual Wadge classes of Borel sets by elements of 
$\omega_1^\omega$ as follows (we sometimes identify $\omega_1^\omega$ with 
$(\omega_1^\omega )^\omega$). The relations ``$u$ \bf is a second type description\rm " and 
``$u$ \bf describes\rm\ $\bf\Gamma$" (written $u\!\in\! {\cal D}$ and 
${\bf\Gamma}_u\! =\! {\bf\Gamma}$~-~ambiguously) are the least relations satisfying the following properties.\medskip

(a) If $u\! =\! 0^\infty$, then $u\!\in\! {\cal D}$ and ${\bf\Gamma}_u\! =\!\{\emptyset\}$.\smallskip

(b) If $u\! =\!\xi^\frown 1^\frown v$, with $v\!\in\! {\cal D}$ and $v(0)\! =\!\xi$, then 
$u\!\in\! {\cal D}$ and ${\bf\Gamma}_u\! =\!\check {\bf\Gamma}_v$.\smallskip

(c) If $u\! =\!\xi^\frown 2^\frown\!\!\! <\! u_p\! >$ satisfies $\xi\!\geq\! 1$, $u_p\!\in\! {\cal D}$, and 
$u_p(0)\!\geq\!\xi$ or $u_p(0)\! =\! 0$, then $u\!\in\! {\cal D}$ and 
${\bf\Gamma}_u\! =\! S_\xi(\bigcup_{p\geq 1}\ {\bf\Gamma}_{u_p},{\bf\Gamma}_{u_0})$.\medskip

 They prove that $\bf\Gamma$ is a non self-dual Wadge class of Borel sets exactly when there is $u\!\in\! {\cal D}$ such that 
${\bf\Gamma}(\omega^\omega )\! =\! {\bf\Gamma}_u(\omega^\omega )$.

\vfill\eject

\noindent $\bullet$ In [L8], the elements of $\cal D$ are coded by elements of $\omega^\omega$. An inductive operator $\mathfrak{H}$ over $\omega^\omega$ is defined and there is a partial function $c\! :\!\omega^\omega\!\rightarrow\!\omega_1^\omega$ with 
$c[\mathfrak{H}^\infty ]\! =\! {\cal D}$ (see  Lemma 6.2 in [L8]). Another operator $\mathfrak{J}$ on 
$(\omega^\omega )^3$ is defined in [L8] to code the non self-dual Wadge classes of Borel sets and their elements (see Lemma 6.5 in [L8]). We will need a last inductive operator $\mathfrak{K}$, on $(\omega^\omega )^6$, to code the sets that will play the role of the $\Ana$ sets 
$F^\varepsilon_{\eta ,\xi}$'s, via a universal set $\cal U$ for the class $\Ca (\omega^\omega\!\times\!\omega^\omega )$. More precisely, if 
$(\alpha ,a_0,a_1,b_0,b_1,r)\!\in\! \mathfrak{K}^\infty$, then $b_0,b_1$ and $r$ are completely determined by $(\alpha ,a_0,a_1)$ and in practice $\alpha$ will be in 
$\mathfrak{H}^\infty$, so that we will write $r\! =\! r(\alpha ,a_0,a_1)\! =\! r(u,a_0,a_1)$ if 
$u\! =\! c(\alpha )$. Our $\Ana$ sets $A_0,A_1$ are coded by $a_0,a_1$, in the sense that 
$A_\varepsilon\! =\!\neg {\cal U}_{a_\varepsilon}$. By Lemma 6.6 in [L8], there is a recursive map ${\cal A}\! :\! (\omega^\omega )^2\!\rightarrow\!\omega^\omega$ such that 
$\neg {\cal U}_{{\cal A}(\alpha ,r)}\! =\! (\neg {\cal U}_{(r)_0})\cap
\bigcap_{p\geq 1}~\overline{\neg {\cal U}_{(r)_p}}^{\tau_{\vert\alpha\vert}}$ if $\alpha\!\in\!\Borel$ codes a wellordering, where $r\!\mapsto\!\big( (r)_p\big)_{p\in\omega}$ is a bijection from 
$\omega^\omega$ onto $(\omega^\omega )^\omega$. In the sequel, all the closures will be for 
$\tau_\xi$.\bigskip

\noindent $\bullet$ We argue by induction on $\eta$. As $D_0(\boraxi )\! =\!\{\emptyset\}$, $A_0$ is separable from $A_1$ by a $D_0(\boraxi )$ set when $A_0\! =\!\emptyset$, which is equivalent to $F^0_{0,\xi}\! =\!\overline{A_0}\! =\!\emptyset$. As $D_1(\boraxi )\! =\!\boraxi$, $A_0$ is separable from $A_1$ by a $D_1(\boraxi )$ set when $A_0\cap\overline{A_1}\! =\!\emptyset$ by Theorem 4.1 in [L6], which is equivalent to $F^1_{1,\xi}\! =\!\overline{A_0\cap\overline{A_1}}\! =\!\emptyset$.\bigskip

 Let us do these two basic cases in the spirit of the material from [L8] previously described, which will be done also for the other more complex cases.\bigskip
 
\noindent - Note that $D_0(\boraxi )\! =\!\{\emptyset\}\! =\! {\bf\Gamma}_{0^\infty}$. Let 
$\alpha\!\in\!\Borel$ such that $(\alpha )_n$ codes a wellordering of order type $0$ for each 
$n\!\in\!\omega$. A look at the definition of $\mathfrak{H}$ shows that 
$\alpha\!\in\!\mathfrak{H}^\infty$. Another look at Definition 6.3 in [L8] shows that $\alpha$ is normalized (this will never be a problem in the sequel as well). Lemma 6.5 in [L8] gives 
$\beta ,\gamma\!\in\!\omega^\omega$ with 
$(\alpha ,\beta ,\gamma )\!\in\!\mathfrak{J}^\infty$. Lemma 6.7 in [L8] gives 
$b_0,b_1,r\!\in\!\omega^\omega$ with $(\alpha ,a_1,a_0,b_0,b_1,r)\!\in\!\mathfrak{K}^\infty$. By Theorem 6.10 in [L8], $A_1$ is separable from $A_0$ by a 
$\mbox{pot}\big(\check D_0(\boraxi )\big)$ set if and only if $\neg {\cal U}_r\! =\!\emptyset$. A look at the definition of $\mathfrak{K}$ shows that $r\! =\! a_0$, so that $\neg {\cal U}_r\! =\! A_0$.\bigskip

\noindent - Now $D_1(\boraxi )\! =\!\boraxi\! =\! S_\xi (\check\{\emptyset\} ,\{\emptyset\} )\! =\! 
S_\xi ({\bf\Gamma}_{010^\infty},{\bf\Gamma}_{0^\infty})\! =\! 
S_\xi (\bigcup_{p\geq 1}~{\bf\Gamma}_{010^\infty},{\bf\Gamma}_{0^\infty})\! =\! 
{\bf\Gamma}_{v_1}$, where $v_1\! :=\!\xi 2<0^\infty ,010^\infty ,010^\infty ,...>$. As above, $A_1$ is separable from $A_0$ by a $\mbox{pot}\big(\check D_1(\boraxi )\big)$ set if and only if 
$\neg {\cal U}_r\! =\!\emptyset$. A look at the definition of $\mathfrak{K}$ shows that 
$r\! =\! b_0\! =\! {\cal A}(\alpha_1,<a_0,a_1,a_1,...>)$, where $\vert\alpha_1\vert\! =\!\xi$. Thus 
$\neg {\cal U}_r\! =\! A_0\cap\overline{A_1}$.\bigskip

 In the general case, there is $v_\eta\!\in\! {\cal D}$ such that 
$D_\eta (\boraxi )\! =\! {\bf\Gamma}_{v_\eta}$ and $A_1$ is separable from $A_0$ by a 
$\mbox{pot}\big(\check D_\eta (\boraxi )\big)$ set if and only if 
$\neg {\cal U}_{r(v_\eta ,a_1,a_0)}\! =\!\emptyset$. Moreover,\medskip

(a) if $v_\eta\! =\! 0^\infty$, then $r(v_\eta ,a_1,a_0)\! =\! a_0$,\smallskip

(b) if $v_\eta\! =\!\xi^\frown 1^\frown v$, then $r(v_\eta ,a_1,a_0)\! =\! a_1$,\smallskip

(c) if $v_\eta\! =\!\xi^\frown 2^\frown\!\!\! <\! u_p\! >$ and $r_p\! =\! r(u_p,a_1,a_0)$, then 
$r(v_\eta ,a_1,a_0)\! =\! r(u_0,b_1,b_0)$, where by definition 
$b_i\! :=\! {\cal A}(\alpha_1,<a_i,r_1,r_2,...>)$.\medskip

 It is enough to prove that 
$F^\varepsilon_{\eta ,\xi}\! =\!\overline{\neg {\cal U}_{r(v_\eta ,a_1,a_0)}}$, and we may assume that $\eta\!\geq\! 2$ by the previous discussion.\bigskip
 
\noindent $\bullet$ If $\eta$ is a limit ordinal, then fix a sequence $(\eta_p)_{p\in\omega}$ of even ordinals cofinal in $\eta$. Note that 
$$D_\eta (\boraxi )\! =\! S_\xi (\bigcup_{p\geq 1}~D_{\eta_p}(\boraxi ),\{\emptyset\} )\! =\! 
S_\xi (\bigcup_{p\geq 1}~{\bf\Gamma}_{u_p},{\bf\Gamma}_{u_0})\! =\! 
{\bf\Gamma}_{v_\eta}\mbox{,}$$ 
where $v_\eta\! =\!\xi^\frown 2^\frown <u_p>$.

\vfill\eject

 Therefore, if $r_p\! :=\! r(u_p,a_1,a_0)$, then 
$F_{\theta_p,\xi}^\varepsilon\! =\!\overline{\neg {\cal U}_{r_p}}$ if $p\!\geq\! 1$, by the induction hypothesis. On the other hand, $r(u_0,b_1,b_0)=b_0$. But 
$b_0\! =\!\mathcal{A}(\alpha_1,<a_0,r_1,r_2,...>)$, so that 
$$\neg\mathcal{U}_{b_0}\! =\! 
(\neg\mathcal{U}_{a_0})\cap\bigcap_{p\geq1}~\overline{\neg\mathcal{U}_{r_p}}\mbox{,}$$ 
as required.\bigskip

\noindent $\bullet$ If $\eta\! =\!\theta\! +\! 1$, then 
$$D_\eta (\boraxi )\! =\! S_\xi (\check D_\theta (\boraxi ),\boraxi )\! =\! 
S_\xi (\bigcup_{p\geq 1}~{\bf\Gamma}_{u_p},{\bf\Gamma}_{u_0})\! =\! 
{\bf\Gamma}_{v_\eta}\mbox{,}$$ 
where $v_\eta\! =\!\xi^\frown 2^\frown <u_p>$. Therefore, if $r_p\! :=\! r(u_p,a_1,a_0)$, then 
$F^\varepsilon_{\theta ,\xi}=\overline{\neg\mathcal U_{r_p}}$ if $p\!\geq\! 1$, by the induction hypothesis (there is a double inversion of the superscript, one because the parity of $\theta$ is different from that of $\eta$, and the other one because there is a complement, so that the roles of $A_0,A_1$ are exchanged). By the case $\eta\! =\! 1$ applied to $b_0$ and $b_1$, 
$\neg {\cal U}_{r(u_0,b_1,b_0)}\! =\!\neg\mathcal{U}_{b_0}\cap\overline{\neg\mathcal U_{b_1}}$. Note that
$$\neg\mathcal{U}_{b_i}\! =\! (\neg\mathcal U_{a_i})\cap\bigcap_{p\geq 1}~
\overline{\neg\mathcal{U}_{r_p}}\! =\! (\neg\mathcal U_{a_i})\cap F^\varepsilon_{\theta ,\xi}$$
since $b_i\! =\!\mathcal{A}(\alpha_1,<a_i,r_1,r_2,\ldots >)$. If $r\! :=\! r(v_\eta ,a_1,a_0)$, then 
$$\neg U_r\! =\! (\neg\mathcal U_{a_0})\cap F^\varepsilon_{\theta ,\xi}\cap
\overline{\neg\mathcal U_{a_1}\cap F^{\varepsilon}_{\theta ,\xi}}\! =\! 
A_0\cap F^\varepsilon_{\theta ,\xi}\mbox{,}$$
because $F^\varepsilon_{\theta ,\xi}\! =\!\overline{A_1\cap
\bigcap_{\rho <\theta}~F^\varepsilon_{\rho ,\xi}}
\!\subseteq\!\overline{A_1\cap\overline{A_1\cap\bigcap_{\rho <\theta}~F^\varepsilon_{\rho ,\xi}}}
\!\subseteq\!\overline{A_1\cap F^\varepsilon_{\theta ,\xi}}$ (since the parity of $\theta$ is different from $\varepsilon$). Finally, 
$\overline{\neg U_r}\! =\!\overline{A_0\cap F^\varepsilon_{\theta ,\xi}}\! =\! F^\varepsilon_{\eta ,\xi}$, as required.\hfill{$\square$}\bigskip

\noindent $\underline{\mbox{\bf The main result}}$\bigskip\indent

 We set, for $\eta\! <\!\omega_1$ and $\varepsilon\!\in\! 2$, $\mathbb{B}^\eta_\varepsilon\! :=\!
\{ (0\alpha ,1\beta )\mid (\alpha ,\beta )\!\in\!\mathbb{N}^\eta_\varepsilon\}$.

\begin{thm} \label{checkD2Sigma01} Let $\eta\!\geq\! 1$ be a countable ordinal, $X$ be a Polish space, and $A_0,A_1$ be disjoint analytic relations on $X$ such that $A_0\cup A_1$ is quasi-acyclic. The following are equivalent:\smallskip  

\noindent (1) the set $A_0$ is not separable from $A_1$ by a $\mbox{pot}\big( D_\eta (\boraone )\big)$ set,\smallskip  

\noindent (2) there is $(\mathbb{A}_0,\mathbb{A}_1)\!\in\!
\{ (\mathbb{N}^\eta_0,\mathbb{N}^\eta_1),(\mathbb{B}^\eta_0,\mathbb{B}^\eta_1)\}$ such that 
$(2^\omega ,2^\omega ,\mathbb{A}_0,\mathbb{A}_1)\sqsubseteq (X,X,A_0,A_1)$, via a square map,\smallskip

\noindent (3) $(2^\omega ,2^\omega ,\mathbb{N}^\eta_0,\mathbb{N}^\eta_1)\sqsubseteq (X,X,A_0,A_1)$.\end{thm}

\noindent\bf Proof.\rm ~(1) $\Rightarrow$ (2) Let $\varepsilon\! :=\!\mbox{parity}(\eta )$, and 
$(C_p)_{p\in\omega}$ be a witness for the quasi-acyclicity of $A_0\cup A_1$. We may assume that $X\! =\!\omega^\omega$. Indeed, we may assume that $X$ is zero-dimensional, and thus a closed subset of $\omega^\omega$. As $A_0$ is not separable from $A_1$ by a 
$\mbox{pot}\big( D_\eta (\boraone )\big)$ set in $X^2$, it is also the case in 
$(\omega^\omega )^2$, which gives $f\! :\! 2^\omega\!\rightarrow\!\omega^\omega$. As 
$\Delta (2^\omega )\!\subseteq\!\mathbb{N}^\eta_0$ and 
$\{ (0\alpha ,1\alpha )\mid\alpha\!\in\! 2^\omega\}\!\subseteq\!\mathbb{B}^\eta_0$, the range of 
$\Delta (2^\omega )$ by $f\!\times\! f$ is a subset of $X^2$, so that $f$ takes values in $X$. We may also assume that $A_0,A_1$ are $\Ana$, and that the relation ``$(x,y)\!\in\! C_p$" is $\Borel$ in $(x,y,p)$. By Theorem \ref{kernel1}, 
$$F^\varepsilon_\eta\! =\!\overline{A_0\cap\bigcap_{\theta <\eta}~F^\varepsilon_\theta}^{\tau_1}$$ 
is a nonempty $\Ana$ relation on $X$ (where $F^\varepsilon_\eta\! :=\! F^\varepsilon_{\eta ,1}$, for simplicity).

\vfill\eject

 We set, for $\theta\!\leq\!\eta$, $F_\theta\! :=\! A_{\vert\mbox{parity}(\theta )-\varepsilon\vert}\cap
\bigcap_{\theta'<\theta}~F^\varepsilon_{\theta'}$, so that 
$F^\varepsilon_\theta\! =\!\overline{F_\theta}^{\tau_1}$. We put, for $\theta\leq\eta$, 
$$D_\theta\! :=\!\big\{ (t^0_sw,t^1_sw)\!\in\! {\cal T}^\eta\mid\varphi (s)\! =\!\theta\big\}\mbox{,}$$ 
so that $(D_\theta )_{\theta\leq\eta}$ is a partition of ${\cal T}^\eta$. As 
$D_\eta\! =\!\Delta (2^{<\omega})$, 
$G_{l+1}\! :=\! s\big( (\bigcup_{\theta <\eta}~D_\theta )\cap (2^{l+1}\!\times\! 2^{l+1})\big)$ is a connected acyclic graph on $2^{l+1}$ for each $l\!\in\!\omega$.\bigskip

\noindent\bf Case 1\rm\ $F_\eta\!\not\subseteq\!\Delta (X)$.\bigskip

Let $(x,y)\!\in\! F_\eta\!\setminus\!\Delta (X)$, and $O_0,O_1$ be disjoint $\Borone$ sets with 
$(x,y)\!\in\! O_0\!\times\! O_1$. We can replace $F_\eta$, $A_0$ and $A_1$ with their intersection with $O_0\!\times\! O_1$ if necessary and assume that they are contained in $O_0\!\times\! O_1$.\bigskip

\noindent $\bullet$ We construct the following objects:\bigskip

- sequences $(x_s)_{s\in 2^{<\omega}}$, $(y_s)_{s\in 2^{<\omega}}$ of points of $X$,\smallskip

- sequences $(X_s)_{s\in 2^{<\omega}}$, $(Y_s)_{s\in 2^{<\omega}}$ of $\Ana$ subsets of $X$,\smallskip

- a sequence $(U_{s,t})_{(s,t)\in {\cal T}^\eta}$ of $\Ana$ subsets of $X^2$, and 
$\Phi\! :\! {\cal T}^\eta\!\rightarrow\!\omega$.\bigskip

\noindent We want these objects to satisfy the following conditions:
$$\begin{array}{ll}
& (1)\ x_s\!\in\! X_s\ \wedge\ y_s\!\in\! Y_s\ \wedge\ (x_s,y_t)\!\in\! U_{s,t}\cr
& (2)\ X_{s\varepsilon}\!\subseteq\! X_s\!\subseteq\!\Omega_X\cap O_0\ \wedge\ 
Y_{s\varepsilon}\!\subseteq\! Y_s\!\subseteq\!\Omega_X\cap O_1\ \wedge\ 
U_{s,t}\!\subseteq\! C_{\Phi (s,t)}\cap\Omega_{X^2}\cap (X_s\!\times\! Y_t)\cr
& (3)\ \mbox{diam}_{\mbox{GH}}(X_s)\mbox{, }\mbox{diam}_{\mbox{GH}}(Y_s)\mbox{, }
\mbox{diam}_{\mbox{GH}}(U_{s,t})\!\leq\! 2^{-\vert s\vert}\cr
& (4)\ X_{s0}\cap X_{s1}\! =\! Y_{s0}\cap Y_{s1}\! =\!\emptyset\cr
& (5)\ U_{s\varepsilon ,t\varepsilon}\!\subseteq\! U_{s,t}\cr
& (6)\ U_{s,t}\!\subseteq\! F_\theta\mbox{ if }(s,t)\!\in\! D_\theta
\end{array}$$
$\bullet$ Assume that this has been done. Let $\alpha\!\in\! 2^\omega$. The sequence 
$(X_{\alpha\vert n})_{n\in\omega}$ is a decreasing sequence of nonempty clopen subsets of 
$\Omega_X$ with vanishing diameters, which defines 
$f_0(\alpha )\!\in\!\bigcap_{n\in\omega}~X_{\alpha\vert n}$. As the Gandy-Harrington topology is finer than the original topology, $f_0\! :\! 2^\omega\!\rightarrow\! O_0$ is continuous. By (4), $f_0$ is injective. Similarly, we define $f_1\! :\! 2^\omega\!\rightarrow\! O_1$ injective continuous. Finally, we define $f\! :\! 2^\omega\!\rightarrow\! X$ by $f(\varepsilon\alpha )\! :=\! f_\varepsilon (\alpha )$, so that $f$ is also injective continuous since $O_0,O_1$ are disjoint.\bigskip

 If $(0\alpha ,1\beta )\!\in\!\mathbb{B}^\eta_0$, then there is $\theta\!\leq\!\eta$ of the same parity as $\eta$ such that $(\alpha ,\beta )\vert n\!\in\! D_\theta$ if $n\!\geq\! n_0$. In this case, by (1)-(3) and (5)-(6), $\big( U_{(\alpha ,\beta )\vert n}\big)_{n\geq n_0}$ is a decreasing sequence of nonempty clopen subsets of $A_0\cap\Omega_{X^2}$ with vanishing diameters, so that its intersection is a singleton $\{ F(\alpha ,\beta )\}\!\subseteq\! A_0$. As 
$(x_{\alpha\vert n},y_{\beta\vert n})$ converges (for ${\it\Sigma}_{X^2}$, and thus for 
${\it\Sigma}_X^2$) to $F(\alpha ,\beta )$, 
$\big( f(0\alpha ),f(1\beta )\big)\! =\! F(\alpha ,\beta )\!\in\! A_0$. If 
$(0\alpha ,1\beta )\!\in\!\mathbb{B}^\eta_1$, then the parity of $\theta$ is opposite to that of $\eta$ and, similarly, $\big( f(0\alpha ),f(1\beta )\big)\!\in\! A_1$.\bigskip

\noindent $\bullet$ So let us prove that the construction is possible. Note that 
$(t^0_\emptyset ,t^1_\emptyset )\! =\! (\emptyset ,\emptyset )$, 
${\cal T}^\eta\cap (2^0\!\times\! 2^0)\! =\!\{ (\emptyset ,\emptyset )\}$ and 
$(\emptyset ,\emptyset )\!\in\! D_\eta$. Let 
$(x_\emptyset ,y_\emptyset )\!\in\! F_\eta\cap\Omega_{X^2}$, and 
$\Phi (\emptyset ,\emptyset )\!\in\!\omega$ such that 
$(x_\emptyset ,y_\emptyset )\!\in\! C_{\Phi (\emptyset ,\emptyset )}$. As 
$\Omega_{X^2}\!\subseteq\!\Omega_X^2$, $x_\emptyset ,y_\emptyset\!\in\!\Omega_X$. We choose $\Ana$ subsets $X_\emptyset ,Y_\emptyset$ of $X$ with GH-diameter at most $1$ such that 
$$(x_\emptyset ,y_\emptyset )\!\in\! X_\emptyset\!\times\! Y_\emptyset\!\subseteq\! 
(\Omega_X\cap O_0)\!\times\! (\Omega_X\cap O_1)\mbox{,}$$ 
as well as a $\Ana$ subset $U_{\emptyset ,\emptyset}$ of $X^2$ with GH-diameter at most $1$ such that 
$$(x_\emptyset ,y_\emptyset )\!\in\! U_{\emptyset ,\emptyset}\!\subseteq\! F_\eta\cap 
C_{\Phi (\emptyset ,\emptyset )}\cap\Omega_{X^2}\cap (X_\emptyset\!\times\! Y_\emptyset )
\mbox{,}$$ 
which completes the construction for the length $l\! =\! 0$.

\vfill\eject

 Assume that we have constructed our objects for the sequences of length $l$. Let 
$u\!\in\!\omega^{<\omega}$ and $q\!\in\!\omega$ with $l\! +\! 1\! =<uq>_\eta$, which gives 
$w\!\in\!\omega^{<\omega}$ with $(t^0_{uq},t^1_{uq})\! =\! (t^0_uw0,t^1_uw1)$. We set\bigskip

\leftline{$U\! :=\!\{ x\!\in\! X\mid\exists (x'_s)_{s\in 2^l}\!\in\!\Pi_{s\in 2^l}~X_s~~
\exists (y'_s)_{s\in 2^l}\!\in\!\Pi_{s\in 2^l}~Y_s~~x\! =\! x'_{t^0_uw}~\wedge$}\smallskip

\rightline{$\forall (s,t)\!\in\! {\cal T}^\eta\cap (2^l\times\! 2^l)~~(x'_s,y'_t)\!\in\! U_{s,t}\}\mbox{,}$}\bigskip
 
\leftline{$V\! :=\!\{ y\!\in\! X\mid\exists (x'_s)_{s\in 2^l}\!\in\!\Pi_{s\in 2^l}~X_s~~
\exists (y'_s)_{s\in 2^l}\!\in\!\Pi_{s\in 2^l}~Y_s~~y\! =\! y'_{t^1_uw}~\wedge$}\smallskip

\rightline{$\forall (s,t)\!\in\! {\cal T}^\eta\cap (2^l\times\! 2^l)~~(x'_s,y'_t)\!\in\! U_{s,t}\} .$}\bigskip

\noindent Note that $U,V$ are $\Ana$ and $(x_{t^0_uw},y_{t^1_uw})\!\in\! 
F_{\varphi (u)}\cap (U\!\times\! V)\!\subseteq\!\bigcap_{\theta <\varphi (u)}~
\overline{F_\theta}^{\tau_1}\cap (U\!\times\! V)$. This gives 
$(x_{t^0_uw0},y_{t^1_uw1})\!\in\! F_{\varphi (uq)}\cap (U\!\times\! V)\cap\Omega_{X^2}$. Let 
$(x_{s0})_{s\in 2^l\setminus\{ t^0_uw\}}$ be witnesses for the fact that 
$x_{t^0_uw0}\!\in\! U$, and $(x_{s1})_{s\in 2^l\setminus\{ t^1_uw\}}$ be witnesses for the fact that 
$x_{t^1_uw1}\!\in\! V$.\bigskip

 We need to show that $x_{s0}\!\not=\! x_{s1}$ (and similarly for $y_{s0}$ and $y_{s1}$). First observe that if $s\!\not=\! t\!\in\! 2^l$, then $x_{s\varepsilon}\!\in\! X_s$ and 
$x_{t\varepsilon'}\!\in\! X_t$, so that $x_{s\varepsilon}\!\not=\! x_{t\varepsilon'}$ by condition 4. Similarly, $y_{s\varepsilon}\!\neq\! y_{t\varepsilon'}$. As $\varphi (u)$ and $\varphi (uq)$ do not have the same parity, there is $\epsilon\!\in\! 2$ such that 
$(x_{t^0_uw0},y_{t^1_uw1})\!\in\! A_\epsilon$ and 
$$(x_{t^0_uw1},y_{t^1_uw1})\!\in\! U_{t^0_uw,t^1_uw}\!\subseteq\! A_{1-\epsilon}.$$ 
As $A_0$ and $A_1$ are disjoint, $x_{t^0_uw0}\!\not=\! x_{t^0_uw1}$. Similarly, 
 $y_{t^0_uw0}\!\not=\! y_{t^0_uw1}$.\bigskip 

 So we may assume that $l\!\geq\! 1$ and $s\!\not=\! t^0_uw$. The fact that $G_l$ is a connected graph provides a $G_l$-path from $s$ to $t^0_uw$. This path gives us two $s(A_0\cup A_1)$-paths by the definition of $U$ and $V$, one from $y_{s0}$ to $x_{t^0_uw0}$, and another one from $y_{s1}$ to $x_{t^0_uw1}$. Moreover, the same $\Phi (s',t')$'s are involved in these two pathes since they are induced by the same $G_l$-path. Observe that 
$(x_{t^0_uw0},y_{t^1_uw1}),(x_{t^0_uw1},y_{t^1_uw1})$ are in $s(A_0\cup A_1)$. Also, since 
$x_{s\varepsilon}\in O_0$ and $y_{t\varepsilon'}\in O_1$, no ``$x$'' is equal to no ``$y$''. Thus, by quasi-acyclicity, $y_{s0}\!\neq\! y_{s1}$. Similarly, one can prove that $x_{s0}\!\neq\! x_{s1}$. The following picture illustrates the situation when $l\! =\! 1$:
$$\xymatrixrowsep{0.3in}
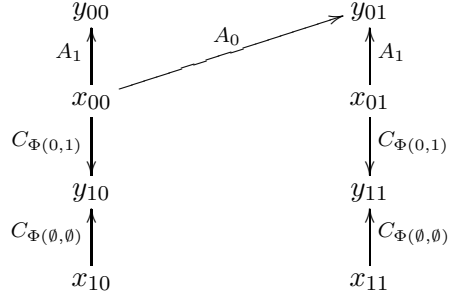
\xymatrix{ 
                       & y_{00} & & & y_{01} \\
                       & x_{00} \ar[u]^{A_1} \ar[d]_{C_{\Phi (0,1)}} \ar[urrr]^[@]{\hbox to 0pt{\hss $\scriptstyle{A_0}$\hss}} & & & x_{01} \ar[d]^{C_{\Phi (0,1)}}
                       \ar[u]_{A_1}  \\  
                       & y_{10} & & & y_{11} \\  
                       & x_{10} \ar[u]^{C_{\Phi (\emptyset ,\emptyset )}} & & & x_{11} 
                       \ar[u]_{C_{\Phi (\emptyset ,\emptyset )}} \\   }$$
Let $\Phi (t^0_uw0,t^1_uw1)\!\in\!\omega$ such that 
$(x_{t^0_uw0},y_{t^1_uw1})\!\in\! C_{\Phi (t^0_uw0,t^1_uw1)}$, and 
$\Phi (s\varepsilon ,t\varepsilon )\! :=\!\Phi (s,t)$ if $(s,t)$ is in ${\cal T}^\eta\cap (2^l\times\! 2^l)$ and $\varepsilon\!\in\! 2$. It remains to take disjoint $\Ana$ sets $X_{s0},X_{s1}\!\subseteq\! X_s$ (respectively $Y_{s0},Y_{s1}\!\subseteq\! Y_s$) with the required properties, as well as 
$V_{s\varepsilon,t\varepsilon'}$, accordingly.\bigskip

\noindent\bf Case 2\rm\ $F_\eta\!\subseteq\!\Delta (X)$.\bigskip

 Let us indicate the differences with Case 1. We set $S\! :=\!\{ x\!\in\! X\mid (x,x)\!\in\! F_\eta\}$, which is a nonempty $\Ana$ set by our assumption.\bigskip

\noindent $\bullet$ We construct the following objects:\bigskip

- a sequence $(x_s)_{s\in 2^{<\omega}}$ of points of $S$,\smallskip

- a sequence $(X_s)_{s\in 2^{<\omega}}$ of $\Ana$ subsets of $X$,\smallskip

- a sequence $(U_{s,t})_{(s,t)\in {\cal T}^\eta}$ of $\Ana$ subsets of $X^2$, and 
$\Phi\! :\! {\cal T}^\eta\!\rightarrow\!\omega$.\bigskip

\noindent We want these objects to satisfy the following conditions:
$$\begin{array}{ll}
& (1)\ x_s\!\in\! X_s\ \wedge\ (x_s,x_t)\!\in\! U_{s,t}\cr
& (2)\ X_{s\varepsilon}\!\subseteq\! X_s\!\subseteq\!\Omega_X\cap S\ \wedge\ 
U_{s,t}\!\subseteq\! C_{\Phi (s,t)}\cap\Omega_{X^2}\cap (X_s\!\times\! X_t)\cr
& (3)\ \mbox{diam}_{\mbox{GH}}(X_s)\mbox{, }\mbox{diam}_{\mbox{GH}}(U_{s,t})\!\leq\! 
2^{-\vert s\vert}\cr
& (4)\ X_{s0}\cap X_{s1}\! =\!\emptyset\cr
& (5)\ U_{s\varepsilon ,t\varepsilon}\!\subseteq\! U_{s,t}\cr
& (6)\ U_{s,t}\!\subseteq\! F_\theta\mbox{ if }(s,t)\!\in\! D_\theta
\end{array}$$
$\bullet$ Assume that this has been done. As in Case 1, we get $f\! :\! 2^\omega\!\rightarrow\! X$ injective continuous such that 
$\mathbb{N}^\eta_\epsilon\!\subseteq\! (f\!\times\! f)^{-1}(A_\epsilon )$ for each $\epsilon\!\in\! 2$.\bigskip

\noindent $\bullet$ So let us prove that the construction is possible. Let 
$(x_\emptyset ,y_\emptyset )\!\in\! F_\eta\cap\Omega_{X^2}$. As $F_\eta\!\subseteq\!\Delta (X)$, $y_\emptyset\! =\! x_\emptyset\!\in\! S$. Let $\Phi (\emptyset ,\emptyset )\!\in\!\omega$ with 
$(x_\emptyset ,x_\emptyset )\!\in\! C_{\Phi (\emptyset ,\emptyset )}$. As 
$\Omega_{X^2}\!\subseteq\!\Omega_X^2$, $x_\emptyset\!\in\!\Omega_X$. We choose a $\Ana$ subset $X_\emptyset$ of $X$ with GH-diameter at most $1$ such that 
$x_\emptyset\!\in\! X_\emptyset\!\subseteq\!\Omega_X\cap S$, as well as a $\Ana$ subset 
$U_{\emptyset ,\emptyset}$ of $X^2$ with GH-diameter at most $1$ such that 
$(x_\emptyset ,x_\emptyset )\!\in\! U_{\emptyset ,\emptyset}\!\subseteq\! F_\eta\cap 
C_{\Phi (\emptyset ,\emptyset )}\cap\Omega_{X^2}\cap (X_\emptyset\!\times\! X_\emptyset )$, which completes the construction for the length $l\! =\! 0$.\bigskip

 For the inductive step, we set\bigskip

\leftline{$U\! :=\!\{ x\!\in\! X\mid\exists (x'_s)_{s\in 2^l}\!\in\!\Pi_{s\in 2^l}~X_s~~x\! =\! x'_{t^0_uw}
\wedge\forall (s,t)\!\in\! {\cal T}^\eta\cap (2^l\times\! 2^l)~~(x'_s,x'_t)\!\in\! U_{s,t}\}\mbox{,}$}\smallskip
 
\leftline{$V\! :=\!\{ x\!\in\! X\mid\exists (x'_s)_{s\in 2^l}\!\in\!\Pi_{s\in 2^l}~X_s~~x\! =\! x'_{t^1_uw}\wedge\forall (s,t)\!\in\! {\cal T}^\eta\cap (2^l\times\! 2^l)~~(x'_s,x'_t)\!\in\! U_{s,t}\} .$}\bigskip

 Again, we need to check that $x_{t^0_q}\!\not=\! x_{t^1_q}$ if $q\!\in\!\omega$. Note first that 
$A_1\cap S^2$ is irreflexive, since otherwise it contains 
$(x,x)\!\in\! A_1\cap F_\eta\!\subseteq\! A_1\cap A_0$. By construction, 
$(x_{t^0_q},x_{t^1_q})\!\in\! F_{\varphi (q)}\!\subseteq\! A_1$, and we are done.\bigskip

\noindent (2) $\Rightarrow$ (3) Note that 
$(2^\omega ,2^\omega ,\mathbb{N}^\eta_0,\mathbb{N}^\eta_1)\sqsubseteq 
(2^\omega ,2^\omega ,\mathbb{B}^\eta_0,\mathbb{B}^\eta_1)$, with witnesses 
$\alpha\!\rightarrow\! 0\alpha$ and $\beta\!\rightarrow\! 1\beta$.\bigskip
 
\noindent (3) $\Rightarrow$ (1) This comes from Lemma \ref{S^2}.\hfill{$\square$}

\begin{prop} \label{squareD2} Let $\eta$ be a countable ordinal. The pairs 
$(\mathbb{N}^\eta_0,\mathbb{N}^\eta_1)$ and $(\mathbb{B}^\eta_0,\mathbb{B}^\eta_1)$ are incomparable for the square reduction.\end{prop}

\noindent\bf Proof.\rm ~There is no map $f\! :\! 2^\omega\!\rightarrow\! 2^\omega$ such that 
$\mathbb{N}^\eta_\varepsilon\!\subseteq\! (f\!\times\! f)^{-1}(\mathbb{B}^\eta_\varepsilon)$ since 
$\Delta (2^\omega )$ is a subset of $\mathbb{N}^\eta_0$.\bigskip

 There is no injection $f\! :\! 2^\omega\!\rightarrow\! 2^\omega$ for which there is $\alpha\!\in\! 2^\omega$ such that $f(0\alpha )\! =\! f(1\alpha )$. Using this fact, assume, towards a contradiction, that there is $f\! :\! 2^\omega\!\rightarrow\! 2^\omega$ injective continuous such that 
$\mathbb{B}^\eta_\varepsilon\!\subseteq\! (f\!\times\! f)^{-1}(\mathbb{N}^\eta_\varepsilon)$. Let $(0t^0_s\gamma,1t^1_s\gamma)\!\in\!\mathbb{B}^\eta_\varepsilon$, so that 
$\big( f(0t^0_s\gamma),f(1t^1_s\gamma)\big)\! =\! 
(t^0_v\gamma',t^1_v\gamma')\!\in\!\mathbb{N}^\eta_{\varepsilon}$.\bigskip

 We claim that $\varphi (s)\leq \varphi (v)$. We proceed by induction on $\varphi (s)$. Notice that is is obvious for $\varphi (s)=0$. Suppose that it holds for all $\theta\! <\!\varphi (s)$. Note that we can find $p_k\!\in\!\omega$ and $\gamma_k\!\in\! 2^\omega$ such that 
$(t^0_{sp_k}\gamma_k,t^1_{sp_k}\gamma_k)\!\in\!\mathbb{N}^\eta_{1-\varepsilon}$ and 
$(t^0_{sp_k}\gamma_k,t^1_{sp_k}\gamma_k)\!\rightarrow\! (t^0_s\gamma,t^1_s\gamma)$. By continuity, 
$$(t^0_{v_k}\gamma',t^1_{v_k}\gamma')\! :=\! 
\big( f(0t^0_{sp_k}\gamma_k),f(1t^1_{sp_k}\gamma_k)\big)\!\rightarrow\! 
(t^0_v\gamma',t^1_v\gamma').$$

 In particular, for $k$ large, $(t^0_v,t^1_v)\!\subseteq\! (t^0_{v_k},t^1_{v_k})$. This implies that the sequence $v_k$ is a strict extension of $v$. Therefore $\varphi (v_k)\! <\!\varphi (v)$. By the induction hypothesis, 
$\varphi (sp_k)\!\leq\!\varphi (v_k)\! <\!\varphi (v)$. If $\varphi(s)\! =\!\theta\! +\! 1$, then 
$\theta\! =\!\varphi (sp_k)\! <\!\varphi (v)$, so we are done. If $\varphi (s)$ is a limit ordinal, then 
$\big(\varphi (sp_k)\big)_{k\in\omega}$ is cofinal in it, so we are done too.\bigskip

 Finally, let $\alpha\in 2^\omega$, so that 
$(0\alpha,1\alpha)\! =\! (0t^0_\emptyset\alpha,1t^1_\emptyset\alpha)\!\in\!\mathbb{B}^\eta_0$. Then $\big( f(0\alpha),f(1\alpha)\big)\! =\! (t^0_v\gamma',t^1_v\gamma')$ with $\varphi (v)\! =\!\eta$, so that $v\! =\!\emptyset$, which contradicts the injectivity of $f$.\hfill{$\square$}\bigskip

\noindent $\underline{\mbox{\bf Consequences}}$

\begin{lem} \label{anaKsigmas} Let $\bf\Gamma$ be a class of sets contained in $\bortwo$ which is either a Wadge class or $\bortwo$, $X$ be a Polish space, and $A,B$ be disjoint analytic relations on $X$. Then exactly one of the following holds:\smallskip  

(a) the set $A$ is separable from $B$ by a $\mbox{pot}({\bf\Gamma})$ set,\smallskip  

(b) there are $K_\sigma$ sets $A'\!\subseteq\! A$ and $B'\!\subseteq\! B$ such that $A'$ is not separable from $B'$ by a $\mbox{pot}({\bf\Gamma})$ set.\end{lem}

\noindent\bf Proof.\rm ~Assume that (a) does not hold. Theorems 1.9 and 1.10 in [L8] give $\boratwo$ relations $\mathbb{S}_0,\mathbb{S}_1$ on $2^\omega$ and $g,h\! :\! 2^\omega\!\rightarrow\! X$ continuous with $\mathbb{S}_0\!\subseteq\! (g\!\times\! h)^{-1}(A)$ and 
$\mathbb{S}_1\!\subseteq\! (g\!\times\! h)^{-1}(B)$. We set $A'\! :=\! (g\!\times\! h)\big[\mathbb{S}_0\big]$ and 
$B'\! :=\! (g\!\times\! h)\big[\mathbb{S}_1\big]$.\hfill{$\square$}

\begin{cor} \label{cor1checkD2Sigma01} Let $\eta\! <\!\omega_1$, $X$ be a Polish space, and $A,B$ be disjoint analytic relations on $X$ such that $A\cup B$ is s-acyclic or locally countable. Then exactly one of the following holds:\smallskip  

(a) the set $A$ is separable from $B$ by a $\mbox{pot}\big( D_\eta (\boraone )\big)$ set,\smallskip  

(b) $(2^\omega ,2^\omega ,\mathbb{N}^\eta_0,\mathbb{N}^\eta_1)\sqsubseteq (X,X,A,B)$ if 
$\eta\!\geq\! 1$ and $(1,1,\mathbb{N}^\eta_0,\mathbb{N}^\eta_1)\sqsubseteq (X,X,A,B)$ if 
$\eta\! =\! 0$.\end{cor}

\noindent\bf Proof.\rm ~By Lemma \ref{S^2}, $\mathbb{N}^\eta_0$ is not separable from 
$\mathbb{N}^\eta_1$ by a $\mbox{pot}\big( D_\eta (\boraone )\big)$ set. This shows that (a) and (b) cannot hold simultaneously. So assume that (a) does not hold. We may assume that 
$\eta\!\geq\! 1$. By Lemma \ref{anaKsigmas}, we may assume that $A,B$ are $\boratwo$. By Lemma \ref{suffqa}, we may also assume that $A\cup B$ is quasi-acyclic. It remains to apply Theorem \ref{checkD2Sigma01}.\hfill{$\square$}

\begin{cor} \label{cor2checkD2Sigma01} Let $\eta$ be a countable ordinal, $X,Y$ be Polish spaces, and $A,B$ be disjoint analytic subsets of $X\!\times\! Y$ such that $A\cup B$ is locally countable. Then exactly one of the following holds:\smallskip  

(a) the set $A$ is separable from $B$ by a $\mbox{pot}\big( D_\eta (\boraone )\big)$ set,\smallskip  

(b) $(2^\omega ,2^\omega ,\mathbb{N}^\eta_0,\mathbb{N}^\eta_1)\sqsubseteq (X,Y,A,B)$ if 
$\eta\!\geq\! 1$ and $(1,1,\mathbb{N}^\eta_0,\mathbb{N}^\eta_1)\sqsubseteq (X,Y,A,B)$ if 
$\eta\! =\! 0$.\end{cor}

\noindent\bf Proof.\rm ~We may assume that $\eta\!\geq\! 1$. As in the proof of Corollary 
\ref{cor1checkD2Sigma01}, (a) and (b) cannot hold simultaneously. So assume that (a) does not hold. We put $Z\! :=\! X\!\oplus\! Y$, $A'\! :=\!\big\{\big( (x,0),(y,1)\big)\!\in\! Z^2\mid (x,y)\!\in\! A\big\}$ and $B'\! :=\!\big\{\big( (x,0),(y,1)\big)\!\in\! Z^2\mid (x,y)\!\in\! B\big\}$. Then $Z$ is  Polish, $A',B'$ are disjoint analytic relations on $Z$, $A'\cup B'$ is locally countable, and $A'$ is not separable from $B'$ by a $\mbox{pot}\big( D_\eta (\boraone )\big)$ set.

\vfill\eject

 Corollary \ref{cor1checkD2Sigma01} gives $f',g'\! :\! 2^\omega\!\rightarrow\! Z$ injective continuous such that ${\mathbb{N}^\eta_0\!\subseteq\! (f'\!\times\! g')^{-1}(A')}$, and also 
$\mathbb{N}^\eta_1\!\subseteq\! (f'\!\times\! g')^{-1}(B')$. We set 
$f(\alpha )\! :=\!\Pi_0[f'(\alpha )]$, and $g(\beta )\! :=\!\Pi_0[g'(\beta )]$. As 
$\Delta (2^\omega )\!\subseteq\!\mathbb{N}^\eta_0$, $f'$ takes values in $X\!\times\!\{ 0\}$ and $g'$ takes values in $Y\!\times\!\{ 1\}$. This implies that $f\! :\! 2^\omega\!\rightarrow\! X$, 
$g\! :\! 2^\omega\!\rightarrow\! Y$ are injective continuous. We are done since 
$\mathbb{N}^\eta_0\!\subseteq\! (f\!\times\! g)^{-1}(A)$ and 
$\mathbb{N}^\eta_1\!\subseteq\! (f\!\times\! g)^{-1}(B)$.\hfill{$\square$}\bigskip

\noindent\bf Notation.\rm ~If $A$ is a relation on $2^\omega$, then we set  
$G_A\! :=\!\{ (0\alpha ,1\beta )\mid (\alpha ,\beta )\!\in\! A\}$.

\begin{lem} \label{suffacy} Let $A$ be an antisymmetric s-acyclic relation on $2^\omega$. Then $G_A$ is s-acyclic.\end{lem}

\noindent\bf Proof.\rm ~We argue by contradiction, which gives $n\!\geq\! 2$ and an injective 
$s(G_A)$-path $(\varepsilon_iz_i)_{i\leq n}$ such that 
$(\varepsilon_0z_0,\varepsilon_nz_n)\!\in\! s(G_A)$. This implies that 
$\varepsilon_i\!\not=\!\varepsilon_{i+1}$ if $i\! <\! n$ and $n$ is odd. Thus $(z_i)_{i\leq n}$ is a 
$s(A)$-path such that $(z_{2j})_{2j\leq n}$ and $(z_{2j+1})_{2j+1\leq n}$ are injective and 
$(z_0,z_n)\!\in\! s(A)$. As $s(A)$ is acyclic, the sequence $(z_i)_{i\leq n}$ is not injective. We erase $z_{2j+1}$ from this sequence if $z_{2j+1}\!\in\!\{ z_{2j},z_{2j+2}\}$ and $2j\! +\!1\!\leq\! n$, which gives a sequence $(z'_i)_{i\leq n'}$ which is still a $s(A)$-path with $(z'_0,z'_{n'})\!\in\! s(A)$, and moreover satisfies $z'_i\!\not=\! z'_{i+1}$ if $i\! <\! n'$.\bigskip

 If $n'\! <\! 2$, then $n\! =\! 3$, $z_0\! =\! z_1$ and $z_2\! =\! z_3$. As $A$ is antisymmetric and 
$\varepsilon_3\! =\!\varepsilon_1\!\not=\!\varepsilon_2\! =\!\varepsilon_0$, we get $z_0\! =\! z_2$, which is absurd. If $n'\!\geq\! 2$, then $(z'_i)_{i\leq n'}$ is not injective again. We choose a subsequence of it with at least three elements,  made of consecutive elements, such that the first and the last elements are equal, and of minimal length with these properties. The acyclicity of 
$s(A)$ implies that this subsequence has exactly three elements, say 
$(z'_i,z'_{i+1},z'_{i+2}\! =\! z'_i)$.\bigskip

 If $z'_i\! =\! z_{2j+1}$, then $z'_{i+1}\! =\! z_{2j+2}$, $z'_{i+2}\! =\! z_{2j+4}$ and 
$z_{2j+3}\! =\! z_{2j+2}$.  As $A$ is antisymmetric and 
$\varepsilon_{2j+3}\! =\!\varepsilon_{2j+1}\!\not=\!\varepsilon_{2j+2}\! =\!\varepsilon_{2j+4}$, we get $z_{2j+2}\! =\! z_{2j+4}$, which is absurd. If $z'_i\! =\! z_{2j}$, then $z'_{i+1}\! =\! z_{2j+2}$, and $z'_{i+2}\! =\! z_{2j+3}$. As $A$ is antisymmetric and 
$\varepsilon_{2j+3}\! =\!\varepsilon_{2j+1}\!\not=\!\varepsilon_{2j+2}\! =\!\varepsilon_{2j}$, we get $z_{2j}\! =\! z_{2j+2}$, which is absurd.\hfill{$\square$}

\begin{cor} \label{caracpartialcheckD2Sigma01} Let $\eta\!\geq\! 1$ be a countable ordinal, $X$ be a Polish space, and $A,B$ be disjoint analytic relations on $X$. The following are equivalent:\smallskip

\noindent (1) there is an s-acyclic relation $R\!\in\!\ana$ such that $A\cap R$ is not separable from $B\cap R$ by a $\mbox{pot}\big( D_\eta (\boraone )\big)$ set,\smallskip

\noindent (2) there is a locally countable relation $R\!\in\!\ana$ such that $A\cap R$ is not separable from $B\cap R$ by a $\mbox{pot}\big( D_\eta (\boraone )\big)$ set,\smallskip

\noindent (3) 
$(2^\omega ,2^\omega ,\mathbb{N}^\eta_0,\mathbb{N}^\eta_1)\sqsubseteq (X,X,A,B)$,\smallskip

\noindent (4) there is $(\mathbb{A}_0,\mathbb{A}_1)\!\in\!
\{ (\mathbb{N}^\eta_0,\mathbb{N}^\eta_1),(\mathbb{B}^\eta_0,\mathbb{B}^\eta_1)\}$ such that 
$(2^\omega ,2^\omega ,\mathbb{A}_0,\mathbb{A}_1)\sqsubseteq (X,X,A,B)$, via a square map.\smallskip

\noindent A similar result holds for $\eta\! =\! 0$ with $1$ instead of $2^\omega$.\end{cor}

\noindent\bf Proof.\rm ~(1) $\Rightarrow$ (3),(4) and (2) $\Rightarrow$ (3),(4) This is a consequence of Corollary \ref{cor1checkD2Sigma01} and its proof.\bigskip

\noindent (4) $\Rightarrow$ (1) By the remarks before Lemma \ref{S^2}, 
$\mathbb{N}^\eta_0\cup\mathbb{N}^\eta_1$ has s-acyclic levels. This implies that 
$\mathbb{N}^\eta_0\cup\mathbb{N}^\eta_1$ is s-acyclic. As 
$\mathbb{N}^\eta_0\cup\mathbb{N}^\eta_1$ is antisymmetric, 
$\mathbb{B}^\eta_0\cup\mathbb{B}^\eta_1$ is s-acyclic too, by Lemma \ref{suffacy}. Thus we can take $R\! :=\! (f\!\times\! f)[\mathbb{A}_0\cup\mathbb{A}_1]$ since the s-acyclicity is preserved by images by the square of an injection, and by Lemma \ref{S^2}.\bigskip

\noindent (4) $\Rightarrow$ (2) We can take $R\! :=\! (f\!\times\! f)[\mathbb{A}_0\cup\mathbb{A}_1]$ since $\mathbb{A}_0\cup\mathbb{A}_1$ is locally countable, by Lemma \ref{S^2}.\bigskip

\noindent (3) $\Rightarrow$ (2) We can take 
$R\! :=\! (f\!\times\! f)[\mathbb{N}^\eta_0\cup\mathbb{N}^\eta_1]$ since 
$\mathbb{N}^\eta_0\cup\mathbb{N}^\eta_1$ is locally countable, by Lemma \ref{S^2}.
\hfill{$\square$}\bigskip

\noindent\bf Remark.\rm ~There is a version of Corollary \ref{caracpartialcheckD2Sigma01} for 
$\check D_\eta (\boraone )$ instead of $D_\eta (\boraone )$, obtained by exchanging the roles of 
$A$ and $B$. This symmetry is also  present in Theorem \ref{checkD2Sigma01}.\bigskip

 We now give some complements when $\eta\! =\! 1$. At the beginning of this section, we mentioned the fact that our examples are in the style of $\mathbb{G}_0$. If $\eta\! =\! 1$, then $\mathbb{G}_0$ itself is involved.
 
\begin{cor} \label{cor1pi01} Let $X$ be a Polish space, and $A,B$ be disjoint analytic relations on $X$ such that\smallskip

- either $A\cup B$ is s-acyclic or locally countable,\smallskip

- or $A$ is contained in a potentially closed s-acyclic or locally countable relation.\smallskip

\noindent Then exactly one of the following holds:\smallskip  

(a) the set $A$ is separable from $B$ by a $\mbox{pot}(\bormone )$ set,\smallskip  

(b) $\big( 2^\omega ,2^\omega ,\mathbb{G}_0,\Delta (2^\omega )\big)\sqsubseteq (X,X,A,B)$.
\end{cor}

\begin{cor} \label{cor2pi01} Let $X,Y$ be Polish spaces, and $A,B$ be disjoint analytic subsets of 
$X\!\times\! Y$ such that $A\cup B$ is locally countable or $A$ is contained in a potentially closed locally countable set. Then exactly one of the following holds:\smallskip  

(a) the set $A$ is separable from $B$ by a $\mbox{pot}(\bormone )$ set,\smallskip  

(b) $\big( 2^\omega ,2^\omega ,\mathbb{G}_0,\Delta (2^\omega )\big)\sqsubseteq (X,Y,A,B)$.
\end{cor}

\section{$\!\!\!\!\!\!$ The class $\Delta\big( D_\eta (\boraone )\big)$}

$\underline{\mbox{\bf Examples}}$\bigskip

\noindent\bf Notation.\rm ~We set, for each countable ordinal $\eta\!\geq\! 1$ and each 
$\varepsilon\!\in\! 2$, 
$$\mathbb{S}^\eta_\varepsilon\! :=\!\Big\{ (t^0_s\gamma ,t^1_s\gamma )\mid 
s\!\in\! T_\eta\!\setminus\!\{\emptyset\}\wedge
\mbox{parity}(\vert s\vert )\! =\! 1\! -\!\big\vert\mbox{parity}\big( s(0)\big)\! -\!\varepsilon\big\vert\wedge\gamma\!\in\! 2^\omega\Big\} .$$ 

\begin{lem} \label{deltaeta1} Let $\eta\!\geq\! 1$ be a countable ordinal, and $C$ be a nonempty clopen subset of $2^\omega$. Then 
$\mathbb{S}^\eta_0\cap C^2$ is not separable from $\mathbb{S}^\eta_1\cap C^2$ by a $\mbox{pot}\Big(\Delta\big( D_\eta (\boraone )\big)\Big)$ set.\end{lem}

\noindent\bf Proof.\rm ~We use the notation in the proof of Lemma \ref{S^2}. We argue by contradiction, which gives $P$ in 
$\mbox{pot}\Big(\Delta\big( D_\eta (\boraone )\big)\Big)$, and a dense $G_\delta$ subset of $2^\omega$ such that 
$P\cap G^2,G^2\!\setminus\! P\!\in\! D_\eta (\boraone )(G^2)$. So let, for each $\varepsilon\!\in\! 2$, $(O^\varepsilon_\theta )_{\theta <\eta}$ be a sequence of open relations on $2^\omega$ such that 
$$P\cap G^2\! =\!\big(\bigcup_{\theta <\eta ,\mbox{parity}(\theta )\not=\mbox{parity}(\eta )}~
O^0_\theta\!\setminus\! (\bigcup_{\theta'<\theta}~O^0_{\theta'})\big)\cap G^2$$ 
and $G^2\!\setminus\! P\! =\!\big(\bigcup_{\theta <\eta ,\mbox{parity}(\theta )\not=\mbox{parity}(\eta )}~
O^1_\theta\!\setminus\! (\bigcup_{\theta'<\theta}~O^1_{\theta'})\big)\cap G^2$.\bigskip

\noindent $\bullet$ Note that $\mathbb{S}^\eta_\varepsilon\! =\!\bigcup_{s\in T_\eta\setminus\{\emptyset\} ,\mbox{parity}(\vert s\vert )=1-\vert\mbox{parity}(s(0))-\varepsilon\vert}~
\mbox{Gr}(f_s)$. Let us show that if $\theta\!\leq\!\eta$, $s\!\in\! T_\eta$ and $\varphi (s)\! =\!\theta$, then $\mbox{Gr}(f_s)\cap (C\cap G)^2\!\subseteq\!\neg O^{1-\mbox{parity}(s(0))}_\theta$ if 
$\theta\! <\!\eta$, and $\mbox{Gr}(f_s)\cap (C\cap G)^2$ is disjoint from 
$\bigcup_{\theta'<\theta}~(O^0_{\theta'}\cup O^1_{\theta'})$ if $\theta\! =\!\eta$. The objects 
$s\! =\!\emptyset$ and $\theta\! =\!\eta$ will give the contradiction.\bigskip

\noindent $\bullet$ We argue by induction on $\theta$. Note that 
$\mbox{Gr}(f_s)\cap (C\cap G)^2\!\subseteq\!
\mathbb{S}^\eta_{1-\vert\mbox{parity}(\vert s\vert )-\mbox{parity}(s(0))\vert}\cap G^2$ if 
$\theta\! =\! 0$ since $s\!\not=\!\emptyset$. As $\mathbb{S}^\eta_\varepsilon\cap G^2\!\subseteq\!\neg O^{\vert\mbox{parity}(\eta )-\varepsilon\vert}_0$ for each $\varepsilon\!\in\! 2$ and 
$\vert s\vert$ has the same parity as $\eta$ if $\theta =0$, we are done.\bigskip

\noindent $\bullet$ Assume that the result has been proved for $\theta'\! <\!\theta$. If $\theta$ is the successor of $\theta'$, then the induction assumption implies that 
$\mbox{Gr}(f_{sq})\cap (C\cap G)^2\!\subseteq\!\neg O^{1-\mbox{parity}((sq)(0))}_{\theta'}$ for each $q$. We set, for each $\varepsilon\!\in\! 2$, 
$C_s^\varepsilon\! :=\!\bigcup_{k\in\omega}~\mbox{Gr}(f_{s(2k+\varepsilon)})$, so that 
$\mbox{Gr}(f_s)\!\subseteq\!\overline{C^\varepsilon_s}$, by the choice of $\psi$. If 
$s\! =\!\emptyset$, then 
$$C^\varepsilon_\emptyset\cap (C\cap G)^2\!\subseteq\!\neg O^{1-\varepsilon}_{\theta'}\mbox{,}$$ 
$\mbox{Gr}(f_s)\cap (C\cap G)^2\!\subseteq\!
\overline{C^\varepsilon_\emptyset}\cap (C\cap G)^2\!\subseteq\!
\overline{C^\varepsilon_\emptyset\cap (C\cap G)^2}\!\subseteq\!\neg O^{1-\varepsilon}_{\theta'}$, which gives the desired inclusion for $\theta\! =\!\eta$.\bigskip

 If $s\!\not=\!\emptyset$, then 
$\mbox{Gr}(f_{sq})\cap (C\cap G)^2\!\subseteq\!\neg O^{1-\mbox{parity}(s(0))}_{\theta'}$ for each $q$, so that 
$$\mbox{Gr}(f_s)\cap (C\cap G)^2\!\subseteq\!\overline{C_s}\cap (C\cap G)^2\!\subseteq\!
\overline{C_s\cap (C\cap G)^2}\!\subseteq\!\neg O^{1-\mbox{parity}(s(0))}_{\theta'}.$$ 
Thus 
$$\mbox{Gr}(f_s)\cap (C\cap G)^2\!\subseteq\! (G^2\!\setminus\! O^{1-\mbox{parity}(s(0))}_{\theta'})\cap\neg (O^{1-\mbox{parity}(s(0))}_\theta\!\setminus\! O^{1-\mbox{parity}(s(0))}_{\theta'})\!\subseteq\!\neg O^{1-\mbox{parity}(s(0))}_\theta$$ 
since $\mbox{parity}(\theta )\! =\!\vert\mbox{parity}(\vert s\vert )\! -\!\mbox{parity}(\eta )\vert$.\bigskip

\noindent $\bullet$ If $\theta$ is limit, then $\big(\varphi (sn)\big)_{n\in\omega}$ is cofinal in 
$\varphi (s)$, and 
$\mbox{Gr}(f_{sn})\cap (C\cap G)^2\!\subseteq\!\neg O^{1-\mbox{parity}((sn)(0))}_{\varphi (sn)}$, by the induction assumption. If $\theta_0\! <\!\varphi (s)$, then there is $n(\theta_0)$ such that 
$\varphi (sn)\! >\!\theta_0$ if $n\!\geq\! n(\theta_0)$. Thus 
$\mbox{Gr}(f_{sn})\cap (C\cap G)^2\!\subseteq\!\neg O^{1-\mbox{parity}((sn)(0))}_{\theta_0}$ if 
$n\!\geq\! n(\theta_0)$. If $s\! =\!\emptyset$, then, for each $\varepsilon\!\in\! 2$,    
$$\begin{array}{ll}
\mbox{Gr}(f_s)\cap (C\cap G)^2\!\!\!\!
& \subseteq\! (C\cap G)^2\cap\overline{C^\varepsilon_s}\!\setminus\! 
C^\varepsilon_s\! =\!\overline{C^\varepsilon_s\cap (C\cap G)^2}\!\setminus\! C^\varepsilon_s\cr
& \subseteq\!\overline{\bigcup_{n\geq n(\theta_0),\mbox{parity}(n)=\varepsilon}~
\mbox{Gr}(f_{sn})\cap (C\cap G)^2}\!\subseteq\!\neg O^{1-\varepsilon}_{\theta_0} .
\end{array}$$ 
Thus $\mbox{Gr}(f_s)\cap (C\cap G)^2\!\subseteq\!\neg\big(\bigcup_{\theta'<\eta}~
(O^0_{\theta'}\cup O^1_{\theta'})\big)$. If $s\!\not=\!\emptyset$, then 
$\mbox{Gr}(f_{sn})\cap (C\cap G)^2\!\subseteq\!\neg O^{1-\mbox{parity}(s(0))}_{\theta_0}$ for each $n$, so that $\mbox{Gr}(f_s)\cap (C\cap G)^2\!\subseteq\!\overline{C_s}\cap (C\cap G)^2
\!\subseteq\!\overline{C_s\cap (C\cap G)^2}\!\subseteq\!\neg O^{1-\mbox{parity}(s(0))}_{\theta_0}$. As $\mbox{parity}(\vert s\vert )\! =\!\mbox{parity}(\eta )$, 
$\mbox{Gr}(f_s)\cap (C\cap G)^2\!\subseteq\!\neg O^{1-\mbox{parity}(s(0))}_\theta$ as above.\hfill{$\square$}\bigskip

\noindent $\underline{\mbox{\bf A topological characterization}}$\bigskip

\noindent\bf Notation.\rm ~We define, for $1\!\leq\!\xi\! <\!\omega_1^{\mbox{CK}}$ and 
$\eta\! <\!\omega_1^{\mbox{CK}}$, $\bigcap_{\theta <0}\ G_{\theta ,\xi}\! :=\! (\omega^\omega )^2$, and, inductively,  
$$G_{\eta ,\xi}\! :=\!\left\{\!\!\!\!\!\!
\begin{array}{ll}
& \bigcap_{\theta <\eta}\ G_{\theta ,\xi}\mbox{ if }\eta\mbox{ is limit (possibly $0$),}\cr
& \overline{A_0\cap G_{\theta ,\xi}}^{\tau_\xi}\cap\overline{A_1\cap G_{\theta ,\xi}}^{\tau_\xi}
\mbox{ if }\eta\! =\!\theta\! +\! 1.
\end{array}
\right.$$

\begin{thm} \label{kernel2} Let $1\!\leq\!\xi\! <\!\omega_1^{\mbox{CK}}$, 
$1\!\leq\!\eta\! =\!\lambda\! +\! 2k\! +\!\varepsilon\! <\!\omega_1^{\mbox{CK}}$ with $\lambda$ limit, 
$k\!\in\!\omega$ and $\varepsilon\!\in\! 2$, and 
$A_0$, $A_1$ be disjoint $\Ana$ relations on $\omega^\omega$. Then the following are equivalent:\smallskip

\noindent (1) the set $A_0$ is not separable from $A_1$ by a 
$\mbox{pot}\Big(\Delta\big( D_\eta (\boraxi )\big)\Big)$ set,\smallskip

\noindent (2) the $\Ana$ set $G_{\eta ,\xi}$ is not empty.\end{thm}

\noindent\bf Proof.\rm ~The proof is in the spirit of that of Theorem \ref{kernel1}. The proof of Theorem 1.10.(2) in [L8] gives $\alpha$ suitable such that $c(\alpha )$ codes the class 
$D_\eta (\boraxi )$. By Theorem 6.26 in [L8] and Theorem \ref{kernel1}, (1) is equivalent to 
$R'(\alpha ,a_0,a_1)\!\not=\!\emptyset$, where 
$$R'(\alpha ,a_0,a_1)\! :=\!\left\{\!\!\!\!\!\!
\begin{array}{ll}
 & F^0_{\theta ,\xi}\cap F^1_{\theta ,\xi}\mbox{ if }\eta\! =\!\theta\! +\! 1\mbox{,}\cr
 & \bigcap_{p\geq 1}~F^0_{\theta_p,\xi}\mbox{ if }\eta\! =\!\mbox{sup}_{p\geq 1}~\theta_p
 \mbox{ is limit }\wedge\theta_p\mbox{ is odd.}
\end{array}
\right.$$ 
So it is enough to prove that 
$$G_{\eta ,\xi}\! =\!\left\{\!\!\!\!\!\!
\begin{array}{ll}
 & F^0_{\theta ,\xi}\cap F^1_{\theta ,\xi}\mbox{ if }\eta\! =\!\theta\! +\! 1\mbox{,}\cr
 & \bigcap_{p\geq 1}~F^0_{\theta_p,\xi}\mbox{ if }\eta\! =\!\mbox{sup}_{p\geq 1}~\theta_p
 \mbox{ is limit }\wedge\theta_p\mbox{ is odd.}
\end{array}
\right.$$
We argue by induction on $\eta$. Note first that 
$G_{1 ,\xi}\! =\!\overline{A_0}\cap\overline{A_1}\! =\! F^0_{0 ,\xi}\cap F^1_{0 ,\xi}$. Then, inductively,
$$\begin{array}{ll}
G_{\theta +2,\xi}\!\!\!
& \! =\!\overline{A_0\cap G_{\theta +1,\xi}}\cap\overline{A_1\cap G_{\theta +1,\xi}}\! =\!
\overline{A_0\cap F^0_{\theta ,\xi}\cap F^1_{\theta ,\xi}}\cap
\overline{A_1\cap F^0_{\theta ,\xi}\cap F^1_{\theta ,\xi}}\cr
& \! =\!\overline{A_0\cap F^{1-\mbox{parity}(\theta )}_{\theta,\xi}}\cap
\overline{A_1\cap F^{\mbox{parity}(\theta )}_{\theta ,\xi}}
\! =\! F^0_{\theta +1,\xi}\cap F^1_{\theta +1,\xi}.
\end{array}$$
If $\lambda$ is limit, then 
$$\begin{array}{ll}
G_{\lambda +1,\xi}\!\!\!
& \! =\!\overline{A_0\cap G_{\lambda ,\xi}}\cap\overline{A_1\cap G_{\lambda ,\xi}}\! =\!
\overline{A_0\cap\bigcap_{\theta <\lambda}~G_{\theta ,\xi}}\cap
\overline{A_1\cap \bigcap_{\theta <\lambda}~G_{\theta ,\xi}}\cr
& \! =\!\overline{A_0\cap\bigcap_{\theta <\lambda}~G_{\theta+1,\xi}}\cap
\overline{A_1\cap\bigcap_{\theta <\lambda}~G_{\theta+1,\xi}}\cr
& \! =\!\overline{A_0\cap\bigcap_{\theta <\lambda}~F^0_{\theta ,\xi}\cap F^1_{\theta ,\xi}}\cap
\overline{A_1\cap\bigcap_{\theta <\lambda}~F^0_{\theta ,\xi}\cap F^1_{\theta ,\xi}}\cr
& \! =\!\overline{A_0\cap\bigcap_{\theta <\lambda}~F^0_{\theta ,\xi}}\cap
\overline{A_1\cap \bigcap_{\theta <\lambda}~F^1_{\theta ,\xi}}\! =\! F^0_{\lambda ,\xi}\cap 
F^1_{\lambda ,\xi} 
\end{array}$$
and $G_{\lambda ,\xi}\! =\!\bigcap_{\theta <\lambda}~G_{\theta ,\xi}\! =\!
\bigcap_{\theta <\lambda}~G_{\theta+1,\xi}\! =\!
\bigcap_{\theta <\lambda}~F^0_{\theta ,\xi}\cap F^1_{\theta ,\xi}\! =\!
\bigcap_{\theta <\lambda}~F^0_{\theta ,\xi}\! =\!
\bigcap_{p\geq 1}~F^0_{\theta_p,\xi}$.\hfill{$\square$}\bigskip

\noindent $\underline{\mbox{\bf The main result}}$\bigskip

 We prove a version of Theorem \ref{checkD2Sigma01} for the class 
$\Delta\big( D_\eta (\boraone )\big)$. We set, for $1\!\leq\!\eta\! <\!\omega_1$ and 
$\varepsilon\!\in\! 2$, $\mathbb{C}^\eta_\varepsilon\! :=\!
\{ (0\alpha ,1\beta )\mid (\alpha ,\beta )\!\in\!\mathbb{S}^\eta_\varepsilon\}$.

\begin{thm} \label{exdeltaeta1} Let $\eta\!\geq\! 1$ be a countable ordinal, $X$ be a Polish space, and $A_0,A_1$ be disjoint analytic relations on $X$ such that $A_0\cup A_1$ is contained in a potentially closed quasi-acyclic relation. The following are equivalent:\smallskip  

\noindent (1) the set $A_0$ is not separable from $A_1$ by a 
$\mbox{pot}\Big(\Delta\big( D_\eta (\boraone )\big)\Big)$ set,\smallskip  

\noindent (2) there is $(\mathbb{A}_0,\mathbb{A}_1)\!\in\!
\{ (\mathbb{N}^\eta_1,\mathbb{N}^\eta_0),(\mathbb{B}^\eta_1,\mathbb{B}^\eta_0),
(\mathbb{N}^\eta_0,\mathbb{N}^\eta_1),(\mathbb{B}^\eta_0,\mathbb{B}^\eta_1),
(\mathbb{S}^\eta_0,\mathbb{S}^\eta_1),(\mathbb{C}^\eta_0,\mathbb{C}^\eta_1)\}$ for which the inequality $(2^\omega ,2^\omega ,\mathbb{A}_0,\mathbb{A}_1)\!\sqsubseteq\! (X,X,A_0,A_1)$ holds, via a square map,\smallskip

\noindent (3) there is $(\mathbb{A}_0,\mathbb{A}_1)\!\in\!
\{ (\mathbb{N}^\eta_1,\mathbb{N}^\eta_0),(\mathbb{N}^\eta_0,\mathbb{N}^\eta_1),
(\mathbb{S}^\eta_0,\mathbb{S}^\eta_1)\}$ such that 
$(2^\omega ,2^\omega ,\mathbb{A}_0,\mathbb{A}_1)\!\sqsubseteq\! (X,X,A_0,A_1)$.\end{thm}

\noindent\bf Proof.\rm ~(1) $\Rightarrow$ (2) The proof is partly similar to that of Theorem 
\ref{checkD2Sigma01}. Let $R$ be a potentially closed quasi-acyclic relation containing 
$A_0\cup A_1$, and $(C_n)_{n\in\omega}$ be a witness for the fact that $R$ is quasi-acyclic. We may assume that $X$ is zero-dimensional (and thus a closed subset of $\omega^\omega$) and $R$ is closed. In fact, we may assume that $X\! =\!\omega^\omega$. Indeed, as $A_0$ is not separable from $A_1$ by a $\mbox{pot}\Big(\Delta\big( D_\eta (\boraone )\big)\Big)$ set in $X^2$, it is also the case in $(\omega^\omega )^2$, which gives 
$f\! :\! 2^\omega\!\rightarrow\!\omega^\omega$. Note that 
$$\mathbb{A}_0\cup\mathbb{A}_1\!\subseteq\! (f\!\times\! f)^{-1}(A_0\cup A_1)\!\subseteq\! 
(f\!\times\! f)^{-1}(X^2)\mbox{,}$$ 
which implies that $\overline{\mathbb{A}_0\cup\mathbb{A}_1}\!\subseteq\! (f\!\times\! f)^{-1}(X^2)$. As $\Delta (2^\omega )\!\subseteq\!
\mathbb{N}^\eta_0\cap\overline{\mathbb{S}^\eta_0\cup\mathbb{S}^\eta_1}$ and 
$$\{ (0\alpha ,1\alpha )\mid\alpha\!\in\! 2^\omega\}\!\subseteq\!
\mathbb{B}^\eta_0\cap\overline{\mathbb{C}^\eta_0\cup\mathbb{C}^\eta_1}\mbox{,}$$ 
the range of $\Delta (2^\omega )$ by $f\!\times\! f$ is a subset of $X^2$, so that $f$ takes values in $X$. We may also assume that $A_0,A_1$ are $\Ana$, and that the relation ``$(x,y)\!\in\! C_p$" is $\Borel$ in $(x,y,p)$. By Theorem \ref{kernel2}, $G_\eta$ is a nonempty $\Ana$ relation on $X$ (we denote $G_\eta\! :=\! G_{\eta ,1}$ and $F^\varepsilon_\eta\! :=\! F^\varepsilon_{\eta ,1}$, for simplicity). We also consider $F_\theta$ with $F^\varepsilon_\theta\! :=\!\overline{F_\theta}^{\tau_1}$. In the sequel, all the closures will refer to the topology $\tau_1$, so that, for example, 
$$G_\eta\cup A_0\cup A_1\!\subseteq\!
\overline{A_0\cup A_1}\!\subseteq\! R\! =\!\bigcup_{n\in\omega}~C_n.$$
$\bullet$ Let us show that 
$A_\epsilon\cap G_\eta\!\subseteq\! F^{\vert\mbox{parity}(\eta )-\epsilon\vert}_\eta$ if 
$\epsilon\!\in\! 2$. We argue by induction on 
$\eta$. If $\eta\! =\! 1$, then $A_\epsilon\cap G_1\!\subseteq\! 
A_\epsilon\cap\overline{A_{1-\epsilon}}\!\subseteq\! F^{1-\epsilon}_1$. If $\eta$ is limit, then 
$A_\epsilon\cap G_\eta\!\subseteq\! A_\epsilon\cap
\bigcap_{\theta <\eta}~F^\epsilon_\theta\!\subseteq\! F^\epsilon_\eta$. Finally, if 
$\eta\! =\!\theta\! +\! 1$, then without loss of generality suppose that $\theta$ is even, so that $\eta$ is odd and 
$$A_\epsilon\cap G_\eta\!\subseteq\! 
A_\epsilon\cap\overline{A_{1-\epsilon}\cap G_{\theta}}\!\subseteq\! 
A_\epsilon\cap F^{1-\epsilon}_\theta .$$ 
Note that this last set is contained in $F^{1-\epsilon}_\eta$, as required.\bigskip

So, if $A_\epsilon\cap G_\eta\!\neq\!\emptyset$ for some $\epsilon\!\in\! 2$ and $e$ is the correct digit, then $F^e_\eta\!\neq\!\emptyset$. Theorem \ref{checkD2Sigma01} gives 
$(\mathbb{A}_0,\mathbb{A}_1)\!\in\!\{ (\mathbb{N}^\eta_1,\mathbb{N}^\eta_0),
(\mathbb{B}^\eta_1,\mathbb{B}^\eta_0),(\mathbb{N}^\eta_0,\mathbb{N}^\eta_1),
(\mathbb{B}^\eta_0,\mathbb{B}^\eta_1)\}$ for which 
$(2^\omega ,2^\omega ,\mathbb{A}_0,\mathbb{A}_1)\!\sqsubseteq\! (X,X,A_0,A_1)$, via a square map.\bigskip

\noindent $\bullet$ Thus, in the sequel, we suppose that 
$G_\eta\cap (A_0\cup A_1)\! =\!\emptyset$. We put 
$$D_\eta\! :=\!\big\{\big( t^0_sw,t^1_sw\big)\!\in\! {\cal T}^\eta\mid s\! =\!\emptyset\big\}\! =\!
\Delta (2^{<\omega})$$ 
and, for $\theta\! <\!\eta$ and $\epsilon\!\in\! 2$,
$$D^\epsilon_\theta\! :=\!\Big\{ (t^0_sw,t^1_sw\big)\!\in\! {\cal T}^\eta\mid 
s\!\in\! T_\eta\!\setminus\!\{\emptyset\} ~\wedge ~\varphi (s)\! =\!\theta ~\wedge ~
\mbox{parity}(\vert s\vert )\! =\! 1\! -\!\big\vert\mbox{parity}\big( s(0)\big)\! -\!\epsilon\big\vert\Big\}\mbox{,}$$
so that $\{ D_\eta\}\cup\{ D^\epsilon_\theta\mid\theta\! <\!\eta ~\wedge ~\epsilon\!\in\! 2\}$ defines a partition of ${\cal T}^\eta$.\bigskip

\noindent\bf Case 1\rm\ $G_\eta\!\not\subseteq\!\Delta (X)$.\bigskip

Let $(x,y)\!\in\! G_\eta\!\setminus\!\Delta (X)$, and $O_0,O_1$ be disjoint $\Borone$ sets with 
$(x,y)\!\in\! O_0\!\times\! O_1$. We can replace $G_\eta$, $A_0$ and $A_1$ with their intersection with $O_0\!\times\! O_1$ if necessary and assume that they are contained in $O_0\!\times\! O_1$. Let us indicate the differences with the proof of Theorem \ref{checkD2Sigma01}.\bigskip

\noindent $\bullet$ Condition (6) is changed as follows:
$$(6)\ U_{s,t}\!\subseteq\!\left\{\!\!\!\!\!\!\!
\begin{array}{ll} 
& G_\eta\mbox{ if }(s,t)\!\in\! D_\eta\cr
& A_\epsilon\cap G_\theta\mbox{ if }(s,t)\!\in\! D^\epsilon_\theta
\end{array}
\right.$$
$\bullet$ If $(0\alpha ,1\beta )\!\in\!\mathbb{C}^\eta_\epsilon$, then there is $\theta\! <\!\eta$ such that $(\alpha ,\beta )\vert n\!\in\! D^\epsilon_\theta$ if $n\!\geq\! n_0$. In this case, 
$\big( U_{(\alpha ,\beta )\vert n}\big)_{n\geq n_0}$ is a decreasing sequence of nonempty clopen subsets of $A_\epsilon\cap\Omega_{X^2}$ with vanishing diameters, so that its intersection is a singleton $\{ F(\alpha ,\beta )\}\!\subseteq\! A_\epsilon$, and 
$\big( f(0\alpha ),f(1\beta )\big)\! =\! F(\alpha ,\beta )\!\in\! A_\epsilon$.\bigskip

\noindent $\bullet$ So let us prove that the construction is possible. Let 
$(x_\emptyset ,y_\emptyset )\!\in\! G_\eta\cap\Omega_{X^2}$. We choose a $\Ana$ subset 
$U_{\emptyset ,\emptyset}$ of $X^2$ such that $(x_\emptyset ,y_\emptyset )\!\in\! 
U_{\emptyset ,\emptyset}\!\subseteq\! G_\eta\cap C_{\Phi (\emptyset ,\emptyset )}\cap
\Omega_{X^2}\cap (X_\emptyset\!\times\! Y_\emptyset )$, which completes the construction for the length $l\! =\! 0$. Assume that we have constructed our objects for the sequences of length $l$. Note that $(x_{t^0_uw},y_{t^1_uw})\!\in\! G_{\varphi (u)}\cap (U\!\times\! V)\!\subseteq\! 
G_{\varphi (uq)+1}\cap (U\!\times\! V)\!\subseteq\!
\overline{A_\epsilon\cap G_{\varphi (uq)}}\cap (U\!\times\! V)$, 
where $\epsilon$ satisfies $(t^0_{uq},t^1_{uq})\!\in\! D^\epsilon_{\varphi (uq)}$. This gives 
$(x_{t^0_uw0},y_{t^1_uw1})\!\in\! 
A_\epsilon\cap G_{\varphi (uq)}\cap (U\!\times\! V)\cap\Omega_{X^2}$. If $u\! =\!\emptyset$, then 
$(t^0_uw1,t^1_uw1)\!\in\! D_\eta$, so that 
$(x_{t^0_uw1},y_{t^1_uw1})\!\in\! U_{t^0_uw,t^1_uw}\!\subseteq\! G_\eta$ and 
$(x_{t^0_uw0},y_{t^1_uw1})\!\in\! A_\epsilon$. As $G_\eta\cap (A_0\cup A_1)\! =\!\emptyset$, 
$x_{t^0_uw0}\!\not=\! x_{t^0_uw1}$. Similarly, $y_{t^0_uw0}\!\not=\! y_{t^0_uw1}$. If 
$u\!\not=\!\emptyset$, then we argue as in the proof of Theorem \ref{checkD2Sigma01} to see that 
$x_{s0}\!\not=\! x_{s1}$ (and similarly for $y_{s0}$ and $y_{s1}$).\bigskip

\noindent\bf Case 2\rm\ $G_\eta\!\subseteq\!\Delta (X)$.\bigskip

 Let us indicate the differences with the proof of Theorem \ref{checkD2Sigma01} and Case 1. We set 
$$S\! :=\!\{ x\!\in\! X\mid (x,x)\!\in\! G_\eta\}\mbox{,}$$ 
which is a nonempty $\Ana$ set by our assumption. We get $f\! :\! 2^\omega\!\rightarrow\! X$ injective continuous such that 
$\mathbb{S}^\eta_\epsilon\!\subseteq\! (f\!\times\! f)^{-1}(A_\epsilon )$ for each $\epsilon\!\in\! 2$. In this case, $A_0\cap S^2$ and $A_1\cap S^2$ are irreflexive.\bigskip

\noindent (2) $\Rightarrow$ (3) Note that 
$(2^\omega ,2^\omega ,\mathbb{N}^\eta_0,\mathbb{N}^\eta_1)\sqsubseteq 
(2^\omega ,2^\omega ,\mathbb{B}^\eta_0,\mathbb{B}^\eta_1)$ and 
$(2^\omega ,2^\omega ,\mathbb{S}^\eta_0,\mathbb{S}^\eta_1)\sqsubseteq 
(2^\omega ,2^\omega ,\mathbb{C}^\eta_0,\mathbb{C}^\eta_1)$, with witnesses 
$\alpha\!\rightarrow\! 0\alpha$ and $\beta\!\rightarrow\! 1\beta$.\bigskip
 
\noindent (3) $\Rightarrow$ (1) This comes from Lemmas \ref{S^2} and \ref{deltaeta1}.
\hfill{$\square$}

\begin{prop} \label{squareDeltaD2} Let $\eta\!\geq\! 1$ be a countable ordinal.\smallskip

(a) If $\eta$ is a successor ordinal, then the pairs $(\mathbb{N}^\eta_1,\mathbb{N}^\eta_0),(\mathbb{B}^\eta_1,\mathbb{B}^\eta_0),(\mathbb{N}^\eta_0,\mathbb{N}^\eta_1),
(\mathbb{B}^\eta_0,\mathbb{B}^\eta_1),(\mathbb{S}^\eta_0,\mathbb{S}^\eta_1)$ and 
$(\mathbb{C}^\eta_0,\mathbb{C}^\eta_1)$ are incomparable for the square reduction.\smallskip

(b) If $\eta$ is a limit ordinal, then $(2^\omega ,2^\omega ,\mathbb{S}^\eta_0,\mathbb{S}^\eta_1)
\!\sqsubseteq\! (2^\omega ,2^\omega ,\mathbb{N}^\eta_1,\mathbb{N}^\eta_0),
(2^\omega ,2^\omega ,\mathbb{N}^\eta_0,\mathbb{N}^\eta_1)$ and 
$$(2^\omega ,2^\omega ,\mathbb{C}^\eta_0,\mathbb{C}^\eta_1)
\!\sqsubseteq\! (2^\omega ,2^\omega ,\mathbb{B}^\eta_1,\mathbb{B}^\eta_0),
(2^\omega ,2^\omega ,\mathbb{B}^\eta_0,\mathbb{B}^\eta_1)\mbox{,}$$ 
via a square map, and the pairs $(\mathbb{S}^\eta_0,\mathbb{S}^\eta_1)$ and 
$(\mathbb{C}^\eta_0,\mathbb{C}^\eta_1)$ are  incomparable for the square reduction.\end{prop}

\noindent\bf Proof.\rm ~(a) We set, for $\theta\!\leq\!\eta$, 
$C_\theta\! :=\!\bigcup_{\varphi (s)\geq\theta}~\mbox{Gr}(f_s)$.\bigskip

\noindent\bf Claim.\it\ Let $\theta\!\leq\!\eta$. Then $C_\theta$ is a closed relation on $2^\omega$.\rm\bigskip

 Indeed, this is inspired by the proof of Theorem 2.3 in [L2].
 
\vfill\eject
 
  We first show that $C^l\! :=\!\bigcup_{s\in\omega^{\leq l},\varphi (s)\geq\theta}~\mbox{Gr}(f_s)$ 
is closed, by induction on $l\!\in\!\omega$. This is clear for $l\! =\! 0$. Assume that the statement is true for $l$. Note that 
$C^{l+1}\! =\! C^l\cup\bigcup_{s\in\omega^{l+1},\varphi (s)\geq\theta}~\mbox{Gr}(f_s)$. Let 
$p_m\!\in\! C^{l+1}$ such that $(p_m)_{m\in\omega}$ converges to $p$. By induction assumption, we may assume that, for each $m$, there is $(s_m,n_m)\!\in\!\omega^l\!\times\!\omega$ such that 
$\varphi (s_mn_m)\!\geq\!\theta$ and $p_m\!\in\!\mbox{Gr}(f_{s_mn_m})$. As the 
$\mbox{Gr}(f_{sn})$'s are closed, we may assume that there is $i\!\leq\! l$ such that the sequence $\big( (s_mn_m)\vert i\big)_{m\in\omega}$ is constant and the sequence 
$\big( (s_mn_m)(i)\big)_{m\in\omega}$ tends to infinity. This implies that 
$p\!\in\!\mbox{Gr}(f_{(s_0n_0)\vert i})\!\subseteq\! C^{l+1}$, which is therefore closed.\bigskip

 Now let $p_m\!\in\! C_\theta$ such that $(p_m)_{m\in\omega}$ converges to $p$. The previous fact  implies that we may assume that, for each $m$, there is $s'_m$ such that 
$\varphi (s'_m)\!\geq\!\theta$ and $p_m\!\in\!\mbox{Gr}(f_{s'_m})$, and that the sequence 
$(\vert s'_m\vert )_{m\in\omega}$ tends to infinity. Note that there is $l$ such that the set of 
$s'_m(l)$'s is infinite. Indeed, assume, towards a contradiction, that this is not the case. Then 
$\{ s\!\in\! T_\eta\mid\exists m\!\in\!\omega ~~s\!\subseteq\! s'_m\}$ is an infinite finitely branching subtree of $T_\eta$. By K\"onig's lemma, it has an infinite branch, which contradicts the wellfoundedness of $T_\eta$. So we may assume that there is $l$ such that the sequence 
$(s'_m\vert l)_{m\in\omega}$ is constant and the sequence $\big( s'_m(l)\big)_{m\in\omega}$ tends to infinity. This implies that $p\!\in\!\mbox{Gr}(f_{s'_0\vert l})\!\subseteq\! C_\theta$.
\hfill{$\diamond$}\bigskip

\noindent $\bullet$ By Lemma \ref{S^2}, $\mathbb{N}^\eta_0$ is not separable from 
$\mathbb{N}^\eta_1$ by a $\mbox{pot}\big( D_\eta (\boraone )\big)$ set, and, by Lemma 
\ref{deltaeta1}, $\mathbb{S}^\eta_0$ is not separable from $\mathbb{S}^\eta_1$ by a 
$\mbox{pot}\Big(\Delta\big( D_\eta (\boraone )\big)\Big)$ set.\bigskip

\noindent $\bullet$ Let us show that $\mathbb{N}^\eta_0$ is separable from $\mathbb{N}^\eta_1$ by a $\check D_\eta (\boraone )$ set. In fact, it is enough to see that  
$\mathbb{N}^\eta_0\!\in\!\check D_\eta (\boraone )$ if $\eta$ is odd and 
$\mathbb{N}^\eta_1\!\in\! D_\eta (\boraone )$ if $\eta$ is even. If $\eta$ is odd, then 
$$\mathbb{N}^\eta_0\! =\!\bigcup_{s\in T_\eta ,\varphi (s)\mbox{ odd}}~\mbox{Gr}(f_s)\! =\! 
C_\eta\cup\bigcup_{\theta <\eta ,\theta\mbox{ odd}}~C_\theta\!\setminus\! C_{\theta +1}.$$ 
We set, for $\theta\! <\!\eta$, $O_\theta\! :=\!\neg C_{\theta +1}$, which defines an increasing sequence of open relations on $2^\omega$ with $\mathbb{N}^\eta_0\! =\!
\neg O_{\eta -1}\cup\bigcup_{\theta <\eta ,\theta\mbox{ odd}}~O_\theta\!\setminus\! O_{\theta -1}$. Thus $\mathbb{N}^\eta_0\!\in\!\check D_\eta (\boraone )$. Similarly, if $\eta$ is even, then 
$\mathbb{N}^\eta_1\! =\!\bigcup_{s\in T_\eta ,f_\eta (s)\mbox{ odd}}~\mbox{Gr}(f_s)\! =\!
\bigcup_{\theta <\eta ,\theta\mbox{ odd}}~C_\theta\!\setminus\! C_{\theta +1}$. We set, for 
$\theta\! <\!\eta$, $O_\theta\! :=\!\neg C_{\theta +1}$, which defines an increasing sequence of open relations on $2^\omega$ with $\mathbb{N}^\eta_1\! =\!
\bigcup_{\theta <\eta ,\theta\mbox{ odd}}~O_\theta\!\setminus\! O_{\theta -1}$. Thus 
$\mathbb{N}^\eta_1\!\in\! D_\eta (\boraone )$. This shows that 
$(2^\omega ,2^\omega ,\mathbb{N}^\eta_1,\mathbb{N}^\eta_0)$ is not $\sqsubseteq$-below 
$(2^\omega ,2^\omega ,\mathbb{N}^\eta_0,\mathbb{N}^\eta_1)$, and consequently that 
$(2^\omega ,2^\omega ,\mathbb{N}^\eta_0,\mathbb{N}^\eta_1)$ is not $\sqsubseteq$-below 
$(2^\omega ,2^\omega ,\mathbb{N}^\eta_1,\mathbb{N}^\eta_0)$.\bigskip

\noindent $\bullet$ Let us show that $\mathbb{S}^\eta_\varepsilon$ is separable from 
$\mathbb{S}^\eta_{1-\varepsilon}$ by a $\check D_\eta (\boraone )$ set if $\varepsilon\!\in\! 2$. We set, for $\theta\!\leq\!\eta$, 
$$C^\varepsilon_\theta\! :=\!\bigcup_{\varphi (s)\geq\theta ,\mbox{ parity}(s(0))=\varepsilon}~
\mbox{Gr}(f_s).$$ 
As in the claim, $(C^\varepsilon_\theta )_{\theta\leq\eta}$ is a decreasing sequence of closed sets. 

\vfill\eject

 Note that 
$$\begin{array}{ll}
\mathbb{S}^\eta_\varepsilon\!\!\!\!\!
& =\!\bigcup_{s\in T_\eta\setminus\{\emptyset\} ,
\mbox{ parity}(\vert s\vert )=1-\vert\mbox{parity}(s(0))-\varepsilon\vert}~\mbox{Gr}(f_s)\cr
& =\!\bigcup_{s\in T_\eta\setminus\{\emptyset\} ,
\vert\mbox{ parity}(\varphi (s))-\mbox{parity}(\eta )\vert =1-\vert\mbox{parity}(s(0))-\varepsilon\vert}~\mbox{Gr}(f_s)\cr
& =\!\bigcup_{s\in T_\eta\setminus\{\emptyset\} ,\mbox{ parity}(s(0))=
\vert 1-\vert\vert\mbox{parity}(\varphi (s))-\mbox{parity}(\eta )\vert -\varepsilon\vert\vert}~
\mbox{Gr}(f_s)\cr
& =\!\bigcup_{\theta <\eta ,\varphi (s)=\theta}~\bigcup_{\mbox{ parity}(s(0))=
\vert 1-\vert\vert\mbox{parity}(\theta )-\mbox{parity}(\eta )\vert -\varepsilon\vert\vert}~
\mbox{Gr}(f_s)\cr
& =\!\bigcup_{\theta <\eta}~\big(\bigcup_{\varphi (s)\geq\theta ,\mbox{ parity}(s(0))=
\vert 1-\vert\vert\mbox{parity}(\theta )-\mbox{parity}(\eta )\vert -\varepsilon\vert\vert}~
\mbox{Gr}(f_s)\big)\setminus\cr
& \hfill{\big(\bigcup_{\varphi (s)\geq\theta +1,\mbox{ parity}(s(0))=
\vert 1-\vert\vert\mbox{parity}(\theta )-\mbox{parity}(\eta )\vert -\varepsilon\vert\vert}~
\mbox{Gr}(f_s)\big)}\cr
& =\!\bigcup_{\theta <\eta}~
C^{1-\vert\vert\mbox{parity}(\theta )-\mbox{parity}(\eta )\vert -\varepsilon\vert}_\theta\!\setminus\! 
C^{1-\vert\vert\mbox{parity}(\theta )-\mbox{parity}(\eta )\vert -\varepsilon\vert}_{\theta +1}.
\end{array}$$
Assume first that $\eta\! =\!\theta_0\! +\! 1$ is a successor ordinal. We define an increasing  sequence 
$(O_\theta )_{\theta <\eta}$ of open sets as follows:
$$O_\theta\! :=\!\left\{\!\!\!\!\!\!\!\!
\begin{array}{ll}
& \neg (C^{1-\varepsilon}_{\theta+1}\cup C^\varepsilon_\theta )\mbox{ if }\theta\! <\!\theta_0
\mbox{,}\cr
& \neg C^\varepsilon_\theta\mbox{ if }\theta\! =\!\theta_0\mbox{,}
\end{array}
\right.$$
so that 
$D\! :=\!\neg D\big( (O_\theta )_{\theta <\eta}\big)\!\in\!\check D_\eta (\boraone )$.\bigskip

 We now check that $D$ separates $\mathbb{S}^\eta_\varepsilon$ from 
$\mathbb{S}^\eta_{1-\varepsilon}$. If $\theta\! <\!\eta$ has a parity opposite to that of $\eta$, then either $\theta\! =\!\theta_0$ and 
$C_\theta^\varepsilon\!\setminus\! C^\varepsilon_{\theta +1}\!\subseteq\! C^\varepsilon_{\theta_0} \!\subseteq\!\neg (\bigcup_{\theta'<\eta}~O_{\theta'})\!\subseteq\! D$. Or $\theta\! <\!\theta_0$,  
$\theta\! +\! 1\! <\!\theta_0\! <\!\eta$ has the same parity as $\eta$, and 
$C_\theta^\varepsilon\!\setminus\! C^\varepsilon_{\theta +1}\!\subseteq\! O_{\theta +1}
\!\setminus\! (\bigcup_{\theta'\leq\theta}~O_{\theta'})\!\subseteq\! D$. If now $\theta\! <\!\eta$ has the same parity as $\eta$, then 
$C_\theta^{1-\varepsilon}\!\setminus\! C^{1-\varepsilon}_{\theta +1}\!\subseteq\! O_\theta
\!\setminus\! (\bigcup_{\theta'<\theta}~O_{\theta'})\!\subseteq\! D$. Thus 
$\mathbb{S}^\eta_\varepsilon\!\subseteq\! D$. Similarly, 
$\mathbb{S}^\eta_{1-\varepsilon}\!\subseteq\!\neg D$. If $\eta$ is a limit ordinal, then we set 
$O_\theta\! :=\!\neg (C^{1-\varepsilon}_{\theta+1}\cup C^\varepsilon_\theta )$ and argue similarly. This shows that 
$(2^\omega ,2^\omega ,\mathbb{N}^\eta_\varepsilon ,\mathbb{N}^\eta_{1-\varepsilon})$ is not 
$\sqsubseteq$-below $(2^\omega ,2^\omega ,\mathbb{S}^\eta_0,\mathbb{S}^\eta_1)$ for each 
$\varepsilon\!\in\! 2$.\bigskip

\noindent $\bullet$ Let us prove that 
$(2^\omega ,2^\omega ,\mathbb{S}^\eta_0,\mathbb{S}^\eta_1)$ is not $\sqsubseteq$-below 
$(2^\omega ,2^\omega ,\mathbb{N}^\eta_\varepsilon ,\mathbb{N}^\eta_{1-\varepsilon})$ if 
$\varepsilon\!\in\! 2$ and $\eta$ is a successor ordinal. Let us do it for $\varepsilon\! =\! 0$, the other case being similar. We argue by contradiction, which gives $f,g$ injective continuous with 
$\mathbb{S}^\eta_\varepsilon\!\subseteq\! (f\!\times\! g)^{-1}(\mathbb{N}^\eta_\varepsilon)$ 
for each $\varepsilon\!\in\! 2$. We set, for $\theta\! <\!\eta$ and $\varepsilon\!\in\! 2$, 
$$U^\varepsilon_\theta\! :=\!\bigcup_{\theta\leq\theta'<\eta ,\varphi (s)=\theta',
\mbox{ parity}(s(0))=
\vert 1-\vert\vert\mbox{parity}(\theta')-\mbox{parity}(\eta )\vert -\varepsilon\vert\vert}~
\mbox{Gr}(f_s).$$ 
Note that the sequence $(U^\varepsilon_\theta )_{\theta <\eta}$ is decreasing, 
$\mathbb{S}^\eta_\varepsilon\! =\! U^\varepsilon_0$,  
$$\overline{U^0_\theta\cup U^1_\theta}\! =\! C^0_\theta\cup C^1_\theta\! =\! 
U^0_\theta\cup U^1_\theta\cup\Delta (2^\omega )\! =\! C_\theta\mbox{,}$$ 
and $C^0_{\theta +1}\cup C^1_{\theta +1}\! =\!\overline{U^0_\theta}\cap\overline{U^1_\theta}$ if 
$\theta\! <\!\eta$ since 
$$\overline{U^\varepsilon_\theta}\! =\! C^0_{\theta +1}\cup C^1_{\theta +1}\cup
\bigcup_{\varphi (s)=\theta ,\mbox{ parity}(s(0))=
\vert 1-\vert\vert\mbox{parity}(\theta )-\mbox{parity}(\eta )\vert -\varepsilon\vert\vert}~
\mbox{Gr}(f_s)\mbox{,}$$ 
as in the claim. Let us prove that 
$U^0_\theta\cup U^1_\theta\!\subseteq\! (f\!\times\! g)^{-1}(C_\theta )$ if $\theta\! <\!\eta$. We argue by induction on $\theta$, and the result is clear for $\theta\! =\! 0$. If $\theta\! =\!\theta'\! +\! 1$ is a successor ordinal, then 
$$U^0_\theta\cup U^1_\theta\!\subseteq\! C^0_\theta\cup C^1_\theta\! =\!
\overline{U^0_{\theta'}}\cap\overline{U^1_{\theta'}}\!\subseteq\! (f\!\times\! g)^{-1}
(\overline{\mathbb{N}^\eta_0\cap C_{\theta'}}\cap\overline{\mathbb{N}^\eta_1\cap C_{\theta'}})
\!\subseteq\! (f\!\times\! g)^{-1}(C_\theta ).$$ 
If $\theta$ is a limit ordinal, then 
$U^0_\theta\cup U^1_\theta\!\subseteq\!\bigcap_{\theta'<\theta}~(U^0_{\theta'}\cup U^1_{\theta'})
\!\subseteq\! (f\!\times\! g)^{-1}(\bigcap_{\theta'<\theta}~C_{\theta'})
\! =\! (f\!\times\! g)^{-1}(C_\theta )$. This implies that 
$C^0_\eta\cup C^1_\eta\!\subseteq\! (f\!\times\! g)^{-1}(C_\eta )$. In particular, 
$\Delta (2^\omega )$ is sent into itself by $f\!\times\! g$ and $f\! =\! g$. As $\eta\! =\!\theta\! +\! 1$ is a successor ordinal, $U_\theta^0\!\subseteq\! (f\!\times\! f)^{-1}
(\mathbb{N}^\eta_0\cap C_\theta )\!\subseteq\! (f\!\times\! f)^{-1}\big(\Delta (2^\omega )\big)$, which contradicts the injectivity of $f$.

\vfill\eject

\noindent $\bullet$ So we proved that ${\cal A}\! :=\!\{ (\mathbb{N}^\eta_1,\mathbb{N}^\eta_0),
(\mathbb{N}^\eta_0,\mathbb{N}^\eta_1),(\mathbb{S}^\eta_0,\mathbb{S}^\eta_1)\}$ is a  
$\sqsubseteq$-antichain if $\eta$ is a successor ordinal. For the same reasons, 
${\cal B}\! :=\!\{ (\mathbb{B}^\eta_1,\mathbb{B}^\eta_0),(\mathbb{B}^\eta_0,\mathbb{B}^\eta_1),(\mathbb{C}^\eta_0,\mathbb{C}^\eta_1)\}$ is a $\sqsubseteq$-antichain if $\eta$ is a successor ordinal. Moreover, no pair in $\cal A$ is below a pair in $\cal B$ for the square reduction since 
$\Delta (2^\omega )\!\subseteq\!\mathbb{N}^\eta_0\cap
\overline{\mathbb{S}^\eta_0\cup\mathbb{S}^\eta_1}$ and the element of the pairs in $\cal B$ are contained in the clopen set $N_0\!\times\! N_1$.\bigskip

 It remains to prove that we cannot find 
$(\mathbb{A},\mathbb{B}),(\mathbb{A}',\mathbb{B}')\!\in\! {\cal A}$ and a continuous injection 
${f\! :\! 2^\omega\!\rightarrow\! 2^\omega}$ such that 
$G_\mathbb{A}\!\subseteq\! (f\!\times\! f)^{-1}(\mathbb{A}')$ and 
$G_\mathbb{B}\!\subseteq\! (f\!\times\! f)^{-1}(\mathbb{B}')$. We argue by contradiction. If 
${(\mathbb{A},\mathbb{B})\!\not=\! (\mathbb{A}',\mathbb{B}')}$ and $\varepsilon\!\in\! 2$, then we define continuous injections $f_\varepsilon\! :\! 2^\omega\!\rightarrow\! 2^\omega$ by 
$f_\varepsilon (\alpha )\! :=\! f(\varepsilon\alpha )$. Note that ${f_0\!\times\! f_1}$ reduces 
$(\mathbb{A},\mathbb{B})$ to $(\mathbb{A}',\mathbb{B}')$, which contradicts the fact that $\cal A$ is a $\sqsubseteq$-antichain. Thus ${(\mathbb{A},\mathbb{B})\! =\! (\mathbb{A}',\mathbb{B}')}$, and $(\mathbb{A},\mathbb{B})\! =\! (\mathbb{S}^\eta_0,\mathbb{S}^\eta_1)$ by Proposition 
\ref{squareD2}. As in the proof of Proposition \ref{squareD2}, $\varphi (s)\!\leq\!\varphi (v)$. If 
$\alpha\!\in\! 2^\omega$, then $(0\alpha ,1\alpha )$ is the limit of 
$(0t^0_{p_k}\gamma_k,1t^1_{p_k}\gamma_k)$. Note that 
$\big( f(0t^0_{p_k}\gamma_k),f(1t^1_{p_k}\gamma_k)\big)\! =\! 
(t^0_{v_k}\gamma'_k,t^1_{v_k}\gamma'_k)$ and $\varphi (p_k)\!\leq\!\varphi (v_k)$. As 
$\big(\varphi (p_k)\big)_{k\in\omega}$ is cofinal in $\varphi (\emptyset )\! =\!\eta$, so is 
$\big(\varphi (v_k)\big)_{k\in\omega}$. This implies that 
$\big( f(0\alpha ),f(1\alpha )\big)\!\in\!\Delta (2^\omega )$, which contradicts the injectivity of $f$.\bigskip

\noindent (b) Let us prove that 
$(2^\omega ,2^\omega ,\mathbb{S}^\eta_0,\mathbb{S}^\eta_1)\sqsubseteq 
(2^\omega ,2^\omega ,\mathbb{N}^\eta_\varepsilon ,\mathbb{N}^\eta_{1-\varepsilon})$ with a square map if $\varepsilon\!\in\! 2$. Let us do it for $\varepsilon\! =\! 0$, the other case being similar. We construct a map $\phi\! :\! 2^{<\omega}\!\rightarrow\! 2^{<\omega}$ satisfying the following:
$$\begin{array}{ll}
& (1)~\forall l\!\in\!\omega ~~\exists k_l\!\in\!\omega ~~\phi [2^l]\!\subseteq\! 2^{k_l}\cr
& (2)~\phi (s)\!\subsetneqq\!\phi (s\varepsilon )\cr
& (3)~\phi (s0)\!\not=\!\phi (s1)\cr
& (4)~\forall s\!\in\! T_\eta\!\setminus\!\{\emptyset\} ~~\Big(
\mbox{parity}(\vert s\vert )\! =\! 1\! -\!\big\vert\mbox{parity}\big( s(0)\big)\! -\!\varepsilon\big\vert
\Big)\Rightarrow\exists v_s\!\in\! T_\eta ~~\mbox{parity}(\vert v_s\vert )\! =\!\varepsilon\ \wedge\cr
& \begin{array}{ll}
& (a)~\forall w\!\in\! 2^{<\omega}~~\exists w'\!\in\! 2^{<\omega}~~
\big(\phi (t^0_sw),\phi (t^1_sw)\big)\! =\! (t^0_{v_s}w',t^1_{v_s}w')\cr
& (b)~\varphi (s)\!\leq\!\varphi (v_s)
\end{array}
\end{array}$$
Assume that this is done. Then the map 
$f\! :\!\alpha\!\mapsto\!\mbox{lim}_{n\rightarrow\infty}~\phi (\alpha\vert n)$ is as desired. So let us check that the construction of $\phi$ is possible. We construct $\phi (s)$ by induction on the length of $s$.\bigskip

- We set $k_0\! :=\! 0$ and $\phi (\emptyset )\! :=\!\emptyset$.\bigskip

- Note that $<0>_\eta\ =\! 1$ and $(t^0_0,t^1_0)\! =\! (0,1)$. As $\eta\!\geq\! 1$ is limit, 
$\varphi (1)\! >\!\varphi (0)$ are odd ordinals, so that $\varphi (10)\!\geq\!\varphi (0)$ is an even ordinal. We set $k_1\! :=<10>_\eta$, $\phi (\varepsilon )\! :=\! t^\varepsilon_{10}$ and 
$v_0\! :=\! 10$. This completes the construction of $\phi [2^1]$, and our conditions are satisfied since $k_1\! >\! 0$.\bigskip

- We next want to construct $\phi (s)$ for $s\in 2^{l+1}$, with $l\!\geq\! 1$, assuming that we have constructed $\phi (s)$ if $|s|\!\leq\! l$. Note that there is exactly one sequence $u$ such that 
$(t_u^0,t_u^1)\!\in\! 2^{l+1}$. We first define simultaneously $\phi (t_u^0)$ and $\phi (t_u^1)$, and then extend the definition to the other sequences in $2^{l+1}$.\bigskip

 If $\vert u\vert\!\geq\! 2$, then there are $u_0\!\in\!\omega^{<\omega}$ and $w\!\in\! 2^{<\omega}$ such that $t_u^\varepsilon=t^\varepsilon_{u_0}w\varepsilon$. By condition (4), 
$\big(\phi(t^0_{u_0}w),\phi(t^0_{u_0}w)\big)\! =\! (t^0_vw',t^1_vw')$ for some 
$v\!\in\!\omega^{<\omega}$ and $w'\!\in\! 2^{<\omega}$. Let $q\!\in\!\omega$ such that 
$w'\!\subseteq\!\psi (q)$ and $\varphi (u)\!\leq\!\varphi (vq)$. We can find such a $q$ because if 
$\varphi (v)\! =\!\nu\! +\! 1$, then $\varphi (vq)\! =\!\nu$, but 
$\varphi (u)\! <\!\varphi (u_0)\!\leq\!\nu\! +\! 1$ so that $\varphi (u)\!\leq\!\nu$. If $\varphi (v)$ is limit, then $\big(\varphi (vq)\big)_{q\in\omega}$ is cofinal in $\varphi (v)$ and 
$\varphi (u)\! <\!\varphi(u_0)\!\leq\!\varphi (v)$. We set 
$\phi (t^\varepsilon_{u_0}w\varepsilon )\! :=\! t^\varepsilon_{vq}$. By definition, there is 
$N\!\in\!\omega$ such that $t^\varepsilon_{vq}\! =\! t^\varepsilon_vw'0^N\varepsilon$. We set 
$\phi (s\varepsilon )\! :=\!\phi (s)0^N\varepsilon$, for any $s\in 2^l$. Conditions (1)-(3) clearly hold. So let us check condition (4). First note that 
$\big(\phi (t^0_u),\phi(t^1_u)\big)\! =\! (t^0_{vq},t^1_{vq})$ by definition, so that (4) holds for $u$ since $\vert u\vert\! -\!\vert u_0\vert\! =\! \vert vq\vert\! -\!\vert v\vert\! =\! 1$.

\vfill\eject

 Suppose now that there are $u_1\!\in\!\omega^{<\omega}$, $z\!\in\! 2^{<\omega}$ and $e\!\in\! 2$ such that $(s,t)\! =\! (t^0_{u_1}ze,t^1_{u_1}ze)$. By the induction hypothesis, 
$\big(\phi(t^0_{u_1}ze),\phi(t^1_{u_1}ze)\big)\! =\!
\big(\phi(t^0_{u_1}z)0^Ne,\phi(t^1_{u_1}z)0^Ne\big)\! =\! 
(t^0_{v_{u_1}}z'0^Ne,t^1_{v_{u_1}}z'0^Ne)$. Thus  conditions (4) is checked.\bigskip

 Otherwise, $\vert u\vert\! =\! 1$ and $u\! =<p>$ for some $p\!\in\!\omega\!\setminus\!\{ 0\}$. Let 
$w\! :=\! t^0_u\vert l$. Note there are infinitely many $q$'s such that $\phi (w)\!\subseteq\!\psi (q)$. As $\eta$ is a limit ordinal, $\big(\varphi (q)\big)_{q\in\omega}$ is strictly increasing. Thus $q$ can be chosen so that $\varphi (p)\!\leq\!\varphi (q)$. If $p$ is odd, then we set 
$\phi(t^\varepsilon_u)\! :=\! t^\varepsilon_{<q>}$. If $p$ is even, then we set 
${\phi (t^\varepsilon_u)\! :=\! t^\varepsilon_{q0}}$. Let $w^0$ and $w^1$ be the sequences such that $\phi(t^\varepsilon_u)\! =\! \phi (w)w^\varepsilon\varepsilon$. Note that they are different if $p$ is even. As in the previous case, we define 
$\phi (s\varepsilon)\! :=\!\phi (s)w^\varepsilon\varepsilon$, for any $s\in 2^l$. Notice how the choice of $w^\varepsilon$ only depends on the last coordinate of $s\varepsilon$. The conditions are verified as before for $\big(\phi(t^0_u),\phi(t^1_u)\big)$. For the other cases, 
$$\big(\phi(t^0_{u_1}ze),\phi(t^1_{u_1}ze)\big)\! =\!
 \big(\phi(t^0_{u_1}z)w^ee,\phi(t^1_{u_1}z)w^ee\big)\! =\! 
 (t^0_{v_{u_1}}w'w^ee,t^1_{v_{u_1}}w'w^ee)\mbox{,}$$ 
 by the induction hypothesis. So the conditions are checked.\bigskip
 
It remains to note that $(2^\omega ,2^\omega ,\mathbb{C}^\eta_0,\mathbb{C}^\eta_1)\sqsubseteq 
(2^\omega ,2^\omega ,\mathbb{B}^\eta_\varepsilon ,\mathbb{B}^\eta_{1-\varepsilon})$ with a square map if $\varepsilon\!\in\! 2$, with witness 
$\varepsilon\alpha\!\mapsto\!\varepsilon f(\alpha )$.\hfill{$\square$}\bigskip
 
\noindent $\underline{\mbox{\bf Consequences}}$

\begin{cor} \label{cor1delta01} Let $\eta\!\geq\! 1$ be a countable ordinal, $X$ be a Polish space, and $A,B$ be disjoint analytic relations on $X$ such that $A\cup B$ is contained in a potentially closed s-acyclic or locally countable relation. Then exactly one of the following holds:\smallskip  

(a) the set $A$ is separable from $B$ by a 
$\mbox{pot}\Big(\Delta\big( D_\eta (\boraone )\big)\Big)$ set,\smallskip  

(b) there is $(\mathbb{A}_0,\mathbb{A}_1)\!\in\!\{ (\mathbb{N}^\eta_1,\mathbb{N}^\eta_0),
(\mathbb{N}^\eta_0,\mathbb{N}^\eta_1),(\mathbb{S}^\eta_0,\mathbb{S}^\eta_1)\}$ with 
$(2^\omega ,2^\omega ,\mathbb{A}_0,\mathbb{A}_1)\!\sqsubseteq\! (X,X,A,B)$.\end{cor}

\noindent\bf Proof.\rm ~By Lemmas \ref{S^2} and \ref{deltaeta1}, (a) and (b) cannot hold simultaneously. So assume that (a) does not hold. By Lemma \ref{suffqa}, we may assume that 
$A\cup B$ is contained in a potentially closed quasi-acyclic relation. It remains to apply Theorem \ref{exdeltaeta1}.\hfill{$\square$}

\begin{cor} \label{cor2delta01} Let $\eta\!\geq\! 1$ be a countable ordinal, $X,Y$ be Polish spaces, and $A,B$ be disjoint analytic subsets of $X\!\times\! Y$ such that $A\cup B$ is contained in a potentially closed locally countable set. Then exactly one of the following holds:\smallskip  

(a) the set $A$ is separable from $B$ by a 
$\mbox{pot}\Big(\Delta\big( D_\eta (\boraone )\big)\Big)$ set,\smallskip  

(b) there is $(\mathbb{A}_0,\mathbb{A}_1)\!\in\!\{ (\mathbb{N}^\eta_1,\mathbb{N}^\eta_0),
(\mathbb{N}^\eta_0,\mathbb{N}^\eta_1),(\mathbb{S}^\eta_0,\mathbb{S}^\eta_1)\}$ with 
$(2^\omega ,2^\omega ,\mathbb{A}_0,\mathbb{A}_1)\!\sqsubseteq\! (X,X,A,B)$.\end{cor}

\noindent\bf Proof.\rm ~As in the proof of Corollary \ref{cor1delta01}, (a) and (b) cannot hold simultaneously. Then we argue as in the proof of Corollary \ref{cor2checkD2Sigma01}. $A'\cup B'$ is contained in a potentially closed locally countable relation, and $A'$ is not separable from $B'$ by a $\mbox{pot}\Big(\Delta\big( D_\eta (\boraone )\big)\Big)$ set. Corollary \ref{cor1delta01} gives $f',g'\! :\! 2^\omega\!\rightarrow\! Z$.\hfill{$\square$}

\vfill\eject

\begin{cor} \label{caracpartialdelta01} Let $\eta\!\geq\! 1$ be a countable ordinal, $X$ be a Polish space, and $A,B$ be disjoint analytic relations on $X$. The following are equivalent:\smallskip

\noindent (1) there is a potentially closed s-acyclic relation $R\!\in\!\ana$ such that $A\cap R$ is not separable from $B\cap R$ by a $\mbox{pot}\Big(\Delta\big( D_\eta (\boraone )\big)\Big)$ set,\smallskip

\noindent (2) there is a potentially closed locally countable relation $R\!\in\!\ana$ such that $A\cap R$ is not separable from $B\cap R$ by a $\mbox{pot}\Big(\Delta\big( D_\eta (\boraone )\big)\Big)$ set,\smallskip

\noindent (3) there is $(\mathbb{A}_0,\mathbb{A}_1)\!\in\!
\{ (\mathbb{N}^\eta_1,\mathbb{N}^\eta_0),(\mathbb{N}^\eta_0,\mathbb{N}^\eta_1),
(\mathbb{S}^\eta_0,\mathbb{S}^\eta_1)\}$ with 
$(2^\omega ,2^\omega ,\mathbb{A}_0,\mathbb{A}_1)\!\sqsubseteq\! (X,X,A,B)$,\smallskip

\noindent (4) there is $(\mathbb{A}_0,\mathbb{A}_1)\!\in\!
\{ (\mathbb{N}^\eta_1,\mathbb{N}^\eta_0),(\mathbb{B}^\eta_1,\mathbb{B}^\eta_0),
(\mathbb{N}^\eta_0,\mathbb{N}^\eta_1),(\mathbb{B}^\eta_0,\mathbb{B}^\eta_1),
(\mathbb{S}^\eta_0,\mathbb{S}^\eta_1),(\mathbb{C}^\eta_0,\mathbb{C}^\eta_1)\}$ such that the inequality  $(2^\omega ,2^\omega ,\mathbb{A}_0,\mathbb{A}_1)\!\sqsubseteq\! (X,X,A,B)$ holds, via a square map.\end{cor} 

\noindent\bf Proof.\rm ~(1) $\Rightarrow$ (3),(4) and (2) $\Rightarrow$ (3),(4) This is a consequence of Corollary \ref{cor1delta01} and its proof.\bigskip

\noindent (4) $\Rightarrow$ (1) By the remarks before Lemma \ref{S^2}, 
$\mathbb{N}^\eta_0\cup\mathbb{N}^\eta_1$ has s-acyclic levels. This implies that 
$\mathbb{N}^\eta_0\cup\mathbb{N}^\eta_1$ and $\mathbb{S}^\eta_0\cup\mathbb{S}^\eta_1$ are s-acyclic. As $\mathbb{N}^\eta_0\cup\mathbb{N}^\eta_1$ is antisymmetric, 
$\mathbb{B}^\eta_0\cup\mathbb{B}^\eta_1$ and $\mathbb{C}^\eta_0\cup\mathbb{C}^\eta_1$ are s-acyclic too, by Lemma \ref{suffacy}. Thus we can take 
$R\! :=\! (f\!\times\! f)[\mathbb{A}_0\cup\mathbb{A}_1]$ since the s-acyclicity is preserved by images by the square of an injection, and by Lemmas \ref{S^2} and \ref{deltaeta1}.\bigskip

\noindent (3),(4) $\Rightarrow$ (2) We can take 
$R\! :=\! (f\!\times\! f)[\mathbb{A}_0\cup\mathbb{A}_1]$ since $\mathbb{A}_0\cup\mathbb{A}_1$ is locally countable, by Lemmas \ref{S^2} and \ref{deltaeta1}.\hfill{$\square$}

\section{$\!\!\!\!\!\!$ Background}\indent

 We now give some material to prepare the study of the Borel classes of rank two.\bigskip

\noindent $\underline{\mbox{\bf Potential Wadge classes}}$\bigskip

 In Theorem \ref{motivating}, $\mathbb{S}_0\cup\mathbb{S}_1$ is a subset of the body of a tree 
$T$ on $2^2$ which does not depend on $\bf\Gamma$. We first describe a simple version of $T$, which is sufficient to study the Borel classes (see [L6]). We identify $(2^l)^2$ and $(2^2)^l$, for each 
$l\!\in\!\omega\! +\! 1$.
 
\begin{defi} \label{frame} (1) Let 
${\cal F}\!\subseteq\!\bigcup_{l\in\omega}~(2^l)^2\!\equiv\! (2^2)^{<\omega}$. We say that 
${\cal F}$ is a \bf frame\it\ if\medskip

(a) $\forall l\!\in\!\omega~\exists ! (s_l,t_l)\!\in\! {\cal F}\!\cap\! (2^l)^2$,\smallskip

(b) $\forall p,q\!\in\!\omega~\forall w\!\in\! 2^{<\omega}~\exists N\!\in\!\omega~
(s_q0w0^N,t_q1w0^N)\!\in\! {\cal F}$ and $(|s_q0w0^N|\! -\! 1)_0\! =\! p$,\smallskip

(c) $\forall l\! >\! 0~\exists q\! <\! l~\exists w\!\in\! 2^{<\omega}\ (s_l,t_l)\! =\! (s_q0w,t_q1w)$.\medskip

\noindent (2) If ${\cal F}\! =\!\{ (s_l,t_l)\mid l\!\in\!\omega\}$ is a frame, then we will call $T$ the tree on $2^2$ generated by $\cal F$: 
$$T\! :=\!\big\{ (s,t)\!\in\! (2^2)^{<\omega}\mid s\! =\!\emptyset\vee
\big(\exists q\!\in\!\omega ~\exists w\!\in\! 2^{<\omega}~(s,t)\! =\! (s_q0w,t_q1w)\big)\big\}.$$
\end{defi}

 The existence condition in (a) and the density condition (b) ensure that $\lceil T\rceil$ is big enough to contain sets of arbitrary high potential complexity. The uniqueness condition in (a) and condition (c) ensure that $\lceil T\rceil$ is small enough to make the reduction in Theorem 
\ref{motivating} possible. The last part of condition (b) gives a control on the verticals which is very useful to construct complicated examples.

\vfill\eject

 In the sequel, $T$ will be the tree generated by a fixed frame $\cal F$ (Lemma 3.3 in [L6] ensures the existence of concrete frames). Note that $\lceil T\rceil\!\subseteq\! N_0\!\times\! N_1$, which will be useful in the sequel (recall that $N_s$ is the basic clopen set of sequences beginning with $s\!\in\! 2^{<\omega}$).\bigskip
 
\noindent $\underline{\mbox{\bf Acyclicity}}$\bigskip

 We will use some material from [L6] and [L8], where some possibly different notions of acyclicity of the levels of $T$ are involved. We will check that they coincide in our case.

\begin{defi} Let $X$ be a set, and $A$ be a relation on $X$.\smallskip

(a) An \bf $A$-path\it\ is a finite sequence $(x_i)_{i\leq n}$ of points of $X$ such that 
$(x_i,x_{i+1})\!\in\! A$ if $i\! <\! n$.\smallskip

(b) We say that $A$ is \bf connected\it\ if for any $x,y\!\in\! X$ there is an $A$-path $(x_i)_{i\leq n}$ with $x_0\! =\! x$ and $x_n\! =\! y$.\smallskip

(c) An \bf $A$-cycle\it\ is an $A$-path $(x_i)_{i\leq L}$ with $L\!\geq\! 3$, $(x_i)_{i<L}$ is injective and $x_L\! =\! x_0$ (so that $A$ is acyclic if and only if there is no $A$-cycle).\end{defi}

\begin{lem} \label{tree} Let $l\!\in\!\omega$, and $T_l\! :=\! T\cap (2^l)^2$ be the $l$th 
\bf  level\it\ of $T$.\smallskip

(a) $s(T_l)$ is connected and acyclic. In particular, $\lceil T\rceil$ is s-acyclic.\smallskip

(b) A tree $S$ on $2^2$ has acyclic levels in the sense of [L6] if and only if $S$ has suitable levels in the sense of [L8], and this is the case of $T$.\end{lem}

\noindent\bf Proof.\rm ~(a) We argue by induction on $l$. The statement is clear for $l\! =\! 0$. For the inductive step we use the fact that $T_{l+1}\! =\!\{ (s\varepsilon ,t\varepsilon )\mid (s,t)\!\in\! T_l\wedge\varepsilon\!\in\! 2\}\cup\{ (s_l0,t_l1)\}$. As the map $s\varepsilon\!\mapsto\! s$ defines an isomorphism from $\{ (s\varepsilon ,t\varepsilon )\mid (s,t)\!\in\! T_l\}$ onto $T_l$, we are done. A cycle for $s(\lceil T\rceil )$ gives a cycle for $s(T_l)$, for $l$ big enough to ensure the injectivity of the initial segments.\bigskip

\noindent (b) Assume that $S$ has acyclic levels in the sense of [L6]. This means that, for each 
$l$, the graph $G_{S_l}$ with set of vertices $2^l\!\oplus\! 2^l$ (with typical element 
$\overline{x_\varepsilon}\! :=\! (x_\varepsilon ,\varepsilon )\!\in\! 2^l\!\times\! 2$) and set of edges 
$$\big\{\{\overline{x_0},\overline{x_1}\}\mid\vec x\! :=\! (x_0,x_1)\!\in\! S_l\big\}$$
is acyclic. We have to see that $S$ has suitable levels in the sense of [L8]. This means that, for each $l$, the following hold:\smallskip

- $S_l$ is finite,\smallskip

- $\exists\varepsilon\!\in\! 2~~x^0_\varepsilon\!\not=\! x^1_\varepsilon$ if 
$\vec {x^0}\!\not=\!\vec {x^1}\!\in\! S_l$,\smallskip

- consider the graph $G^{S_l}$ with set of vertices $S_l$ and set of edges 
$$\big\{\{\vec {x^0},\vec {x^1}\}\mid\vec {x^0}\!\not=\!\vec {x^1}\wedge\exists\varepsilon\!\in\! 2~~
x^0_\varepsilon\! =\! x^1_\varepsilon\big\} ;$$ 
then for any $G^{S_l}$-cycle $(\vec {x^n})_{n\leq L}$, there are $\varepsilon\!\in\! 2$ and 
$k\! <\! m\! <\! n\! <\! L$ such that $x^k_\varepsilon\! =\! x^m_\varepsilon\! =\! x^n_\varepsilon$.\bigskip

 The first two properties are obvious. So assume that $(\vec {x^n})_{n\leq L}$ is a $G^{S_l}$-cycle for which we cannot find $\varepsilon\!\in\! 2$ and $k\! <\! m\! <\! n\! <\! L$ such that 
$x^k_\varepsilon\! =\! x^m_\varepsilon\! =\! x^n_\varepsilon$.

\vfill\eject

\noindent\bf Case 1\rm\ $x^0_0\! =\! x^1_0$.\bigskip

\noindent\bf Subcase 1.1\rm\ $L$ is odd.\bigskip

 Note that $L\!\geq\! 5$. Indeed, $L\!\geq\! 3$ since $(\vec {x^n})_{n\leq L}$ is a $G^{S_l}$-cycle. So we just have to see that $L\!\not=\! 3$. As $x^0_0\! =\! x^1_0$ and 
$\vec {x^0}\!\not=\!\vec {x^1}$, $x^0_1\!\not=\! x^1_1$. By the choice of $(\vec {x^n})_{n\leq L}$, $x^1_0\!\not=\! x^2_0$. Thus $x^1_1\! =\! x^2_1$. By the choice of $(\vec {x^n})_{n\leq L}$, 
$x^2_1\!\not=\! x^3_1$. Thus $x^2_0\! =\! x^3_0$ and $x^3_0\!\not=\! x^0_0$. Therefore 
$\vec {x^3}\!\not=\!\vec {x^0}$ and $L\!\not=\! 3$.\bigskip

 Then 
$\overline{x^0_0},\overline{x^1_1},\overline{x^2_0},...,\overline{x^{L-2}_1},\overline{x^{L-1}_0}$ is a $G_{S_l}$-cycle, by the choice of $(\vec {x^n})_{n\leq L}$.\bigskip

\noindent\bf Subcase 1.2\rm\ $L$ is even, in which case $L\!\geq\! 4$.\bigskip

 Then 
$\overline{x^0_0},\overline{x^1_1},\overline{x^2_0},...,\overline{x^{L-1}_1},\overline{x^{L}_0}$ is a $G_{S_l}$-cycle, by the choice of $(\vec {x^n})_{n\leq L}$.\bigskip
 
\noindent\bf Case 2\rm\ $x^0_0\!\not=\! x^1_0$.\bigskip

 The same arguments work, we just have to exchange the indexes.\bigskip
 
\noindent $\bullet$ Conversely, assume that $(\overline{x^n_{\varepsilon_n}})_{n\leq L}$ is a 
$G_{S_l}$-cycle. Then $L$ is even, and actually $L\!\geq\! 4$.\bigskip

\noindent\bf Case 1\rm\ $\varepsilon_0\! =\! 0$.\bigskip

 Then $(x^0_{\varepsilon_0},x^1_{\varepsilon_1}),(x^2_{\varepsilon_2},x^1_{\varepsilon_1}),...,
(x^{L-2}_{\varepsilon_{L-2}},x^{L-1}_{\varepsilon_{L-1}}),
(x^{L}_{\varepsilon_{L}},x^{L-1}_{\varepsilon_{L-1}}),(x^0_{\varepsilon_0},x^1_{\varepsilon_1})$ is a $G^{S_l}$-cycle of length $L\! +\! 1$. If $\varepsilon\!\in\! 2$, then each $\varepsilon$th coordinate appears exactly twice  before the last element of the cycle.\bigskip

\noindent\bf Case 2\rm\ $\varepsilon_0\! =\! 1$.\bigskip

 The same argument works, we just have to exchange the coordinates.\bigskip

\noindent $\bullet$ By Proposition 3.2 in [L6], $T$ has acyclic levels in the sense of [L6].
\hfill{$\square$}

\section{$\!\!\!\!\!\!$ The classes $\bormtwo$ and $\boratwo$}

$\underline{\mbox{\bf Example}}$\bigskip

 We will use an example for ${\bf\Gamma}\! =\!\bormtwo$ different from that in [L6], so that we prove the following.

\begin{lem} \label{E0} $\lceil T\rceil\cap\mathbb{E}_0$ is not separable from 
$\lceil T\rceil\!\setminus\!\mathbb{E}_0$ by a $\mbox{pot}(\bormtwo )$ set.\end{lem}

\noindent\bf Proof.\rm ~We argue by contradiction, which gives $P\!\in\!\mbox{pot}(\bormtwo )$, and also a dense $G_\delta$ subset $G$ of $2^\omega$ such that 
$P\cap G^2\!\in\! \bormtwo (G^2)$. Let $(O_n)_{n\in\omega}$ be a sequence of dense open subsets of $2^\omega$ with intersection $G$. Note that 
$\lceil T\rceil\cap\mathbb{E}_0\cap G^2\! =\!\lceil T\rceil\cap P\cap G^2\!\in\!
\bortwo (\lceil T\rceil\cap G^2)$. By Baire's theorem, it is enough to prove that 
$\lceil T\rceil\cap\mathbb{E}_0\cap G^2$ is dense and co-dense in the nonempty space 
$\lceil T\rceil\cap G^2$. So let $q\!\in\!\omega$ and $w\!\in\! 2^{<\omega}$. Pick 
$u_0\!\in\! 2^\omega$ such that $N_{s_q0wu_0}\!\subseteq\! O_0$, 
$v_0\!\in\! 2^\omega$ such that $N_{t_q1wu_0v_0}\!\subseteq\! O_0$, 
$u_1\!\in\! 2^\omega$ such that $N_{s_q0wu_0v_0u_1}\!\subseteq\! O_1$,  
$v_1\!\in\! 2^\omega$ such that $N_{t_q1wu_0v_0u_1v_1}\!\subseteq\! O_1$, and so on.

\vfill\eject

 Then 
$(s_q0wu_0v_0u_1v_1...,t_q1wu_0v_0u_1v_1...)\!\in\!\lceil T\rceil\cap\mathbb{E}_0\cap G^2$. 
Similarly, pick $N_0\!\in\!\omega$ such that $(s_q0w0^{N_0},t_q1w0^{N_0})\!\in\! {\cal F}$, 
$u_0\!\in\! 2^\omega$ such that $N_{s_q0w0^{N_0}0u_0}\!\subseteq\! O_0$, 
$v_0\!\in\! 2^\omega$ such that $N_{t_q1w0^{N_0}1u_0v_0}\!\subseteq\! O_0$, 
$N_1\!\in\!\omega$ such that 
$(s_q0w0^{N_0}0u_0v_00^{N_1},t_q1w0^{N_0}1u_0v_00^{N_1})\!\in\! {\cal F}$, 
$u_1\!\in\! 2^\omega$ such that 
$$N_{s_q0w0^{N_0}0u_0v_00^{N_1}0u_1}\!\subseteq\! O_1\mbox{,}$$ 
$v_1\!\in\! 2^\omega$ such that $N_{t_q1w0^{N_0}1u_0v_00^{N_1}1u_1v_1}\!\subseteq\! O_1$, and so on. Then 
$$(s_q0w0^{N_0}0u_0v_00^{N_1}0u_1v_1...,
t_q1w0^{N_0}1u_0v_00^{N_1}1u_1v_1...)\!\in\!\lceil T\rceil\cap G^2\!\setminus\!\mathbb{E}_0.$$
This finishes the proof.\hfill{$\square$}\bigskip

\noindent $\underline{\mbox{\bf The main result}}$\bigskip

 We reduce the study of disjoint analytic sets to that of disjoint Borel sets of low complexity, for the first classes we are considering.

\begin{lem} \label{anaKsigma} Let $X$ be a Polish space, and $A,B$ be disjoint analytic relations on $X$. Then exactly one of the following holds:\smallskip  

(a) the set $A$ is separable from $B$ by a $\mbox{pot}(\bormtwo )$ set,\smallskip  

(b) there is a $K_\sigma$ relation $A'\!\subseteq\! A$ which is not $\mbox{pot}(\bormtwo )$ such that $\overline{A'}\!\setminus\! A'\!\subseteq\! B$.\end{lem}

\noindent\bf Proof.\rm ~Theorem 1.10 in [L8] and Lemmas \ref{E0}, \ref{tree} give 
$g,h\! :\! 2^\omega\!\rightarrow\! X$ continuous such that the inclusions 
$\lceil T\rceil\cap\mathbb{E}_0\!\subseteq\! (g\!\times\! h)^{-1}(A)$ and 
$\lceil T\rceil\!\setminus\!\mathbb{E}_0\!\subseteq\! (g\!\times\! h)^{-1}(B)$ hold. We set 
$A'\! :=\! (g\!\times\! h)\big[\lceil T\rceil\cap\mathbb{E}_0\big]$, 
$B'\! :=\! (g\!\times\! h)\big[\lceil T\rceil\!\setminus\!\mathbb{E}_0\big]$ and 
$C'\! :=\! (g\!\times\! h)\big[\lceil T\rceil\big]$. Note that $A'$ is a $K_\sigma$ subset of $A$, 
$B'\!\subseteq\! B$, so that the compact set $C'$ is the disjoint union of $A'$ and $B'$. As 
$\lceil T\rceil\cap\mathbb{E}_0$ is dense in $\lceil T\rceil$, $C'$ is also the closure of $A'$. As 
$\lceil T\rceil\cap\mathbb{E}_0\! =\!\lceil T\rceil\cap (g\!\times\! h)^{-1}(A')$, $A'$ is not 
$\mbox{pot}(\bormtwo )$, by Lemma \ref{E0}.\hfill{$\square$}

\begin{thm} \label{delta02} Let $X$ be a Polish space, and $A,B$ be disjoint analytic relations on $X$ such that $A$ is quasi-acyclic. Then one of the following holds:\smallskip  

(a) the set $A$ is separable from $B$ by a $\mbox{pot}(\bormtwo )$ set,\smallskip  

(b) there is $f\! :\! 2^\omega\!\rightarrow\! X$ injective continuous such that the inclusions 
$\lceil T\rceil\cap\mathbb{E}_0\!\subseteq\! (f\!\times\! f)^{-1}(A)$ and 
$\lceil T\rceil\!\setminus\!\mathbb{E}_0\!\subseteq\! (f\!\times\! f)^{-1}(B)$ hold.\end{thm}

\noindent\bf Proof.\rm ~Assume that (a) does not hold. By Lemma \ref{anaKsigma}, we may assume that $B$ is the complement of $A$. Let $(C_n)_{n\in\omega}$ be a witness for the fact that $A$ is quasi-acyclic. Note that there are disjoint Borel subsets $O_0,O_1$ of $X$ such that 
$A\cap (O_0\!\times\! O_1)$ is not $\mbox{pot}(\bormtwo )$. We may assume that $X$ is zero-dimensional, the $C_n$'s are closed, and $O_0,O_1$ are clopen, refining the topology if necessary. We can also replace $A$ and the $C_n$'s with their intersection with 
$O_0\!\times\! O_1$ and assume that they are contained in $O_0\!\times\! O_1$.\bigskip

\noindent $\bullet$ We may assume that $X$ is recursively presented, $O_0,O_1\!\in\!\Borel$ and the relation ``$(x,y)\!\in\! C_n$" is $\Borel$ in $(x,y,n)$. As ${\it\Delta}_X$ is Polish finer than the topology on $X$, $A\!\notin\!\bormtwo (X^2,\tau_1)$. We now perform the following derative on 
$A$. We set, for $F\!\in\!\bormone (X^2,\tau_1)$, 
$F'\! :=\!\overline{F\cap A}^{\tau_1}\cap\overline{F\!\setminus\! A}^{\tau_1}$ (see 22.30 in [K]).

\vfill\eject

 Then we inductively define, for any ordinal $\xi$, $F_\xi$ by
$$\left\{\!\!\!\!\!\!
\begin{array}{ll}
& F_0\! :=\! X^2\cr
& F_{\xi +1}\! :=\! F'_\xi\cr
& F_\lambda\! :=\!\bigcap_{\xi <\lambda} F_\xi\mbox{ if }\lambda\mbox{ is limit}
\end{array}\right.$$
(see 22.27 in [K]). As $(F_\xi )$ is a decreasing sequence of closed subsets of the Polish space 
$(X^2,\tau_1)$, there is $\theta\! <\!\omega_1$ such that $F_\theta\! =\! F_{\theta +1}$. In particular, 
$F_\theta\! =\! F_{\theta +1}\! =\! F'_\theta\! =\!
\overline{F_\theta\cap A}^{\tau_1}\cap\overline{F_\theta\!\setminus\! A}^{\tau_1}$, so that 
$F_\theta\cap A$ and $F_\theta\!\setminus\! A$ are $\tau_1$-dense in $F_\theta$.\bigskip

\noindent $\bullet$ Let us prove that $F_\theta$ is not empty. We argue by contradiction:
$$X^2\! =\!\neg F_\theta\! =\!\bigcup_{\xi\leq\theta}~\neg F_\xi\! =\!
\bigcup_{\xi\leq\theta}~(\neg F_\xi\cap\bigcap_{\eta <\xi}~F_\eta )\! =\!
\bigcup_{\xi <\theta}~F_\xi\!\setminus\! F_{\xi +1}\mbox{,}$$
so that $A\! =\!\bigcup_{\xi <\theta}~A\cap F_\xi\!\setminus\! F_{\xi +1}$. But 
$A\cap F_\xi\!\setminus\! F_{\xi +1}\! =\! A\cap F_\xi\!\setminus\! 
(\overline{F_\xi\cap A}^{\tau_1}\cap\overline{F_\xi\!\setminus\! A}^{\tau_1})\! =\! 
F_\xi\!\setminus\!\overline{F_\xi\!\setminus\! A}^{\tau_1}$. This means that 
$(F_\xi\!\setminus\! F_{\xi +1})_{\xi <\theta}$ is a countable partition of $(X^2,\tau_1)$ into 
$\bortwo$ sets, and that $A$ is $\bortwo$ on each piece of the partition. This implies that $A$ is 
$\bortwo (X^2,\tau_1)$, which is absurd.\bigskip

\noindent $\bullet$ Let us prove that $F_\theta$ is $\Ana$. We use 7C in [Mo]. We define a set relation by 
$$\varphi (x,y,P)\Leftrightarrow (x,y)\!\notin\! (\neg P)'.$$
Note that $\varphi$ is monotone, and thus operative. It is also $\Ca$ on $\Ca$. By 3E.2, 3F.6 and 4B.2 in [Mo], we can apply 7C.8 in [Mo], so that $\varphi^\infty (x,y)$ is $\Ca$. An induction shows that $\varphi^\xi (x,y)$ is equivalent to ``$(x,y)\!\notin\! F_{\xi +1}$". Thus $(x,y)\!\notin\! F_\theta$ is equivalent to $(x,y)\!\notin\!\bigcap_\xi~F_\xi\! =\!\bigcap_\xi~F_{\xi +1}$, 
$(x,y)\!\in\!\bigcup_\xi\neg F_{\xi +1}$ and $\varphi^\infty (x,y)$.\bigskip

\noindent $\bullet$ We are ready to prove the following key property:
$$\forall q\!\in\!\omega ~~\forall U,V\!\in\!\Ana (X)~~F_\theta\cap (U\!\times\! V)\!\not=\!\emptyset\Rightarrow\exists n\!\geq\! q~~F_\theta\cap C_n\cap (U\!\times\! V)\!\not=\!\emptyset .$$

 Indeed, this property says that $I\! :=\! F_\theta\cap (\bigcup_{n\geq q}~C_n)$ is 
${\it\Sigma}_X^2$-dense in $F_\theta$ for each $q\!\in\!\omega$. We fix $q\!\in\!\omega$, and prove first that $I$ is $\tau_1$-dense in $F_\theta$. So let $U,V\!\in\!\Borel$ such that 
$F_\theta\cap (U\!\times\! V)$ is nonempty. As $F_\theta\!\setminus\! A$ is $\tau_1$-dense in 
$F_\theta$, we get $(x,y)\!\in\! (F_\theta\!\setminus\! A)\cap (U\!\times\! V)$. As 
$F_\theta\cap A$ is $\tau_1$-dense in $F_\theta$, we get $(x_k,y_k)\!\in\! F_\theta\cap A$ converving to $(x,y)$ for $\tau_1$. Pick $n_k\!\in\!\omega$ such that $(x_k,y_k)\!\in\! C_{n_k}$. As $C_{n_k}$ is closed, and thus $\tau_1$-closed, we may assume that the sequence 
$(n_k)_{k\in\omega}$ is strictly increasing. Now $(x_k,y_k)\!\in\! I\cap (U\!\times\! V)$ if $k$ is big enough. In order to get the statement for ${\it\Sigma}_X^2$, we have to note that $I$ is $\Ana$ since $F_\theta$ is $\Ana$ and the relation ``$(x,y)\!\in\! C_n$" is $\Borel$ in $(x,y,n)$. This implies that $\overline{I}^{\tau_1}\! =\!\overline{I}^{{\it\Sigma}_X^2}$, by a double application of the separation theorem. Therefore 
$F_\theta\!\subseteq\!\overline{I}^{\tau_1}\! =\!\overline{I}^{{\it\Sigma}_X^2}$ and $I$ is 
${\it\Sigma}_X^2$-dense in $F_\theta$.\bigskip

\noindent $\bullet$ We set, for $\vec u\! =\! (u_0,u_1)\!\in\! T\!\setminus\!\{\vec\emptyset\}$,
$$\begin{array}{ll}
n(\vec u)\!\!\!\! 
& := \mbox{Card}\big(\{i\! <\!\vert\vec u\vert\mid u_0(i)\!\not=\! u_1(i)\}\big)\mbox{,}\cr
\vec t(\vec u)\!\!\!\!
& :=\! (s_q0,t_q1)\ \mbox{ if }\ \vec u\! =\! (s_q0w,t_q1w).
\end{array}$$
$\bullet$ We are ready for the construction of $f$. We construct the following objects:\bigskip

- sequences $(x_s)_{s\in 2^{<\omega}\setminus\{\emptyset\} , s(0)=0}$, 
$(y_s)_{s\in 2^{<\omega}\setminus\{\emptyset\} , s(0)=1}$ of points of $X$,\smallskip

- sequences $(X_s)_{s\in 2^{<\omega}\setminus\{\emptyset\} , s(0)=0}$, 
$(Y_s)_{s\in 2^{<\omega}\setminus\{\emptyset\} , s(0)=1}$ of $\Ana$ subsets of $X$,\smallskip

- a map 
$\Phi\! :\!\big\{\vec t(\vec u)\mid\vec u\!\in\! T\!\setminus\!\{\vec\emptyset\}\big\}\!\rightarrow\!\omega$.\bigskip

\noindent We want these objects to satisfy the following conditions:
$$\begin{array}{ll}
& (1)\ x_s\!\in\! X_s\ \wedge\ y_s\!\in\! Y_s\cr
& (2)\ X_{s\varepsilon}\!\subseteq\! X_s\!\subseteq\!\Omega_X\cap O_0\ \wedge\ 
Y_{s\varepsilon}\!\subseteq\! Y_s\!\subseteq\!\Omega_X\cap O_1\cr
& (3)\ \mbox{diam}_{\mbox{GH}}(X_s)\mbox{, }\mbox{diam}_{\mbox{GH}}(Y_s)\!\leq\! 
2^{-\vert s\vert}\cr
& (4)\ (x_{u_0},y_{u_1})\!\in\! F_\theta\cap C_{\Phi(\vec t(\vec u))}\cr
& (5)\ (X_{u_0}\!\times\! Y_{u_1})\cap (\bigcup_{n<n(\vec u)}\ C_n)\! =\!\emptyset\cr
& (6)\ X_{s0}\cap X_{s1}\! =\! Y_{s0}\cap Y_{s1}\! =\!\emptyset
\end{array}$$
$\bullet$ Assume that this has been done. As in the proof of Theorem \ref{checkD2Sigma01}, we get $f\! :\! N_\varepsilon\!\rightarrow\! O_\varepsilon$ injective continuous, so that 
$f\! :\! 2^\omega\!\rightarrow\! X$ is injective continuous. If 
$(\alpha ,\beta )\!\in\!\lceil T\rceil\cap\mathbb{E}_0$, then 
$\Phi (\vec t\big( (\alpha ,\beta )\vert n\big) )\! =\! N$ if $n$ is big enough. In this case, by (4), 
$(x_{\alpha\vert n},y_{\beta\vert n})\!\in\! C_N$ which is closed, so that 
${\big( f(\alpha ),g(\beta )\big)\!\in\! C_N\!\subseteq\! A}$. If 
$(\alpha ,\beta )\!\in\!\lceil T\rceil\!\setminus\!\mathbb{E}_0$, then the sequence 
$(n\big( (\alpha ,\beta )\vert n\big) )_{n>0}$ tends to infinity. Thus 
$\big( f(\alpha ),g(\beta )\big)$ is not in $\bigcup_{n\in\omega}~C_n\! =\! A$ by (5).\bigskip

\noindent $\bullet$ So let us prove that the construction is possible. The key property gives 
$\Phi (0,1)\!\geq\! 1$ and $(x_0,y_1)$ in $F_\theta\cap C_{\Phi (0,1)}\cap\Omega_{X^2}$. As 
$\Omega_{X^2}\!\subseteq\!\Omega_X^2$, $x_0,y_1\!\in\!\Omega_X$. We choose $\Ana$ subsets $X_0,Y_1$ of $X$ with GH-diameter at most $2^{-1}$ such that 
$(x_0,y_1)\!\in\! X_0\!\times\! Y_1\!\subseteq\!
\big( (\Omega_X\cap O_0)\!\times\! (\Omega_X\cap O_1)\big)\!\setminus\! C_0$, which completes the construction for the length $l\! =\! 1$.\bigskip

 Let $l\!\geq\! 1$. We now want to build $x_s,X_s,y_s,Y_s$ for $s\!\in\! 2^{l+1}$, as well as 
$\Phi (s_l0,t_l1)$. Note that 
$(x_{s_l},y_{t_l})\!\in\! F_\theta\cap (U\!\times\! V)$, where 
$$\begin{array}{ll}
& U\! :=\!\{ x'_{s_l}\!\in\! X_{s_l}\mid
\exists (x'_s)_{s\in 2^l\setminus\{ s_l\} ,s(0)=0}\!\in\!\Pi_{s\in 2^l\setminus\{ s_l\} ,s(0)=0}~X_s~~
\exists (y'_s)_{s\in 2^l,s(0)=1}\!\in\!\Pi_{s\in 2^l,s(0)=1}~Y_s\cr
& \hfill{\forall\vec u\!\in\! T\cap (2^l\!\times\! 2^l)~~
(x'_{u_0},y'_{u_1})\!\in\! F_\theta\cap C_{\Phi (\vec t(\vec u))}\}\mbox{,}}\cr
& V\! :=\!\{ y'_{t_l}\!\in\! Y_{t_l}\mid
\exists (x'_s)_{s\in 2^l,s(0)=0}\!\in\!\Pi_{s\in 2^l,s(0)=0}~X_s~~
\exists (y'_s)_{s\in 2^l\setminus\{ t_l\} ,s(0)=1}\!\in\!\Pi_{s\in 2^l\setminus\{ t_l\} ,s(0)=1}~Y_s\cr
& \hfill{\forall\vec u\!\in\! T\cap (2^l\!\times\! 2^l)~~
(x'_{u_0},y'_{u_1})\!\in\! F_\theta\cap C_{\Phi (\vec t(\vec u))}\} .}
\end{array}$$
The key property gives 
$\Phi (s_l0,t_l1)\! >\!\mbox{max}\big( n(s_l0,t_l1),\mbox{max}_{q<l}~\Phi (s_q0,t_q1)\big)$ and 
$$(x_{s_l0},y_{t_l1})\!\in\! F_\theta\cap C_{\Phi (s_l0,t_l1)}\cap (U\!\times\! V).$$ 
The fact that $x_{s_l0}\!\in\! U$ gives witnesses $(x_{s0})_{s\in 2^l\setminus\{ s_l\} ,s(0)=0}$ and 
$(y_{s0})_{s\in 2^l,s(0)=1}$. Similarly, the fact that 
$y_{t_l1}\!\in\! V$ gives $(x_{s1})_{s\in 2^l,s(0)=0}$ and 
$(y_{s1})_{s\in 2^l\setminus\{ t_l\} ,s(0)=1}$. Note that $x_{s_l0}\!\not=\! x_{s_l1}$ because 
$$(x_{s_l0},y_{t_l1})\!\in\! C_{\Phi (s_l0,t_l1)}\mbox{,}$$ 
$(x_{s_l1},y_{t_l1})\!\in\! C_{\Phi (\vec t (s_l1,t_l1))}$, and 
$\Phi (s_l0,t_l1)\! >\!\Phi\big(\vec t (s_l1,t_l1)\big)$. Similarly, $y_{t_l0}\!\not=\! y_{t_l1}$. If 
$s\!\in\! 2^l$, then the connectedness of $s(T_l)$ gives an injective $s(T)$-path $p_s$ from $s$ to $s_l$. This gives a $s(A)$-path from $x_{s0}$ to $x_{s1}$ if $s(0)\! =\! 0$, and a $s(A)$-path from $y_{s0}$ to $y_{s1}$ if $s(0)\! =\! 1$. Using the quasi-acyclicity of $A$, we see, by induction on the length of $p_s$, that $x_{s0}\!\not=\! x_{s1}$ and $y_{s0}\!\not=\! y_{s1}$.

\vfill\eject

 The following picture illustrates the situation when $l\! =\! 2$.
$$\xymatrixrowsep{0.3in}
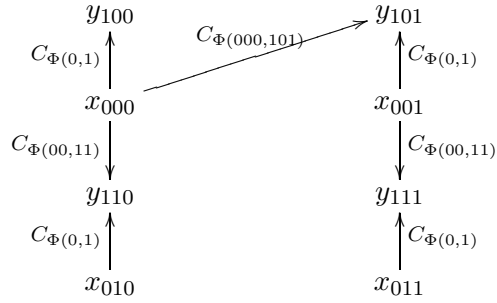
\xymatrix{ 
                       & y_{100} & & & y_{101} \\
                       & x_{000} \ar[u]^{C_{\Phi (0,1)}} \ar[d]_{C_{\Phi (00,11)}} \ar[urrr]^[@]{\hbox to 0pt{\hss $\scriptstyle{C_{\Phi (000,101)}}$\hss}} & & & x_{001} 
                       \ar[u]_{C_{\Phi (0,1)}} \ar[d]^{C_{\Phi (00,11)}} \\ 
                       & y_{110} & & & y_{111} \\  
                       & x_{010} \ar[u]^{C_{\Phi (0,1)}} & & & x_{011} \ar[u]_{C_{\Phi (0,1)}} \\     }$$ 
Then we take small enough $\Ana$ neighborhoods of the $x_{s\varepsilon}$'s and $y_{s\varepsilon}$'s to complete the construction.\hfill{$\square$}\bigskip

\noindent $\underline{\mbox{\bf Consequences}}$

\begin{cor} \label{cor1pi02} Let $X$ be a Polish space, and $A,B$ be disjoint analytic relations on $X$ such that $A$ is either s-acyclic, or locally countable. Then exactly one of the following holds:\smallskip  

(a) the set $A$ is separable from $B$ by a $\mbox{pot}(\bormtwo )$ set,\smallskip  

(b) there is $f\! :\! 2^\omega\!\rightarrow\! X$ injective continuous such that the inclusions 
$\lceil T\rceil\cap\mathbb{E}_0\!\subseteq\! (f\!\times\! f)^{-1}(A)$ and 
$\lceil T\rceil\!\setminus\!\mathbb{E}_0\!\subseteq\! (f\!\times\! f)^{-1}(B)$ hold.\end{cor}

\noindent\bf Proof.\rm ~By Lemma \ref{E0}, $\lceil T\rceil\cap\mathbb{E}_0$ is not separable from $\lceil T\rceil\!\setminus\!\mathbb{E}_0$ by a $\mbox{pot}(\bormtwo )$ set. This shows that (a) and (b) cannot hold simultaneously. So assume that (a) does not hold. By Lemma \ref{anaKsigma}, we may assume that $A$ is $\boratwo$ and $B$ is the complement of $A$. By Lemma \ref{suffqa}, we may also assume that $A$ is quasi-acyclic. It remains to apply Theorem \ref{delta02}.
\hfill{$\square$}

\begin{cor} \label{cor2pi02} Let $X,Y$ be Polish spaces, and $A,B$ be disjoint analytic subsets of $X\!\times\! Y$ such that $A$ is locally countable. Then exactly one of the following holds:\smallskip  

(a) the set $A$ is separable from $B$ by a $\mbox{pot}(\bormtwo )$ set,\smallskip  

(b) $(2^\omega ,2^\omega ,\lceil T\rceil\cap\mathbb{E}_0,\lceil T\rceil\!\setminus\!\mathbb{E}_0)
\sqsubseteq (X,Y,A,B)$.\end{cor}

\noindent\bf Proof.\rm ~As in the proof of Corollary \ref{cor1pi02}, (a) and (b) cannot hold simultaneously. So assume that (a) does not hold. We argue as in the proof of Corollary 
\ref{cor2checkD2Sigma01}. Corollary \ref{cor1pi02} gives $f'\! :\! 2^\omega\!\rightarrow\! Z$.
\hfill{$\square$}

\begin{cor} \label{caracpartialpi02} Let $X$ be a Polish space, and $A,B$ be disjoint analytic relations on $X$. The following are equivalent:\smallskip

\noindent (1) there is an s-acyclic relation $R\!\in\!\ana$ such that $A\cap R$ is not separable from $B\cap R$ by a $\mbox{pot}(\bormtwo )$ set,\smallskip

\noindent (2) there is $f\! :\! 2^\omega\!\rightarrow\! X$ injective continuous with 
$\lceil T\rceil\cap\mathbb{E}_0\!\subseteq\! (f\!\times\! f)^{-1}(A)$ and 
$\lceil T\rceil\!\setminus\!\mathbb{E}_0\!\subseteq\! (f\!\times\! f)^{-1}(B)$.\end{cor} 

\noindent\bf Proof.\rm ~(1) $\Rightarrow$ (2) We apply Corollary \ref{cor1pi02}.\bigskip

\noindent (2) $\Rightarrow$ (1) We can take $R\! :=\! (f\!\times\! f)\big[\lceil T\rceil\big]$.\hfill{$\square$}\bigskip

\noindent\bf Remark.\rm ~There is a version of Corollary \ref{caracpartialpi02} for $\boratwo$ instead of $\bormtwo$, obtained by exchanging the roles of $A$ and $B$. This symmetry is not present in Theorem \ref{delta02}.

\begin{cor} \label{partialsigma02} Let $X$ be a Polish space, and $A,B$ be disjoint analytic relations on $X$ such that $A$ is contained in a pot$(F_\sigma )$ s-acyclic relation, or $A\cup B$ is s-acyclic. Then exactly one of the following holds:\smallskip  

(a) the set $A$ is separable from $B$ by a $\mbox{pot}(\boratwo )$ set,\smallskip  

(b) there is $f\! :\! 2^\omega\!\rightarrow\! X$ injective continuous such that the inclusions 
$\lceil T\rceil\!\setminus\!\mathbb{E}_0\!\subseteq\! (f\!\times\! f)^{-1}(A)$ and 
$\lceil T\rceil\cap\mathbb{E}_0\!\subseteq\! (f\!\times\! f)^{-1}(B)$ hold.\end{cor}

\noindent\bf Proof.\rm ~Let $R$ be a pot$(F_\sigma )$ s-acyclic relation containing $A$. Then 
there is no $\mbox{pot}(\boratwo )$ set $P$ separating $A\cap R\! =\! A$ from $B\cap R$, since otherwise $P\cap R\!\in\!\mbox{pot}(\boratwo )$ and separates $A$ from $B$. Corollary 
\ref{caracpartialpi02} gives $f\! :\! 2^\omega\!\rightarrow\! X$ injective continuous with 
$\lceil T\rceil\cap\mathbb{E}_0\!\subseteq\! (f\!\times\! f)^{-1}(B)$ and 
$\lceil T\rceil\!\setminus\!\mathbb{E}_0\!\subseteq\! (f\!\times\! f)^{-1}(A)$.\bigskip

 If $A\cup B$ is s-acyclic, then we apply Corollary \ref{cor1pi02}.\hfill{$\square$}\bigskip

\noindent\bf Remarks.\rm\ (1) Corollary \ref{partialsigma02} also holds when $A\cup B$ is locally countable, but we did not mention it in the statement since (a) always holds in this case. Indeed, by reflection, $A\cup B$ is contained in a locally countable Borel set $C$. As $A,B$ are disjoint analytic sets, there is a Borel set $D$ separating $A$ from $B$. Thus $C\cap D$ is a locally countable Borel set separating $A$ from $B$. But a locally countable Borel set has $\boratwo$ vertical sections, and is therefore $\mbox{pot}(\boratwo )$ (see [Lo2]).\bigskip

\noindent (2) There is a version of Corollary \ref{partialsigma02} for ${\bf\Gamma}\! =\!\boraone$, where we replace the class $F_\sigma$ with the class of open sets. We do not state it since (a) always holds in this case. Indeed, a potentially open s-acyclic relation is a countable union of Borel rectangles for which at least one side is a singleton, so that this union is potentially clopen, just like any of its Borel subsets.

\section{$\!\!\!\!\!\!$ The class $\bortwo$}

$\underline{\mbox{\bf Example}}$\bigskip

 We set, for each $\varepsilon\!\in\! 2$, 
$$\mathbb{E}^\varepsilon_0\! :=\!\{ (\alpha ,\beta )\!\in\! 2^\omega\!\times\! 2^\omega\mid
\exists m\! >\! 0~~\alpha (m)\!\not=\!\beta (m)\wedge\forall n\! >\! m~~\alpha (n)\! =\!\beta (n)\wedge 
(m\! -\! 1)_0\!\equiv\!\varepsilon\ (\mbox{mod }2)\} .$$

\begin{lem} \label{E00} $\lceil T\rceil\cap\mathbb{E}^0_0$ is not separable from 
$\lceil T\rceil\cap\mathbb{E}^1_0$ by a $\mbox{pot}(\bortwo )$ set.\end{lem}

\noindent\bf Proof.\rm ~The proof is similar to that of Lemma \ref{E0}. We argue by contradiction, which gives $D$ in $\mbox{pot}(\bortwo )$, and also a dense $G_\delta$ subset $G$ of 
$2^\omega$ such that $D\cap G^2\!\in\! \bortwo (G^2)$. Let $(O_n)_{n\in\omega}$ be a sequence of dense open subsets of $2^\omega$ with intersection $G$. Note that 
$\lceil T\rceil\cap\mathbb{E}^0_0\cap G^2\!\subseteq\!\lceil T\rceil\cap D\cap G^2$, 
$\lceil T\rceil\cap\mathbb{E}^1_0\cap G^2\!\subseteq\!\lceil T\rceil\cap G^2\!\setminus\! D$ and 
$\lceil T\rceil\cap D\cap G^2\!\in\!\bortwo (\lceil T\rceil\cap G^2)$. By Baire's theorem, it is enough to prove that ${\lceil T\rceil\cap\mathbb{E}^0_0\cap G^2}$ and 
$\lceil T\rceil\cap\mathbb{E}^1_0\cap G^2$ are dense in $\lceil T\rceil\cap G^2$. Let us do it for 
$\lceil T\rceil\cap\mathbb{E}^0_0\cap G^2$, the other case being similar. So let $q\!\in\!\omega$ and $w\!\in\! 2^{<\omega}$. Pick $N\!\in\!\omega$ such that 
$(s_q0w0^{N_0},t_q1w0^N)$ is in ${\cal F}$ and ${(\vert s_q0w0^N\vert\! -\! 1)_0\! =\! 0}$. Then we argue as in the proof of of Lemma \ref{E0}: pick $u_0\!\in\! 2^\omega$ with 
$N_{s_q0w0^N0u_0}\!\subseteq\! O_0$, 
$v_0\!\in\! 2^\omega$ with $N_{t_q1w0^N1u_0v_0}\!\subseteq\! O_0$,  
$u_1\!\in\! 2^\omega$ with $N_{s_q0w0^N0u_0v_0u_1}\!\subseteq\! O_1$, 
$v_1\!\in\! 2^\omega$ with $N_{t_q1w0^N1u_0v_0u_1v_1}\!\subseteq\! O_1$, and so on. Then 
$(s_q0w0^N0u_0v_0u_1v_1...,t_q1w0^N1u_0v_0u_1v_1...)$ is in 
$\lceil T\rceil\cap\mathbb{E}^0_0\cap G^2$.\hfill{$\square$}\bigskip

\noindent $\underline{\mbox{\bf The main result}}$\bigskip

 We will prove a version of Theorem \ref{delta02} for the class $\bortwo$.

\begin{thm} \label{delt02} Let $X$ be a Polish space, and $A,B$ be disjoint analytic relations on $X$ such that $A\cup B$ is quasi-acyclic. Then one of the following holds:\smallskip  

(a) the set $A$ is separable from $B$ by a $\mbox{pot}(\bortwo )$ set,\smallskip  

(b) there is $f\! :\! 2^\omega\!\rightarrow\! X$ injective continuous such that the inclusions 
$\lceil T\rceil\cap\mathbb{E}^0_0\!\subseteq\! (f\!\times\! f)^{-1}(A)$ and 
$\lceil T\rceil\cap\mathbb{E}^1_0\!\subseteq\! (f\!\times\! f)^{-1}(B)$ hold.\end{thm}

\noindent\bf Proof.\rm ~The proof is similar to that of of Theorem \ref{delta02}. Assume that (a) does not hold. By Lemma \ref{anaKsigmas}, we may assume that $A,B$ are $\boratwo$. Let 
$(C_n)_{n\in\omega}$ be a witness for the fact that $A\cup B$ is quasi-acyclic. As $A,B$ are 
$\boratwo$, we may assume that each $C_n$ is either contained in $A$, or contained in $B$. Note that there are disjoint Borel subsets $O_0,O_1$ of $X$ such that $A\cap (O_0\!\times\! O_1)$ is not separable from $B\cap (O_0\!\times\! O_1)$ by a $\mbox{pot}(\bortwo )$ set. We may assume that $X$ is zero-dimensional, the $C_n$'s are closed, and $O_0,O_1$ are clopen, refining the topology if necessary. We can also replace $A,B$ and the $C_n$'s with their intersection with $O_0\!\times\! O_1$ and assume that they are contained in $O_0\!\times\! O_1$. This gives a sequence $(C^0_n)_{n\in\omega}$ (resp., $(C^1_n)_{n\in\omega}$) of pairwise disjoint closed relations on $X$ with union $A$ (resp., $B$).\bigskip

\noindent $\bullet$ We may assume that $X$ is recursively presented, $O_0,O_1$ are $\Borel$ and the relation ``$(x,y)\!\in\! C^\varepsilon_n$" is $\Borel$ in $(x,y,\varepsilon ,n)$. As ${\it\Delta}_X$ is Polish finer than the topology on $X$, $A$ is not separable from $B$ by a $\bortwo (X^2,\tau_1)$ set. We set, for $F\!\in\!\bormone (X^2,\tau_1)$, 
$F'\! :=\!\overline{F\cap A}^{\tau_1}\cap\overline{F\cap B}^{\tau_1}$ (see 22.30 in [K]). Then 
$$F_\theta\! =\! F_{\theta +1}\! =\! F'_\theta\! =\!\overline{F_\theta\cap A}^{\tau_1}\cap\overline{F_\theta\cap B}^{\tau_1}\mbox{,}$$ 
so that $F_\theta\cap A$ and $F_\theta\cap B$ are $\tau_1$-dense in $F_\theta$.\bigskip

\noindent $\bullet$ Let us prove that $F_\theta$ is not empty. We argue by contradiction, so that 
$A\! =\!\bigcup_{\xi <\theta}~A\cap F_\xi\!\setminus\! F_{\xi +1}$. But 
$A\cap F_\xi\!\setminus\! F_{\xi +1}\! =\! A\cap F_\xi\!\setminus\! 
(\overline{F_\xi\cap A}^{\tau_1}\cap\overline{F_\xi\cap B}^{\tau_1})\!\subseteq\! 
F_\xi\!\setminus\!\overline{F_\xi\cap B}^{\tau_1}\!\subseteq\!\neg B$. This means that 
$(F_\xi\!\setminus\! F_{\xi +1})_{\xi <\theta}$ is a countable partition of $(X^2,\tau_1)$ into $\bortwo$ sets, and that $A$ is separable from $B$ by a $\bortwo$ set on each piece of the partition. This implies that $A$ is separable from $B$ by a $\bortwo (X^2,\tau_1)$ set, which is absurd.\bigskip

\noindent $\bullet$ As in the proof of Theorem \ref{delta02}, $F_\theta$ is $\Ana$, and the following key property holds:
$$\forall\varepsilon\!\in\! 2~~\forall q\!\in\!\omega ~~\forall U,V\!\in\!\Ana (X)~~
F_\theta\cap (U\!\times\! V)\!\not=\!\emptyset\Rightarrow\exists n\!\geq\! q~~
F_\theta\cap C^\varepsilon_n\cap (U\!\times\! V)\!\not=\!\emptyset .$$
$\bullet$ We construct again sequences $(x_s)$, $(y_s)$, $(X_s)$, $(Y_s)$ and $\Phi$ satisfying the following conditions:
$$\begin{array}{ll}
& (1)\ x_s\!\in\! X_s\ \wedge\ y_s\!\in\! Y_s\cr
& (2)\ X_{s\varepsilon}\!\subseteq\! X_s\!\subseteq\!\Omega_X\cap O_0\ \wedge\ 
Y_{s\varepsilon}\!\subseteq\! Y_s\!\subseteq\!\Omega_X\cap O_1\cr
& (3)\ \mbox{diam}_{\mbox{GH}}(X_s)\mbox{, }\mbox{diam}_{\mbox{GH}}(Y_s)\!\leq\! 
2^{-\vert s\vert}\cr
& (4)\ (x_{u_0},y_{u_1})\!\in\! F_\theta\cap C^\varepsilon_{\Phi(\vec t(\vec u))}
\mbox{ if }(\vert\vec t (\vec u)\vert\! -\! 2)_0\!\equiv\!\varepsilon\ (\mbox{mod }2)
\mbox{, with the convention }(-1)_0\! =\! 0\cr
& (5)\ X_{s0}\cap X_{s1}\! =\! Y_{s0}\cap Y_{s1}\! =\!\emptyset
\end{array}$$
$\bullet$ Assume that this has been done. If $(\alpha ,\beta )\!\in\!\lceil T\rceil\cap\mathbb{E}^0_0$, then $\Phi (\vec t\big( (\alpha ,\beta )\vert n\big) )\! =\! N$ if $n$ is big enough. In this case, by (4), 
$(x_{\alpha\vert n},y_{\beta\vert n})\!\in\! C^0_N$ which is closed, so that 
$\big( f(\alpha ),g(\beta )\big)\!\in\! C^0_N\!\subseteq\! A$. Similarly, if 
$(\alpha ,\beta )\!\in\!\lceil T\rceil\cap\mathbb{E}^1_0$, then 
$\big( f(\alpha ),g(\beta )\big)\!\in\! C^1_N\!\subseteq\! B$.\bigskip

\noindent $\bullet$ So let us prove that the construction is possible. The key property gives 
$\Phi (0,1)\!\in\!\omega$ and $(x_0,y_1)$ in $F_\theta\cap C^0_{\Phi (0,1)}\cap\Omega_{X^2}$. We choose $\Ana$ subsets $X_0,Y_1$ of $X$ with GH-diameter at most $2^{-1}$ such that 
$(x_0,y_1)\!\in\! X_0\!\times\! Y_1\!\subseteq\! (\Omega_X\cap O_0)\!\times\! (\Omega_X\cap O_1)$, which completes the construction for the length $l\! =\! 1$.\bigskip

 Let $l\!\geq\! 1$. We now want to build $x_s,X_s,y_s,Y_s$ for $s\!\in\! 2^{l+1}$, as well as 
$\Phi (s_l0,t_l1)$. Fix $\eta\!\in\! 2$ such that $(l\! -\! 1)_0\!\equiv\!\eta\ (\mbox{mod }2)$. Note that 
$(x_{s_l},y_{t_l})\!\in\! F_\theta\cap (U\!\times\! V)$, where 
$$\begin{array}{ll}
& U\! :=\!\{ x'_{s_l}\!\in\! X_{s_l}\mid
\exists (x'_s)_{s\in 2^l\setminus\{ s_l\} ,s(0)=0}\!\in\!\Pi_{s\in 2^l\setminus\{ s_l\} ,s(0)=0}~X_s~~
\exists (y'_s)_{s\in 2^l,s(0)=1}\!\in\!\Pi_{s\in 2^l,s(0)=1}~Y_s\cr
& \hfill{\forall\vec u\!\in\! T\cap (2^l\!\times\! 2^l)~~
(x'_{u_0},y'_{u_1})\!\in\! F_\theta\cap C^\varepsilon_{\Phi (\vec t(\vec u))}
\mbox{ if }(\vert\vec t (\vec u)\vert\! -\! 2)_0\!\equiv\!\varepsilon\ (\mbox{mod }2)\}\mbox{,}}\cr\cr
& V\! :=\!\{ y'_{t_l}\!\in\! Y_{t_l}\mid
\exists (x'_s)_{s\in 2^l,s(0)=0}\!\in\!\Pi_{s\in 2^l,s(0)=0}~X_s~~
\exists (y'_s)_{s\in 2^l\setminus\{ t_l\} ,s(0)=1}\!\in\!\Pi_{s\in 2^l\setminus\{ t_l\} ,s(0)=1}~Y_s\cr
& \hfill{\forall\vec u\!\in\! T\cap (2^l\!\times\! 2^l)~~
(x'_{u_0},y'_{u_1})\!\in\! F_\theta\cap C^\varepsilon_{\Phi (\vec t(\vec u))}
\mbox{ if }(\vert\vec t (\vec u)\vert\! -\! 2)_0\!\equiv\!\varepsilon\ (\mbox{mod }2)\} .}
\end{array}$$
The key property gives $\Phi (s_l0,t_l1)\! >\!\mbox{max}_{q<l}~\Phi (s_q0,t_q1)$ and 
$$(x_{s_l0},y_{t_l1})\!\in\! F_\theta\cap C^\eta_{\Phi (s_l0,t_l1)}\cap (U\!\times\! V).$$ 
Note that $x_{s_l0}\!\not=\! x_{s_l1}$ because $(x_{s_l0},y_{t_l1})\!\in\! C^\eta_{\Phi (s_l0,t_l1)}$,  
$(x_{s_l1},y_{t_l1})\!\in\! C^\varepsilon_{\Phi (\vec t (s_l1,t_l1))}$ if 
$$(\vert\vec t (s_l1,t_l1)\vert\! -\! 2)_0\!\equiv\!\varepsilon\ (\mbox{mod }2)\mbox{,}$$ 
and $\Phi (s_l0,t_l1)\! >\!\Phi\big(\vec t (s_l1,t_l1)\big)$. Similarly, $y_{t_l0}\!\not=\! y_{t_l1}$. If $s\!\in\! 2^l$, then there is an injective $s(T)$-path $p_s$ from $s$ to $s_l$. This gives a $s(A\cup B)$-path from $x_{s0}$ to $x_{s1}$ if $s(0)\! =\! 0$, and a $s(A\cup B)$-path from $y_{s0}$ to $y_{s1}$ if $s(0)\! =\! 1$. Using the quasi-acyclicity of $s(A\cup B)$, we see, by induction on the length of $p_s$, that $x_{s0}\!\not=\! x_{s1}$ and $y_{s0}\!\not=\! y_{s1}$.\hfill{$\square$}\bigskip

\noindent $\underline{\mbox{\bf Consequences}}$\bigskip

\begin{cor} \label{cor1delta02} Let $X$ be a Polish space, and $A,B$ be disjoint analytic relations on $X$ such that\smallskip

- either $A\cup B$ is either s-acyclic or locally countable\smallskip

- or $A$ is contained in a $\mbox{pot}(\bortwo )$ s-acyclic or locally countable relation.\smallskip

\noindent Then exactly one of the following holds:\smallskip  

(a) the set $A$ is separable from $B$ by a $\mbox{pot}(\bortwo )$ set,\smallskip  

(b) there is $f\! :\! 2^\omega\!\rightarrow\! X$ injective continuous such that the inclusions 
$\lceil T\rceil\cap\mathbb{E}^0_0\!\subseteq\! (f\!\times\! f)^{-1}(A)$ and 
$\lceil T\rceil\cap\mathbb{E}^1_0\!\subseteq\! (f\!\times\! f)^{-1}(B)$ hold.\end{cor}

\noindent\bf Proof.\rm ~By Lemma \ref{E00}, $\lceil T\rceil\cap\mathbb{E}^0_0$ is not separable from $\lceil T\rceil\cap\mathbb{E}^1_0$ by a $\mbox{pot}(\bortwo )$ set. This shows that (a) and (b) cannot hold simultaneously. So assume that (a) does not hold.\bigskip
 
\noindent - If $A\cup B$ is s-acyclic or locally countable, then by Lemma \ref{anaKsigmas}, we may assume that $A,B$ are $\boratwo$. By Lemma \ref{suffqa}, we may also assume that $A\cup B$ is quasi-acyclic. It remains to apply Theorem \ref{delt02}.\bigskip
 
\noindent - Assume that $R$ is $\mbox{pot}(\bortwo )$ and contains $A$. Then there is no 
$\mbox{pot}(\bortwo )$ set $P$ separating $A\cap R\! =\! A$ from $B\cap R$, since otherwise 
$P\cap R\!\in\!\mbox{pot}(\bortwo )$ separates $A$ from $B$. It remains to apply the first point. This finishes the proof.\hfill{$\square$}

\begin{cor} \label{cor2delta02} Let $X,Y$ be Polish spaces, and $A,B$ be disjoint analytic subsets of $X\!\times\! Y$ such that $A\cup B$ is locally countable or $A$ is contained in a 
$\mbox{pot}(\bortwo )$ locally countable set. Then exactly one of the following holds:\smallskip  

(a) the set $A$ is separable from $B$ by a $\mbox{pot}(\bortwo )$ set,\smallskip  

(b) $(2^\omega ,2^\omega ,\lceil T\rceil\cap\mathbb{E}^0_0,\lceil T\rceil\cap\mathbb{E}^1_0)
\sqsubseteq (X,Y,A,B)$.\end{cor}

\noindent\bf Proof.\rm ~As in the proof of Corollary \ref{cor1delta02}, (a) and (b) cannot hold simultaneously. Then we argue as in the proof of Corollary \ref{cor2checkD2Sigma01}. The set $A'\cup B'$ is locally countable or $A'$  is contained in a $\mbox{pot}(\bortwo )$ locally countable set, and $A'$ is not separable from $B'$ by a $\mbox{pot}(\bortwo )$ set. Corollary 
\ref{cor1delta02} gives $f'\! :\! 2^\omega\!\rightarrow\! Z$.\hfill{$\square$}

\begin{cor} \label{caracpartialdelta02} Let $X$ be a Polish space, and $A,B$ be disjoint analytic relations on $X$. The following are equivalent:\smallskip

\noindent (1) there is an s-acyclic or locally countable relation $R\!\in\!\ana$ such that $A\cap R$ is not separable from $B\cap R$ by a $\mbox{pot}(\bortwo )$ set,\smallskip

\noindent (2) there is $f\! :\! 2^\omega\!\rightarrow\! X$ injective continuous with 
$\lceil T\rceil\cap\mathbb{E}^0_0\!\subseteq\! (f\!\times\! f)^{-1}(A)$ and 
$\lceil T\rceil\cap\mathbb{E}^1_0\!\subseteq\! (f\!\times\! f)^{-1}(B)$.\end{cor} 

\noindent\bf Proof.\rm ~(1) $\Rightarrow$ (2) We apply Corollary \ref{cor1delta02}.\bigskip

\noindent (2) $\Rightarrow$ (1) We can take $R\! :=\! (f\!\times\! f)\big[\lceil T\rceil\cap\mathbb{E}_0\big]$.
\hfill{$\square$}

\section{$\!\!\!\!\!\!$ The classes $D_n(\boratwo )$ and $\check D_n(\boratwo )$}

$\underline{\mbox{\bf Examples}}$\bigskip

\noindent\bf Notation.\rm\ Let $\eta\!\geq\! 1$ be a countable ordinal, and 
$S_\eta\! :\!\omega\!\rightarrow\!\eta$ be onto. We set 
$$C_0\! :=\!\{\alpha\!\in\! 2^\omega\mid\exists m\!\in\!\omega ~~
\forall p\!\geq\! m~~\alpha (p)\! =\! 0\}$$ 
and, for $1\!\leq\!\theta\! <\!\eta$, 
$C_\theta\! :=\!\big\{\alpha\!\in\! 2^\omega\mid\exists m\!\in\!\omega ~~
\forall p\!\in\!\omega ~~\alpha (<m,p>)\! =\! 0\wedge S_\eta\big( (m)_0\big)\!\leq\!\theta\big\}$, so that 
$(C_\theta )_{\theta <\eta}$ is an increasing sequence of $\boratwo$ subsets of $2^\omega$. We then set $D_\eta\! :=\! D\big( (C_\theta )_{\theta <\eta}\big)$. 

\begin{lem} \label{ExCn} The set $D_\eta$ is $D_\eta (\boratwo )$-complete.\end{lem} 

\noindent\bf Proof.\rm ~By 21.14 in [K], it is enough to see that $D_\eta$ is not 
$\check D_\eta (\boratwo )$ since it is $D_\eta (\boratwo )$. We will prove more. Let us say that a pair $(\theta ,F)$ is {\bf suitable} if $\theta\!\leq\!\eta$, $F$ is a chain of finite binary sequences, $I_F\! :=\!\bigcap_{s\in F}~\{\alpha\!\in\! N_s\mid (\alpha )_{\vert s\vert}\! =\! 0^\infty\}$ is not empty and $S_\eta\big( (\vert s\vert )_0\big)\!\geq\!\theta$ for each $s\!\in\! F$. Let us prove that 
$I_F\cap D\big( (C_{\theta'})_{\theta'<\theta}\big)$ is not $\check D_\theta (\boratwo )$ if 
$(\theta ,F)$ is suitable. This will give the result since $(\eta ,\emptyset )$ is suitable.\bigskip

 We argue by induction on $\theta$. If $\theta\! =\! 1$, then the $\boratwo$ set $I_F\cap C_0$ is dense and co-dense in the closed set $I_F$, so that it is not $\bormtwo$, by Baire's theorem. Assume the result proved for ${\theta'\! <\!\theta}$. We argue by contradiction, which gives an increasing sequence $(H_{\theta'})_{\theta' <\theta}$ of $\boratwo$ sets with 
$${I_F\cap D\big( (C_{\theta'})_{\theta'<\theta}\big)\! =\!
\neg D\big( (H_{\theta'})_{\theta' <\theta}\big)}.$$ 

 As $\neg (\bigcup_{\theta'<\theta}~C_{\theta'})$ is comeager in $I_F$, 
$I_F\cap\bigcup_{\theta'<\theta}~H_{\theta'}$ too, which gives $\theta'\! <\!\theta$ with parity opposite to that of $\theta$ and $s'\!\supseteq\!\mbox{max}_{s\in F}~s$ such that 
$S_\eta\big( (\vert s'\vert)_0\big)\! =\!\theta'$ and 
$\emptyset\!\not=\! I_F\cap N_{s'}\!\subseteq\! H_{\theta'}$. We set $F'\! :=\! F\cup\{ s'\}$, so that 
$(\theta' ,F')$ is suitable. By induction assumption, 
$I_{F'}\cap D\big( (C_{\theta''})_{\theta''<\theta'}\big)$ is not $\check D_{\theta'}(\boratwo )$. But 
$I_{F'}\cap D\big( (C_{\theta''})_{\theta''<\theta'}\big)\! =\! 
I_{F'}\!\setminus\! D\big( (H_{\theta''})_{\theta''<\theta'}\big)\!\in\!\check D_{\theta'}(\boratwo )$ since 
$I_{F'}\!\subseteq\! C_{\theta'}$, which is absurd.\hfill{$\square$}\bigskip

\noindent\bf Notation.\rm\ We now fix an effective frame in the sense of Definition 2.1 in [L8], which are frames in the sense of Definition \ref{frame}. Lemma 2.3 in [L8] proves the existence of such an effective frame. Note that $(s_1,t_1)\! =\! (0,1)$, so that $s_1(0)\!\not=\! t_1(0)$. But 
$s_{l+1}(l)\! =\! t_{l+1}(l)$ if $l\!\geq\! 1$. Indeed, it is enough to see that 
$\Big(\big( (l)_1\big)_1\Big)_0\! +\!\Big(\big( (l)_1\big)_1\Big)_1\! <\! l$ in this case, by the proof of Lemma 2.3 in [L8]. As $(q)_0\! +\! (q)_1\!\leq\! q$, and $(q)_0\! +\! (q)_1\! <\! q$ if $q\!\geq\! 2$, we may assume that $\big( (l)_1\big)_1\!\in\! 2$. If $\big( (l)_1\big)_1\! =\! 0$, then we are done since 
$l\!\geq\! 1$. If $\big( (l)_1\big)_1\! =\! 1$, then $l\!\geq\! 2$ and we are done too.\bigskip

\noindent $\bullet$ The {\bf shift map} ${\cal S}\! :\! 2^L\!\rightarrow\! 2^{L-1}$ is defined by 
${\cal S}(\alpha )(m)\! :=\!\alpha (m\! +\! 1)$ when $1\!\leq\! L\!\leq\!\omega$, with the convention $\omega\! -\! 1\! :=\!\omega$.\bigskip

\noindent $\bullet$ The {\bf symmetric difference} $\alpha\Delta\beta$ of 
$\alpha ,\beta\!\in\! 2^L$ is the element of $2^L$ defined by 
$(\alpha\Delta\beta )(m)\! =\! 1$ exactly when $\alpha (m)\!\not=\!\beta (m)$, if $L\!\leq\!\omega$.\bigskip

\noindent $\bullet$ We set $\mathbb{N}_\eta\! :=\!
\{ (\alpha ,\beta )\!\in\!\lceil T\rceil\mid {\cal S}(\alpha\Delta\beta )\!\notin\! D_\eta\}$. 

\begin{lem} \label{ExDnS2} The $\check D_\eta (\boratwo )$ set $\mathbb{N}_\eta$ is not separable from $\lceil T\rceil\!\setminus\!\mathbb{N}_\eta$ by a 
$\mbox{pot}\big( D_\eta (\boratwo )\big)$ set.\end{lem} 

\noindent\bf Proof.\rm ~As $\lceil T\rceil$ is closed, $D_\eta$ is $D_\eta (\boratwo )$ and 
${\cal S},\Delta$ are continuous, $\mathbb{N}_\eta$ is $\check D_\eta (\boratwo )$. By Lemma 2.6 in [L8], it is enough to check that $D_\eta$ is ccs (see Definition 2.5 in [L8]). We just have to check that the $C_\theta$'s are ccs. So let $\alpha ,\alpha_0\!\in\! 2^\omega$ and 
$F\! :\! 2^\omega\!\rightarrow\! 2^\omega$ satisfying the conclusion of Lemma 2.4.(b) in [L8]. Note that $\alpha\!\in\! C_0$ exactly when 
$\{ m\!\in\!\omega\mid\alpha (m)\! =\! 1\}$ is finite, so that $C_0$ is ccs. If $\theta\!\geq\! 1$, then  
$\alpha\!\notin\! C_\theta$ exactly when, for each $m$, $S_\eta\big( (m)_0\big)\!\leq\!\theta$ or there is $p$ with $\alpha (<m,p>)\! =\! 1$. As 
$\big( B_\alpha (<m,p>)\big)_0\! =\! (<m,p>)_0\! =\! m$, $C_\theta$ is ccs too.\hfill{$\square$}\bigskip

\noindent $\underline{\mbox{\bf The main result}}$\bigskip

\noindent\bf Notation.\rm\ From now on, $\eta\! <\!\omega$. We set, for $2\!\leq\!\theta\!\leq\!\eta$ and $(s,t)\!\in\! (2\!\times\! 2)^{<\omega}\!\setminus\!\{ (\emptyset ,\emptyset )\}$, 
$$m^\theta_{s,t}\! :=\!\mbox{min}\big\{ m\!\in\!\omega\mid\big( {\cal S}(s\Delta t)\big)_m\!\subseteq\! 0^\infty\wedge S_\eta\big( (m)_0\big)\! <\!\theta\big\} .$$ 
We also set $s^-\! :=<s(0),...,s(\vert s\vert\! -\! 2)>$ if $s\!\in\! 2^{<\omega}$.

\vfill\eject

\noindent $\bullet$ We define the following relation on $(2\!\times\! 2)^{<\omega}$:\bigskip

\leftline{$(s,t)~R~(s',t')\Leftrightarrow (s,t)\!\subseteq\! (s',t')~\wedge ~
\bigg(\vert s\vert\!\leq\! 1~\vee ~\Big(\vert s\vert\!\geq\! 2~\wedge ~
\exists 2\!\leq\!\theta\!\leq\!\eta ~~m^\theta_{s,t}\!\not=\! m^\theta_{s^-,t^-}~\wedge$}\smallskip

\rightline{$\forall (s,t)\!\subseteq\! (s'',t'')\!\subseteq\! (s',t')~~\forall\theta\! <\!\theta'\!\leq\!\eta ~~
m^{\theta'}_{s,t}\! =\! m^{\theta'}_{s^-,t^-}\! =\! m^{\theta'}_{s'',t''}\Big)~\vee$}\smallskip

\rightline{$\Big(\vert s\vert\!\geq\! 2~\wedge ~s(\vert s\vert\! -\! 1)\!\not=\! t(\vert s\vert\! -\! 1)~
\wedge$}\smallskip

\rightline{$\forall (s,t)\!\subseteq\! (s'',t'')\!\subseteq\! (s',t')~~\forall 2\!\leq\!\theta\!\leq\!\eta ~~
m^\theta_{s,t}\! =\! m^\theta_{s^-,t^-}\! =\! m^\theta_{s'',t''}\Big)~\vee$}\smallskip

\rightline{$\Big(\vert s\vert\!\geq\! 2~\wedge ~\forall (s,t)\!\subseteq\! (s'',t'')\!\subseteq\! (s',t')~~
\big(\forall 2\!\leq\!\theta\!\leq\!\eta ~~m^\theta_{s,t}\! =\! m^\theta_{s^-,t^-}\! =\! m^\theta_{s'',t''}\big)~\wedge$}\smallskip

\rightline{$s''(\vert s''\vert\! -\! 1)\! =\! t''(\vert s''\vert\! -\! 1)\Big)\bigg).$}\bigskip

\noindent Note that $R$ is a {\bf tree relation}, which means that it is a partial order (it contains the diagonal, is antisymmetric and transitive) with minimum  element $(\emptyset ,\emptyset )$, the set of predecessors of any sequence is finite and lineary ordered by $R$. Moreover, $R$ is 
{\bf distinguished} in $\subseteq$, which means that $(s,t)~R~(s',t')$ if 
$(s,t)\!\subseteq\! (s',t')\!\subseteq\! (s'',t'')$ and $(s,t)~R~(s'',t'')$ (see [D-SR]).\bigskip

\noindent $\bullet$ We set
$$\begin{array}{ll}
& D_\eta\! :=\!\{ (s,t)\!\in\! T\mid\vert s\vert\!\geq\! 2\Rightarrow m^\eta_{s,t}\!\not=\! m^\eta_{s^-,t^-}\}\mbox{ if }\eta\!\geq\! 2\mbox{,}\cr\cr
& D_\theta\! :=\!\{ (s,t)\!\in\! T\mid\vert s\vert\!\geq\! 2~\wedge ~m^\theta_{s,t}\!\not=\! m^\theta_{s^-,t^-}~\wedge ~\forall\theta\! <\!\theta'\!\leq\!\eta ~~m^{\theta'}_{s,t}\! =\! m^{\theta'}_{s^-,t^-}\}
\mbox{ if }2\!\leq\!\theta\! <\!\eta\mbox{,}\cr\cr
& D_1\! :=\!\{ (s,t)\!\in\! T\mid\vert s\vert\!\geq\! 2~\wedge ~\forall 2\!\leq\!\theta\!\leq\!\eta ~~
m^\theta_{s,t}\! =\! m^\theta_{s^-,t^-}~\wedge ~s(\vert s\vert\! -\! 1)\!\not=\! t(\vert s\vert\! -\! 1)\}
\mbox{,}\cr\cr
& D_0\! :=\!\{ (s,t)\!\in\! T\mid\vert s\vert\!\geq\! 2~\wedge ~s(\vert s\vert\! -\! 1)\! =\! t(\vert s\vert\! -\! 1)\}\mbox{,}\cr
\end{array}$$
so that the $(D_\theta )_{\theta\leq\eta}$ is a partition of $T$.

\begin{thm} \label{checkDnSigma02} Let $1\!\leq\!\eta\! <\!\omega$. Let $X$ be a Polish space, and $A_0,A_1$ be disjoint analytic relations on $X$ such that $A_0\cup A_1$ is s-acyclic. Then exactly one of the following holds:\smallskip  

(a) the set $A_0$ is separable from $A_1$ by a $\mbox{pot}\big( D_\eta (\boratwo )\big)$ set,\smallskip  

(b) $(2^\omega ,2^\omega ,\mathbb{N}_\eta ,\lceil T\rceil\!\setminus\!\mathbb{N}_\eta )\sqsubseteq (X,X,A_0,A_1)$, via a square map.\end{thm}

\noindent\bf Proof.\rm ~By Lemma \ref{ExDnS2}, (a) and (b) cannot hold simultaneously. So assume that (a) does not hold. Note first that we may assume that 
$A_0\cup A_1$ is compact and $A_1$ is $D_\eta (\boratwo )$. Indeed, Theorems 1.9 and 1.10 in [L8] give $\mathbb{S}\!\in\! D_\eta (\boratwo )(\lceil T\rceil )$ and 
$f',g'\! :\! 2^\omega\!\rightarrow\! X$ continuous such that the inclusions 
$\mathbb{S}\!\subseteq\! (f'\!\times\! g')^{-1}(A_1)$ and 
$\lceil T\rceil\!\setminus\!\mathbb{S}\!\subseteq\! (f'\!\times\! g')^{-1}(A_0)$ hold. Let 
$(\Sigma_\theta )_{\theta <\eta}$ be an increasing sequence of $\boratwo (\lceil T\rceil )$ sets with 
$\mathbb{S}\! =\! D\big( (\Sigma_\theta )_{\theta <\eta}\big)$, 
$K\! :=\! (f'\!\times\! g')\big[\lceil T\rceil\big]$, and 
${R_\theta\! :=\! (f'\!\times\! g')\big[\Sigma_\theta\big]}$. Note that $K$ is compact, $R_\theta$ is 
$K_\sigma$, $D\big( (R_\theta )_{\theta <\eta}\big)\!\subseteq\! A_1$, 
$K\!\setminus\! D\big( (R_\theta )_{\theta <\eta}\big)\!\subseteq\! A_0$, 
${D\big( (R_\theta )_{\theta <\eta}\big)\! =\! K\cap A_1}$,  
${K\!\setminus\! D\big( (R_\theta )_{\theta <\eta}\big)\! =\! K\cap A_0}$, so that 
$D\big( (R_\theta )_{\theta <\eta}\big)$ is not separable from 
$K\!\setminus\! D\big( (R_\theta )_{\theta <\eta}\big)$ by a 
$\mbox{pot}\big(\check D_\eta (\boratwo )\big)$ set. So we can replace $A_1,A_0$ with 
$D\big( (R_\theta )_{\theta <\eta}\big)$, $K\!\setminus\! D\big( (R_\theta )_{\theta <\eta}\big)$, respectively.

\vfill\eject

\noindent $\bullet$ We may also assume that $X$ is zero-dimensional and there are disjoint clopen subsets $O_0,O_1$ of $X$ such that $A_0\cap (O_0\!\times\! O_1)$ is not separable from $A_1\cap (O_0\!\times\! O_1)$ by a $\mbox{pot}\big( D_\eta (\boratwo )\big)$ set. So, without loss of generality, we will assume that $A_0\cup A_1\!\subseteq\! O_0\!\times\! O_1$. We may also assume that $X$ is recursively presented, $A_0,A_1,O_0,O_1,R_\theta$ are $\Borel$, and $R_\theta$ is the union of 
$\Borel\cap\bormone\!\subseteq\!\Ana\cap\bormone (\tau_1)\!\subseteq\!\boraone (\tau_2)$ sets.\bigskip

 We set, for $\theta\! <\!\eta$, 
$N_\theta\! :=\! R_\theta\!\setminus\! (\bigcup_{\theta'<\theta}~R_{\theta'})\cap
\bigcap_{\theta'<\theta}~\overline{N_{\theta'}}^{\tau_2}$. Note that the $N_\theta$'s are pairwise disjoint, which will be useful in the construction to get the injectivity of our reduction maps. We use the notation of Theorem \ref{kernel1}. For simplicity, we set 
$F^\varepsilon_\theta\! :=\! F^\varepsilon_{\theta ,2}$.\bigskip

\noindent\bf Claim.\it\ (a) Assume that $k\! +\! 1\! <\!\eta$. Then 
$F^\varepsilon_k\! =\!\overline{N_k}^{\tau_2}\cup E_k$, where $E_k\!\subseteq\!\neg R_{k+1}$ is 
$\tau_2$-closed.\smallskip

\noindent (b) $A_0\cap\bigcap_{\theta <\eta}~F^\varepsilon_\theta\! =\! N_\eta\! :=\! 
K\!\setminus\! (\bigcup_{\theta <\eta}~R_\theta )\cap
\bigcap_{\theta <\eta}~\overline{N_\theta}^{\tau_2}$.\rm\bigskip

\noindent (a) Indeed, we argue by induction on $k$ to prove (a). In the proof of this claim, all the closures will refer to $\tau_2$. Note first that 
$R_0\!\subseteq\! A_\varepsilon\!\subseteq\! R_0\cup\neg R_1$, so that 
$F^\varepsilon_0\! =\!\overline{A_\varepsilon}\! =\!\overline{R_0}\cup E_0\! =\!
\overline{N_0}\cup E_0$. Then, inductively, 
$$\begin{array}{ll}
F^\varepsilon_{k+1}
& \! =\!\overline{A_{1-\vert\mbox{parity}(k)-\varepsilon\vert}\cap F^\varepsilon_k}
\! =\!\overline{A_{1-\vert\mbox{parity}(k)-\varepsilon\vert}\cap (\overline{N_k}\cup E_k)}\cr\cr
& \! =\!\overline{\big( (R_{k+1}\!\setminus\! R_k)\cup (R_{k+3}\!\setminus\! R_{k+2})...\big)
\cap (\overline{N_k}\cup E_k)}\! =\!\overline{N_{k+1}}\cup E_{k+1}.
\end{array}$$ 
(b) Note then that $F^\varepsilon_{\eta -1}
\! =\!\overline{A_1\cap\bigcap_{k+1<\eta}~F^\varepsilon_k}
\! =\!\overline{A_1\cap\bigcap_{k+1<\eta}~(\overline{N_k}\cup E_k)}
\! =\!\overline{N_{\eta -1}}$, so that 
$$A_0\cap\bigcap_{\theta <\eta}~F^\varepsilon_\theta\! =\! 
K\!\setminus\! (\bigcup_{\theta <\eta}~R_\theta )\cap\bigcap_{\theta <\eta}~\overline{N_\theta}.$$ 
This proves the claim.\hfill{$\diamond$}\bigskip
 
\noindent $\bullet$ We construct the following objects:\bigskip

\noindent - sequences $(x_s)_{s\in 2^{<\omega}, 0\subseteq s}$, 
$(y_s)_{s\in 2^{<\omega}, 1\subseteq s}$ of points of $X$,\bigskip

\noindent - sequences $(X_s)_{s\in 2^{<\omega}, 0\subseteq s}$, 
$(Y_s)_{s\in 2^{<\omega}, 1\subseteq s}$ of $\Ana$ subsets of $X$,\bigskip

\noindent - a sequence $(U_{s,t})_{(s,t)\in T\setminus\{ (\emptyset ,\emptyset )\}}$ of $\Ana$ subsets of $X^2$.\bigskip

 We want these objects to satisfy the following conditions:
$$\begin{array}{ll}
& (1)\ x_s\!\in\! X_s\ \wedge\ y_s\!\in\! Y_s\ \wedge\ (x_s,y_t)\!\in\! U_{s,t}\cr
& (2)\ X_{s\varepsilon}\!\subseteq\! X_s\!\subseteq\!\Omega_X\cap O_0\ \wedge\ 
Y_{s\varepsilon}\!\subseteq\! Y_s\!\subseteq\!\Omega_X\cap O_1\ \wedge\ U_{s,t}\!\subseteq\!\Omega_{X^2}\cap (X_s\!\times\! Y_t)\cr
& (3)\ \mbox{diam}_{\mbox{GH}}(X_s)\mbox{, }\mbox{diam}_{\mbox{GH}}(Y_s)\mbox{, }
\mbox{diam}_{\mbox{GH}}(U_{s,t})\!\leq\! 2^{-\vert s\vert}\cr
& (4)\ X_{s0}\cap X_{s1}\! =\! Y_{s0}\cap Y_{s1}\! =\!\emptyset\cr
& (5)\ \big( ~(s,t)~R~(s',t')\ \wedge\ \exists\theta\!\leq\! 2~~(s,t),(s',t')\!\in\! D_\theta ~\big)\Rightarrow U_{s',t'}\!\subseteq\! U_{s,t}\cr
& (6)\ U_{s,t}\!\subseteq\! N_\theta\mbox{ if }(s,t)\!\in\! D_\theta\cr
& (7)\ (s,t)~R~(s',t')\Rightarrow U_{s',t'}\!\subseteq\!\overline{U_{s,t}}^{\tau_1}
\end{array}$$
$\bullet$ Assume that this has been done. As in the proof of Theorem \ref{checkD2Sigma01}, we get $f\! :\! 2^\omega\!\rightarrow\! X$ injective continuous. If $(\alpha ,\beta )\!\in\!\mathbb{N}_\eta$, then we can find $\theta\! <\!\eta$ of parity opposite to that of $\eta$ and $(n_k)_{k\in\omega}$ strictly increasing such that $(\alpha ,\beta )\vert n_k\!\in\! D_\theta$ and 
$(\alpha ,\beta )\vert n_k~R~(\alpha ,\beta )\vert n_{k+1}$ for each $k\!\in\!\omega$. In this case, by (1)-(3) and (5)-(6), $\big( U_{(\alpha ,\beta )\vert n_k}\big)_{k\in\omega}$ is a decreasing sequence of nonempty clopen subsets of $A_0\cap\Omega_{X^2}$ with vanishing diameters, so that its intersection is a singleton $\{ F(\alpha ,\beta )\}\!\subseteq\! A_0$. As 
$(x_{\alpha\vert n},y_{\beta\vert n})$ converges (for ${\it\Sigma}_{X^2}$ and thus for 
${\it\Sigma}_X^2$) to $F(\alpha ,\beta )$, 
$\big( f(\alpha ),f(\beta )\big)\! =\! F(\alpha ,\beta )\!\in\! A_0$. If 
$(\alpha ,\beta )\!\in\!\lceil T\rceil\!\setminus\!\mathbb{N}_\eta$, then we argue similarly to see 
that $\big( f(\alpha ),f(\beta )\big)\!\in\! A_1$.\bigskip

\noindent $\bullet$ So let us prove that the construction is possible. Let 
$(x_0,y_1)\!\in\! N_\eta\cap\Omega_{X^2}$, $X_0,Y_1$ be $\Ana$ subsets of $X$ with diameter at most $2^{-1}$ such that $x_0\!\in\! X_0\!\subseteq\!\Omega_X\cap O_0$ and 
$y_1\!\in\! Y_1\!\subseteq\!\Omega_X\cap O_1$, and $U_{0,1}$ be a $\Ana$ subset of $X^2$ with diameter at most $2^{-1}$ such that 
$(x_0,y_1)\!\in\! U_{0,1}\!\subseteq\! N_\eta\cap\Omega_{X^2}\cap (X_0\!\times\! Y_1)$. This completes the construction for $l\! =\! 1$ since $(0,1)\!\in\! D_\eta$.\bigskip

\noindent - Note that $(0^2,1^2)\!\in\! D_\eta$ since $m^\eta_{0,1}\! =\! 0$ and 
$m^\eta_{0^2,1^2}\! =\! 1$ if $\eta\!\geq\! 2$. We set 
$S_0\! :=\!\overline{U_{0,1}}^{\tau_1}\cap (X_0\!\times\! Y_1)$ and 
$S_1\! :=\! S_0\cap N_0\cap\Omega_{X^2}$. As $U_{0,1}\!\subseteq\!\overline{N_0}^{\tau_2}$, 
$S_0\!\subseteq\!\overline{S_1}^{\tau_1}$. In particular, $\Pi_\varepsilon [S_1]$ is ${\it\Sigma}_X$-dense in $\Pi_\varepsilon [S_0]$ for each $\varepsilon\!\in\! 2$, by continuity of the projections. As $(x_0,y_1)\!\in\! U_{0,1}\cap (\Pi_0[S_0]\!\times\!\Pi_1[S_0])$, this implies that 
$U_{0,1}\cap (\Pi_0[S_1]\!\times\!\Pi_1[S_1])$ is not empty and contains some $(x_{0^2},y_{1^2})$ (the projections maps are open). This gives $y_{10}\!\in\! X$ with $(x_{0^2},y_{10})\!\in\! S_1$, and 
$x_{01}\!\in\! X$ with $(x_{01},y_{1^2})\!\in\! S_1$. As $U_{0,1}\!\subseteq\! N_\eta$ and 
$S_1\!\subseteq\! N_0$, $x_{0^2}\!\not=\! x_{01}$ and $y_{10}\!\not=\! y_{1^2}$. It remains to choose $\Ana$ subsets $X_{0^2},X_{01},Y_{10},Y_{1^2}$ of $X$ with diameter at most $2^{-2}$ such that $(x_{0\varepsilon},y_{1\varepsilon})\!\in\! X_{0\varepsilon}\!\times\! Y_{1\varepsilon}
\!\subseteq\! X_0\!\times\! Y_1$ and $X_{0^2}\cap X_{01}\! =\! Y_{10}\cap Y_{1^2}\! =\!\emptyset$, as well as $\Ana$ subsets $U_{0^2,1^2},U_{0^2,10},U_{01,1^2}$ of $X^2$ with diameter at most $2^{-2}$ such that 
$(x_{0^2},y_{1^2})\!\in\! U_{0^2,1^2}\!\subseteq\! U_{0,1}\cap (X_{0^2}\!\times\! Y_{1^2})$ and 
$(x_{0\varepsilon},y_{1\varepsilon})\!\in\! U_{0\varepsilon ,1\varepsilon}\!\subseteq\!
\overline{U_{0,1}}^{\tau_1}\cap N_0\cap\Omega_{X^2}\cap  (X_{0\varepsilon}\!\times\! Y_{1\varepsilon})$. This completes the construction for $l\! =\! 2$.\bigskip 

\noindent - Assume that our objects are constructed for the level $l\!\geq\! 2$, which is the case for 
$l\! =\! 2$. Note that $(s_l0,t_l1)\!\notin\! D_0$, and we already noticed that 
$s_l(l\! -\! 1)\! =\! t_l(l\! -\! 1)$ since $l\!\geq\! 2$, so that $(s_l,t_l)\!\in\! D_0$. We set 
$(\tilde s,\tilde t)\! :=\! (s_{l-1}0,t_{l-1}1)$ (which is not in $D_0$), and\bigskip 

\leftline{$S_0\! :=\!\big\{
\big( (\overline{x}_s)_{s\in 2^l, 0\subseteq s},(\overline{y}_t)_{t\in 2^l,1\subseteq t}\big)\!\in\! X^{2^l}\mid
\forall (s,t)\!\in\! T\!\cap\! (2^l\!\times\! 2^l)\!\setminus\!\{ (\tilde s,\tilde t)\}~~
(\overline{x}_{s},\overline{y}_{t})\!\in\! U_{s,t}\ \wedge$}\smallskip
  
\rightline{$(\overline{x}_{\tilde s},\overline{y}_{\tilde t})\!\in\!\overline{N_0}^{\tau_2}\cap
\overline{U_{\tilde s,\tilde t}}^{\tau_1}\cap (X_{\tilde s}\!\times\! Y_{\tilde t})\big\}\mbox{,}$}\bigskip

\leftline{$S_1\! :=\!\big\{
\big( (\overline{x}_s)_{s\in 2^l, 0\subseteq s},(\overline{y}_t)_{t\in 2^l,1\subseteq t}\big)\!\in\! S_0\mid
(\overline{x}_{\tilde s},\overline{y}_{\tilde t})\!\in\! N_0\cap\Omega_{X^2}\big\} .$}\bigskip 

\noindent We equip $X^{2^l}$ with the product of the Gandy-Harrington topologies. Let us show that $S_1$ is dense in $S_0$. Let $({\cal U}_s)_{s\in 2^l,0\subseteq s}$ and 
$({\cal V}_t)_{t\in 2^l,1\subseteq t}$ be sequences of $\Ana$ sets with 
$$\big( (\Pi_{s\in 2^l,0\subseteq s}~{\cal U}_s)\!\times\! (\Pi_{t\in 2^l,1\subseteq t}~{\cal V}_t)\big)\cap S_0
\!\not=\!\emptyset$$ 
with witness $\big( (x'_s),(y'_t)\big)$, ${\cal A}_\varepsilon\! :=\{ s\!\in\! 2^l\mid s(l\! -\! 1)\! =\!\varepsilon\}$, and\bigskip  

\leftline{$U\! :=\!\{\overline{x}_{\tilde s}\!\in\! {\cal U}_{\tilde s}\mid
\exists (\overline{x}_s)_{s\in {\cal A}_0\!\setminus\!\{\tilde s\}}\!\in\! 
\Pi_{s\in {\cal A}_0\!\setminus\!\{\tilde s\}}~{\cal U}_s~~
\exists (\overline{y}_t)_{t\in {\cal A}_0}\!\in\!\Pi_{t\in {\cal A}_0}~{\cal V}_t$}\smallskip  

\rightline{$\forall (s,t)\!\in\! T\cap ({\cal A}_0\!\times\! {\cal A}_0)~~
(\overline{x}_s,\overline{y}_t)\!\in\! U_{s,t}\}\mbox{,}$}\bigskip

\leftline{$V\! :=\!\{\overline{y}_{\tilde t}\!\in\! {\cal V}_{\tilde t}\mid
\exists (\overline{x}_s)_{s\in {\cal A}_1}\!\in\!\Pi_{s\in {\cal A}_1}~{\cal U}_s~~
\exists (\overline{y}_t)_{t\in {\cal A}_1\!\setminus\!\{\tilde t\}}\!\in\!
\Pi_{t\in {\cal A}_1\!\setminus\!\{\tilde t\}}~{\cal V}_t$}\smallskip

\rightline{$\forall (s,t)\!\in\! T\cap ({\cal A}_1\!\times\! 
{\cal A}_1)~~(\overline{x}_s,\overline{y}_t)\!\in\! U_{s,t}\} .$}\bigskip

 Then $(x'_{\tilde s},y'_{\tilde t})\!\in\!
\overline{N_0}^{\tau_2}\cap\overline{U_{\tilde s,\tilde t}}^{\tau_1}\cap (U\times V)$. This gives 
$(\overline{x}_{\tilde s},\overline{y}_{\tilde t})$ in 
${N_0\cap\overline{U_{\tilde s,\tilde t}}^{\tau_1}\cap (U\times V)\cap\Omega_{X^2}}$. We choose witnesses 
$(\overline{x}_s)_{s\in {\cal A}_0\setminus\{\tilde s\}}$, 
$(\overline{y}_t)_{t\in {\cal A}_0}$ (resp., $(\overline{x}_s)_{s\in {\cal A}_1}$, 
$(\overline{y}_t)_{t\in {\cal A}_1\setminus\{\tilde t\}}$) for the 
fact that $\overline{x}_{\tilde s}\in U$ (resp., $\overline{y}_{\tilde t}\in V$). 
Then $\big((\overline{x}_s),(\overline{y}_t)\big)\!\in\! 
\big( (\Pi_{s\in 2^l,0\subseteq s}~{\cal U}_t)\times (\Pi_{t\in 2^l,1\subseteq t}~{\cal V}_t)\big)\cap S_1$, as desired.\bigskip 

 The sets $U_\varepsilon\! :=\!\Pi_{s_l}[S_\varepsilon ]$ and $V_\varepsilon\! :=\!\Pi_{t_l}[S_\varepsilon ]$ 
are $\Ana$ sets. As $S_1$ is dense in $S_0$, $U_1$ (resp., $V_1$) is dense in $U_0$ (resp., $V_0$).  Note that ${(x_{s_l},y_{t_l})\in U_{s_l,t_l}\cap (U_0\!\times\! V_0)}$. As $U_1$ (resp,. $V_1$) is dense in $U_0$ (resp., $V_0$), $U_{s_l,t_l}$ meets $U_1\!\times\! V_1$.\bigskip

 Let $(s_l0,t_l1)^R$ be the $R$-predecessor of $(s_l0,t_l1)$. Assume first that 
$(s_l0,t_l1)\!\in\! D_\eta$. Then $(s_l0,t_l1)^R\!\in\! D_\eta$ too. Note that 
${U_{s_l,t_l}\!\subseteq\!\overline{U_{(s_l0,t_l1)^R}}^{\tau_1}}$ since $(s_l0,t_l1)^R~R~(s_l,t_l)$. Thus ${\overline{U_{(s_l0,t_l1)^R}}^{\tau_1}}$ meets $U_1\!\times\! V_1$. This gives 
$(x_{s_l0},y_{t_l1})\!\in\! U_{(s_l0,t_l1)^R}\cap (U_1\!\times\! V_1)$.  
We choose witnesses $(x_{s0})_{s\in 2^l\setminus\{s_l\} ,0\subseteq s}$, 
$(y_{t0})_{t\in 2^l,1\subseteq t}$ (resp., $(x_{s1})_{s\in 2^l,0\subseteq s}$, 
$(y_{t1})_{t\in 2^l\setminus\{t_l\} ,1\subseteq t}$) for the fact that $x_{s_l0}\!\in\! U_1$ (resp., 
$y_{t_l1}\!\in\! V_1$). As $(x_{s_l0},y_{t_l1})\!\in\! U_{(s_l0,t_l1)^R}\!\subseteq\! N_\eta$ and 
$(x_{s_l\varepsilon},y_{t_l\varepsilon})\!\in\! N_0$, ${x_{s_l0}\!\not=\! x_{s_l1}}$ and 
$y_{t_l0}\!\not=\! y_{t_l1}$. As in the proof of Theorem \ref{checkD2Sigma01}, the s-acyclicity of 
$A_0\cup A_1$ and the fact that $O_0,O_1$ are disjoint ensure the fact that 
$x_{s0}\!\not=\! x_{s1}$ and $y_{t0}\!\not=\! y_{t1}$ for $s,t$ arbitrary with the right first coordinate. Then we choose $\Ana$ subsets $X_{s\varepsilon},Y_{t\varepsilon}$ of $X$ with diameter at most $2^{-l-1}$ such that 
$(x_{s\varepsilon},y_{t\varepsilon})\!\in\! X_{s\varepsilon}\!\times\! Y_{t\varepsilon}\!\subseteq\! X_s\!\times\! Y_t$ and $X_{s0}\cap X_{s1}\! =\! Y_{s0}\cap Y_{s1}\! =\!\emptyset$, as well as $\Ana$ subsets 
$U_{s\varepsilon  ,t\varepsilon '}$ of $X^2$, with diameter at most $2^{-l-1}$, containing 
$(x_{s\varepsilon },y_{t\varepsilon '})$ and contained in 
$X_{s\varepsilon}\!\times\! Y_{t\varepsilon}$, such that\bigskip

\noindent - $U_{s_l0,t_l1}\!\subseteq\! U_{(s_l0,t_l1)^R}$,\smallskip

\noindent - $U_{\tilde s\varepsilon ,\tilde t\varepsilon}\!\subseteq\!
\overline{U_{\tilde s,\tilde t}}^{\tau_1}\cap N_0\cap\Omega_{X^2}$,\smallskip

\noindent - $U_{s\varepsilon ,t\varepsilon}\!\subseteq\! U_{s,t}$ if $(s,t)\!\not=\! (\tilde s,\tilde t)$.\bigskip
 
  The argument is the same if $(s_l0,t_l1),(s_l0,t_l1)^R\!\in\! D_\theta$. So it remains to study the case where 
$(s_l0,t_l1)\!\in\! D_{\theta'}$ and $(s_l0,t_l1)^R\!\in\! D_\theta$, and $\theta'\! <\!\theta$. In this case, note that $U_{(s_l0,t_l1)^R}\cap (U_1\!\times\! V_1)$ is not empty and contained in 
$N_\theta\!\subseteq\!\overline{N_{\theta'}}^{\tau_2}$. This gives 
$(x_{s_l0},y_{t_l1})\!\in\! N_{\theta'}\cap\overline{U_{(s_l0,t_l1)^R}}^{\tau_1}\cap\Omega_{X^2}\cap (U_1\!\times\! V_1)$, and we conclude as before.
\hfill{$\square$}\bigskip

\noindent $\underline{\mbox{\bf Consequences}}$\bigskip

\begin{cor} \label{checkDnSigma02uni} Let $1\!\leq\!\eta\! <\!\omega$, $X$ be a Polish space, and $A,B$ be disjoint analytic relations on $X$ such that $A$ is contained in a $\mbox{pot}(\bortwo )$ s-acyclic relation. Then exactly one of the following holds:\smallskip  

(a) the set $A$ is separable from $B$ by a $\mbox{pot}\big( D_\eta (\boratwo )\big)$ set,\smallskip  

(b) $(2^\omega ,2^\omega ,\mathbb{N}_\eta ,\lceil T\rceil\!\setminus\!\mathbb{N}_\eta )\sqsubseteq (X,X,A,B)$, via a square map.\end{cor}

\noindent\bf Proof.\rm ~Let $R$ be a $\mbox{pot}(\bortwo )$ s-acyclic relation containing $A$. By Lemma \ref{ExDnS2}, (a) and (b) cannot hold simultaneously. So assume that (a) does not hold. Then $A$ is not separable from $B\cap R$ by a $\mbox{pot}\big( D_\eta (\boratwo )\big)$ set. This allows us to apply Theorem \ref{checkDnSigma02}.\hfill{$\square$}

\begin{cor} \label{caracpartialcheckD2Sigma02} Let $1\!\leq\!\eta\! <\!\omega$, $X$ be a Polish space, and $A,B$ be disjoint analytic relations on $X$. The following are equivalent:\smallskip

\noindent (1) there is $R\!\in\!\ana$ s-acyclic such that $A\cap R$ is not separable from $B\cap R$ by a $\mbox{pot}\big( D_\eta (\boratwo )\big)$ set,\smallskip

\noindent (2) there is $f\! :\! 2^\omega\!\rightarrow\! X$ injective continuous such that $\mathbb{N}_\eta\!\subseteq\! (f\!\times\! f)^{-1}(A)$ and 
$\lceil T\rceil\!\setminus\!\mathbb{N}_\eta\!\subseteq\! (f\!\times\! f)^{-1}(B)$.\end{cor} 

\noindent\bf Proof.\rm ~(1) $\Rightarrow$ (2) We apply Theorem \ref{checkDnSigma02}.\bigskip

\noindent (2) $\Rightarrow$ (1) We can take $R\! :=\! (f\!\times\! f)\big[\lceil T\rceil\big]$.\hfill{$\square$}

\section{$\!\!\!\!\!\!$ Oriented graphs}

\noindent\bf Proof of Theorem \ref{Zog}.\rm ~Theorem \ref{motivating} provides Borel relations 
$\mathbb{S}_0$, $\mathbb{S}_1$ on $2^\omega$. We saw that $\mathbb{S}_0\cup\mathbb{S}_1$ is a subset of the body of a tree $T$, which does not depend on $\bf\Gamma$, and is contained in 
$N_0\!\times\! N_1$. We set 
$\mathbb{G}_{\bf\Gamma}\! :=\!\mathbb{S}_0\cup (\mathbb{S}_1)^{-1}$, so that 
$\mathbb{G}_{\bf\Gamma}$ is Borel. As 
$\mathbb{S}_0\cup\mathbb{S}_1\!\subseteq\! N_0\!\times\! N_1$ and $\mathbb{S}_0$, 
$\mathbb{S}_1$ are disjoint, $\mathbb{G}_{\bf\Gamma}$ is an oriented graph. If (a) and (b) hold, then $\mathbb{G}_{\bf\Gamma}$ is separable from $\mathbb{G}_{\bf\Gamma}^{-1}$ by a 
$\mbox{pot}({\bf\Gamma})$ set $S$. Note that $S$ also separates 
$\mathbb{S}_0\! =\!\mathbb{G}_{\bf\Gamma}\cap (N_0\!\times\! N_1)$ from 
$\mathbb{S}_1\! =\!\mathbb{G}_{\bf\Gamma}^{-1}\cap (N_0\!\times\! N_1)$, which is absurd. Thus (a) and (b) cannot hold simultaneously.\bigskip

 Assume now that (a) does not hold. Then there are $g,h\! :\! 2^\omega\!\rightarrow\! X$ continuous such that the inclusions $\mathbb{S}_0\!\subseteq\! (g\!\times\! h)^{-1}(G)$ and 
$\mathbb{S}_1\!\subseteq\! (g\!\times\! h)^{-1}(G^{-1})$ hold. It remains to set 
$f(0\alpha )\! :=\! g(0\alpha )$ and $f(1\beta )\! :=\! h(1\beta )$.\hfill{$\square$}\bigskip

\noindent\bf Proof of Theorem \ref{Zog2}.\rm ~We argue as in the proof of Theorem \ref{Zog}. The things to note are the following:\bigskip

- if $G$ is s-acyclic or locally countable, then $s(G)$ too,\smallskip

- as noted in [Lo4], if $G$ is separable from $G^{-1}$ by a $\mbox{pot}({\bf\Gamma})$ set $S$, then $S^{-1}\!\in\!\mbox{pot}({\bf\Gamma})$ separates $G^{-1}$ from $G$, and $\neg S^{-1}\!\in\!\mbox{pot}(\check {\bf\Gamma})$ separates $G$ from $G^{-1}$, so that we can restrict our attention to the classes 
$D_\eta (\boraxi )$ and $\bortwo$.\bigskip

\noindent $\bullet$ If $\bf\Gamma$ has rank two, then Theorem \ref{checkDnSigma02} and Corollary \ref{cor1delta02} provide Borel relations $\mathbb{S}_0$, $\mathbb{S}_1$ on 
$2^\omega$.\bigskip
 
\noindent $\bullet$ If ${\bf\Gamma}\! =\! D_\eta (\boraone )$, then Corollaries 
\ref{cor1checkD2Sigma01} and \ref{caracpartialcheckD2Sigma01} provide 
$f\! :\! 2^\omega\!\rightarrow\! X$ injective continuous such that one of the following holds:\smallskip

(a) $\mathbb{N}^\eta_0\!\subseteq\! (f\!\times\! f)^{-1}(G)$ and 
$\mathbb{N}^\eta_1\!\subseteq\! (f\!\times\! f)^{-1}(G^{-1})$,\smallskip

(b) $\mathbb{B}^\eta_0\!\subseteq\! (f\!\times\! f)^{-1}(G)$ and 
$\mathbb{B}^\eta_1\!\subseteq\! (f\!\times\! f)^{-1}(G^{-1})$.\smallskip

\noindent The case (a) cannot happen since $G^{-1}$ is irreflexive.\hfill{$\square$}\bigskip

\noindent\bf Proof of Theorem \ref{Zog3}.\rm ~Note first that 
$\mathbb{S}^\eta_0\cup (\mathbb{S}^\eta_1)^{-1},\mathbb{C}^\eta_0\cup (\mathbb{C}^\eta_1)^{-1},\mathbb{B}^\eta_0\cup (\mathbb{B}^\eta_1)^{-1}$ and 
$\mathbb{B}^\eta_1\cup (\mathbb{B}^\eta_0)^{-1}$ are Borel oriented graphs with locally countable closure. As in the proof of Theorem \ref{Zog}, $\mathbb{G}$ is not separable from 
$\mathbb{G}^{-1}$ by a $\mbox{pot}\Big(\Delta\big( D_\eta (\boraone )\big)\Big)$ set if 
$\mathbb{G}\!\in\!\{\mathbb{C}^\eta_0\cup (\mathbb{C}^\eta_1)^{-1},
\mathbb{B}^\eta_0\cup (\mathbb{B}^\eta_1)^{-1},
\mathbb{B}^\eta_1\cup (\mathbb{B}^\eta_0)^{-1}\}$. By Lemma \ref{S^2}, 
$\mathbb{S}^\eta_0\cup (\mathbb{S}^\eta_1)^{-1}$ is not separable from 
$(\mathbb{S}^\eta_0)^{-1}\cup\mathbb{S}^\eta_1$ by a 
$\mbox{pot}\Big(\Delta\big( D_\eta (\boraone )\big)\Big)$ set.\bigskip

\noindent $\bullet$ Assume now that (a) does not hold. Corollaries \ref{cor1delta01} and \ref{caracpartialdelta01} provide 
$$(\mathbb{A},\mathbb{B})\!\in\!\{ (\mathbb{N}^\eta_1,\mathbb{N}^\eta_0),
(\mathbb{B}^\eta_1,\mathbb{B}^\eta_0),(\mathbb{N}^\eta_0,\mathbb{N}^\eta_1),
(\mathbb{B}^\eta_0,\mathbb{B}^\eta_1),(\mathbb{S}^\eta_0,\mathbb{S}^\eta_1),
(\mathbb{C}^\eta_0,\mathbb{C}^\eta_1)\}$$ 
and $f\! :\! 2^\omega\!\rightarrow\! X$ injective continuous such that 
$\mathbb{A}\!\subseteq\! (f\!\times\! f)^{-1}(G)$ and 
$\mathbb{B}\!\subseteq\! (f\!\times\! f)^{-1}(G^{-1})$.\bigskip

 The pair $(\mathbb{A},\mathbb{B})$ cannot be in 
$\{ (\mathbb{N}^\eta_1,\mathbb{N}^\eta_0),(\mathbb{N}^\eta_0,\mathbb{N}^\eta_1)\}$ since $G$ and $G^{-1}$ are irreflexive. It is enough to show the existence of $f\! :\! 2^\omega\!\rightarrow\! 2^\omega$ injective continuous such that $\mathbb{B}^\eta_0\cup (\mathbb{B}^\eta_1)^{-1}\!\subseteq\! 
(f\!\times\! f)^{-1}(\mathbb{B}^\eta_1\cup (\mathbb{B}^\eta_0)^{-1})$ to see that (b) holds.\bigskip

- We use the notation of the proof of Proposition \ref{squareDeltaD2}. Let us show that 
$$F^{\mbox{parity}(\eta )}_\theta\! :=\! F^{\mbox{parity}(\eta )}_{\theta ,1}\!\subseteq\! C_\theta$$ 
if $\theta\! <\!\eta$ (where $A_\varepsilon\! =\!\mathbb{N}^\eta_\varepsilon$ and the closures refer to $\tau_1$). We argue by induction on $\theta$. Note first that 
$$F^{\mbox{parity}(\eta )}_0\! =\!\overline{\mathbb{N}^\eta_{\mbox{parity}(\eta )}}\! =\!
\overline{\bigcup_{\mbox{parity}(\varphi (s))=0}~\mbox{Gr}(f_s)}\!\subseteq\!
\overline{C_0}\! =\! C_0\mbox{,}$$ 
by the proof of Proposition \ref{squareDeltaD2}. Then, inductively, 
$$\begin{array}{ll}
F^{\mbox{parity}(\eta )}_\theta\!\!\!\!
& \! =\!\overline{\mathbb{N}^\eta_{\vert\mbox{parity}(\theta )-\mbox{parity}(\eta )\vert}\cap
\bigcap_{\theta'<\theta}~F^{\mbox{parity}(\eta )}_{\theta'}}\cr
& \!\subseteq\!
\overline{\bigcup_{\mbox{parity}(\varphi (s))=\mbox{parity}(\theta )}~
\mbox{Gr}(f_s)\cap\bigcap_{\theta'<\theta}~
\bigcup_{\varphi (s)\geq\theta'}~\mbox{Gr}(f_s)}\! =\!\overline{C_\theta}\! =\! C_\theta\mbox{,}
\end{array}$$
by the proof of Proposition \ref{squareDeltaD2}.\bigskip
  
- From this we deduce that 
$\mathbb{N}^\eta_0\cap\bigcap_{\theta <\eta}~F^{\mbox{parity}(\eta )}_\theta$ is contained in 
$$\big(\bigcup_{\mbox{parity}(\varphi (s))=\mbox{parity}(\eta )}~
\mbox{Gr}(f_s)\big)\cap\bigcap_{\theta <\eta}~C_\theta\!\subseteq\!\mbox{Gr}(f_\emptyset )\! =\!
\Delta (2^\omega ).$$ 
As $\mathbb{N}^\eta_0\cup\mathbb{N}^\eta_1$ is locally countable and 
$\mathbb{N}^\eta_0\cap\bigcap_{\theta <\eta}~F^{\mbox{parity}(\eta )}_\theta\!\subseteq\!
\Delta (2^\omega )$, the proof of Theorem \ref{checkD2Sigma01} gives 
$h\! :\! 2^\omega\!\rightarrow\! 2^\omega$ injective continuous such that 
$\mathbb{N}^\eta_0\!\subseteq\! (h\!\times\! h)^{-1}\big( (\mathbb{N}^\eta_0)^{-1}\big)$ and 
$\mathbb{N}^\eta_1\!\subseteq\! (h\!\times\! h)^{-1}\big( (\mathbb{N}^\eta_1)^{-1}\big)$ (we are in the case 2 of this proof). The map 
$f\! :\!\varepsilon\alpha\!\mapsto\! (1\! -\!\varepsilon )h(\alpha )$ is as desired.\bigskip

\noindent $\bullet$ As $\Delta (2^\omega )$ is contained in the closure of 
$\mathbb{S}^\eta_0\cup (\mathbb{S}^\eta_1)^{-1}$, this last relation is not below the two others.\bigskip

- Assume, towards a contradiction, that $\mathbb{B}^\eta_0\cup (\mathbb{B}^\eta_1)^{-1}$ is below $\mathbb{S}^\eta_0\cup (\mathbb{S}^\eta_1)^{-1}$. This gives $s\!\in\! 2^{<\omega}$ and 
$\varepsilon\!\in\! 2$ such that $\big( N_{0s},N_{1s},
\mathbb{B}^\eta_0\cap (N_{0s}\!\times\! N_{1s}),
\mathbb{B}^\eta_1\cap (N_{0s}\!\times\! N_{1s})\big)\sqsubseteq 
\big( 2^\omega ,2^\omega ,(\mathbb{S}^\eta_\varepsilon )^{1-2\varepsilon},
(\mathbb{S}^\eta_{1-\varepsilon})^{1-2\varepsilon}\big)$. By Lemma \ref{S^2}, 
$\mathbb{N}^\eta_0\cap N_s^2$ is not separable from $\mathbb{N}^\eta_1\cap N_s^2$ by a 
$\mbox{pot}\big( D_\eta (\boraone )\big)$ set. As $\mathbb{N}^\eta_0\cup\mathbb{N}^\eta_1$ is locally countable and $\mathbb{N}^\eta_0\cap\bigcap_{\theta <\eta}~F^{\mbox{parity}(\eta )}_\theta\!\subseteq\!
\Delta (2^\omega )$, the proof of Theorem \ref{checkD2Sigma01} gives 
$h\! :\! 2^\omega\!\rightarrow\! N_s$ injective continuous such that 
$\mathbb{N}^\eta_\epsilon\!\subseteq\! (h\!\times\! h)^{-1}(\mathbb{N}^\eta_\epsilon\cap N_s^2)$ for each $\epsilon\!\in\! 2$ (we are in the case 2 of this proof). This implies that 
$(2^\omega ,2^\omega ,\mathbb{B}^\eta_0,\mathbb{B}^\eta_1)\sqsubseteq
\big( N_{0s},N_{1s},\mathbb{B}^\eta_0\cap (N_{0s}\!\times\! N_{1s}),
\mathbb{B}^\eta_1\cap (N_{0s}\!\times\! N_{1s})\big)$ and 
$$(2^\omega ,2^\omega ,\mathbb{B}^\eta_0,\mathbb{B}^\eta_1)\sqsubseteq 
\big( 2^\omega ,2^\omega ,(\mathbb{S}^\eta_\varepsilon )^{1-2\varepsilon},
(\mathbb{S}^\eta_{1-\varepsilon})^{1-2\varepsilon}\big) .$$ 
By Corollary \ref{caracpartialcheckD2Sigma01}, 
$(2^\omega ,2^\omega ,\mathbb{N}^\eta_0,\mathbb{N}^\eta_1)\sqsubseteq 
(2^\omega ,2^\omega ,\mathbb{B}^\eta_0,\mathbb{B}^\eta_1)$, so that 
$$(2^\omega ,2^\omega ,\mathbb{N}^\eta_0,\mathbb{N}^\eta_1)\sqsubseteq 
\big( 2^\omega ,2^\omega ,(\mathbb{S}^\eta_\varepsilon )^{1-2\varepsilon},
(\mathbb{S}^\eta_{1-\varepsilon})^{1-2\varepsilon}\big) .$$ 
But this contradicts the proof of Proposition \ref{squareDeltaD2}.

\vfill\eject

- We will show that $(2^\omega ,2^\omega ,\mathbb{C}^\eta_0,\mathbb{C}^\eta_1)\sqsubseteq 
(2^\omega ,2^\omega ,\mathbb{S}^\eta_0,\mathbb{S}^\eta_1)$. Using the proof of the previous point, this will show that $\mathbb{B}^\eta_0\cup (\mathbb{B}^\eta_1)^{-1}$ is not below 
$\mathbb{C}^\eta_0\cup (\mathbb{C}^\eta_1)^{-1}$.\bigskip

 We use the notation of the proof of Proposition \ref{squareDeltaD2}. Let us show that 
$G_\theta\! :=\! G_{\theta ,1}\!\subseteq\! C_\theta$ if $1\!\leq\!\theta\!\leq\!\eta$ (where 
$A_\varepsilon\! =\!\mathbb{S}^\eta_\varepsilon$ and the closures refer to $\tau_1$). We argue by induction on $\theta$. Note first that 
$$G_1\! =\!\overline{\mathbb{S}^\eta_0}\cap\overline{\mathbb{S}^\eta_1}\! =\!
\overline{U^0_0}\cap\overline{U^1_0}\! =\! C^0_1\cup C^1_1\! =\! C_1$$
by the proof of Proposition \ref{squareDeltaD2}. Then, inductively, 
$$G_{\theta +1}\! =\!\overline{\mathbb{S}^\eta_0\cap G_\theta}\cap
\overline{\mathbb{S}^\eta_1\cap G_\theta}\!\subseteq\!
\overline{U^0_0\cap C_\theta}\cap\overline{U^1_0\cap C_\theta}\!\subseteq\! C_{\theta +1}$$
and $G_\lambda\! =\!\bigcap_{\theta <\lambda}~G_\theta\!\subseteq\!
\bigcap_{\theta <\lambda}~C_\theta\! =\! C_\lambda$ if $\lambda$ is limit.\bigskip
  
 From this we deduce that 
$G_\eta\!\subseteq\! C_\eta\! =\!\mbox{Gr}(f_\emptyset )\! =\!\Delta (2^\omega )$.  
As $\mathbb{S}^\eta_0\cup\mathbb{S}^\eta_1$ is locally countable and 
$G_\eta\!\subseteq\!\Delta (2^\omega )$, the proof of Theorem \ref{exdeltaeta1} gives 
$h\! :\! 2^\omega\!\rightarrow\! N_s$ injective continuous such that the inclusion 
$\mathbb{S}^\eta_\epsilon\!\subseteq\! (h\!\times\! h)^{-1}(\mathbb{S}^\eta_\epsilon\cap N_0^2)$ holds for each $\epsilon\!\in\! 2$ (we are in the case 2 of this proof). The maps defined by 
$f(0\alpha )\! :=\! h(\alpha )$, $f (1\alpha )\! :=\! 1\alpha$, $g(1\beta )\! :=\! h(\beta )$ and 
$g (0\beta )\! :=\! 1\beta$,  are as desired.\bigskip

- Assume, towards a contradiction, that $\mathbb{C}^\eta_0\cup (\mathbb{C}^\eta_1)^{-1}$ is below $\mathbb{S}^\eta_0\cup (\mathbb{S}^\eta_1)^{-1}$, with witness $f$. This gives 
$s\!\in\! 2^{<\omega}\!\setminus\!\{\emptyset\}$ and $\varepsilon\!\in\! 2$ such that 
$\mathbb{C}^\eta_\epsilon\cap (N_{0s}\!\times\! N_{1s})\!\subseteq\! 
(f\!\times\! f)^{-1}\big( (\mathbb{S}^\eta_{\vert\epsilon -\varepsilon\vert})^{1-2\varepsilon}\big)$ for each $\epsilon\!\in\! 2$. As in the previous point, there is 
$h\! :\! 2^\omega\!\rightarrow\! N_s$ injective continuous such that 
$$\mathbb{S}^\eta_\epsilon\!\subseteq\! (h\!\times\! h)^{-1}(\mathbb{S}^\eta_\epsilon\cap N_s^2)$$ for each $\epsilon\!\in\! 2$. This implies that if we set $k(\epsilon\alpha )\! :=\!\epsilon h(\alpha )$ and $l\! :=\! f\circ k$, then 
$$\mathbb{C}^\eta_\epsilon\!\subseteq\! 
(k\!\times\! k)^{-1}\big(\mathbb{C}^\eta_\epsilon\cap (N_{0s}\!\times\! N_{1s})\big)$$ 
and $\mathbb{C}^\eta_\epsilon\!\subseteq\! 
(l\!\times\! l)^{-1}\big( (\mathbb{S}^\eta_{\vert\epsilon -\varepsilon\vert})^{1-2\varepsilon}\big)$. As in the proof of Proposition \ref{squareDeltaD2}, we see that the image of 
$$\{ (0\alpha ,1\alpha )\mid\alpha\!\in\! 2^\omega\}$$ 
by $l\!\times\! l$ is contained in the diagonal of $2^\omega$, which is not possible by injectivity of $l$.\bigskip

- Assume that $\eta$ is a successor ordinal. The previous points show that if 
$\mathbb{C}^\eta_0\cup (\mathbb{C}^\eta_1)^{-1}$ is below 
$\mathbb{B}^\eta_0\cup (\mathbb{B}^\eta_1)^{-1}$, then 
$(2^\omega ,2^\omega ,\mathbb{C}^\eta_0,\mathbb{C}^\eta_1)\sqsubseteq 
\big( 2^\omega ,2^\omega ,(\mathbb{B}^\eta_\varepsilon )^{1-2\varepsilon},
(\mathbb{B}^\eta_{1-\varepsilon})^{1-2\varepsilon}\big)$ for some 
$\varepsilon\!\in\! 2$. We saw that there is $h\! :\! 2^\omega\!\rightarrow\! N_0$ injective continuous such that $\mathbb{N}^\eta_\epsilon\!\subseteq\! 
(h\!\times\! h)^{-1}(\mathbb{N}^\eta_\epsilon\cap N_0^2)$ for each $\epsilon\!\in\! 2$. The maps defined by $f(0\alpha )\! :=\! h(\alpha )$, $f (1\alpha )\! :=\! 1\alpha$, $g(1\beta )\! :=\! h(\beta )$ and 
$g (0\beta )\! :=\! 1\beta$ are witnesses for the fact that 
$(2^\omega ,2^\omega ,\mathbb{B}^\eta_0,\mathbb{B}^\eta_1)\sqsubseteq 
(2^\omega ,2^\omega ,\mathbb{N}^\eta_0,\mathbb{N}^\eta_1)$, so that 
$(2^\omega ,2^\omega ,\mathbb{C}^\eta_0,\mathbb{C}^\eta_1)\sqsubseteq 
\big( 2^\omega ,2^\omega ,(\mathbb{N}^\eta_\varepsilon )^{1-2\varepsilon},
(\mathbb{N}^\eta_{1-\varepsilon})^{1-2\varepsilon}\big)$. The maps $\alpha\!\mapsto\! 0\alpha$ and $\beta\!\mapsto\! 1\beta$ are witnesses for the fact that 
$(2^\omega ,2^\omega ,\mathbb{S}^\eta_0,\mathbb{S}^\eta_1)\sqsubseteq 
(2^\omega ,2^\omega ,\mathbb{C}^\eta_0,\mathbb{C}^\eta_1)$. Thus 
$(2^\omega ,2^\omega ,\mathbb{S}^\eta_0,\mathbb{S}^\eta_1)\sqsubseteq 
\big( 2^\omega ,2^\omega ,(\mathbb{N}^\eta_\varepsilon )^{1-2\varepsilon},
(\mathbb{N}^\eta_{1-\varepsilon})^{1-2\varepsilon}\big)$, which contradicts the proof of 
Proposition \ref{squareDeltaD2}.\bigskip

- Assume that $\eta$ is a limit ordinal. Let us show that 
$\mathbb{C}^\eta_0\cup (\mathbb{C}^\eta_1)^{-1}$ is below 
$\mathbb{B}^\eta_0\cup (\mathbb{B}^\eta_1)^{-1}$. The proof of Proposition \ref{squareDeltaD2} provides $h\! :\! 2^\omega\!\rightarrow\! 2^\omega$ injective continuous such that 
$\mathbb{S}^\eta_\varepsilon\!\subseteq\! (h\!\times\! h)^{-1}(\mathbb{N}^\eta_\varepsilon )$ for each $\varepsilon\!\in\! 2$. It remains to set $f(\varepsilon\alpha )\! :=\!\varepsilon h(\alpha )$.
\hfill{$\square$}

\section{$\!\!\!\!\!\!$ Negative results}

- By Theorem 15 in [L4], we cannot completely remove the assumption that $A$ is s-acyclic or locally countable in Corollary \ref{cor1pi02}. We can wonder whether there is an antichain 
basis if this assumption is removed (for this class $\bormtwo$ or any other one appearing in this section). This also shows that we cannot simply assume the disjointness of the analytic sets $A,B$ in Theorem \ref{delta02} and Corollaries \ref{cor2pi02}, \ref{partialsigma02}.\bigskip

\noindent - We can use the proof of the previous fact to get a negative result for the 
class $\bortwo$.

\begin{thm} \label{negdelta02} There is no tuple $(\mathbb{X},\mathbb{Y},\mathbb{A},\mathbb{B})$, where 
$\mathbb{X},\mathbb{Y}$ are Polish and $\mathbb{A},\mathbb{B}$ are disjoint analytic subsets of $\mathbb{X}\!\times\!\mathbb{Y}$, such that for any tuple $({\cal X},{\cal Y},{\cal A},{\cal B})$ of this type, exactly one of the following holds:\smallskip

 (a) $\cal A$ is separable from $\cal B$ by a $\mbox{pot}(\bortwo )$ set,\smallskip
 
 (b) $(\mathbb{X},\mathbb{Y},\mathbb{A},\mathbb{B})\sqsubseteq ({\cal X},{\cal Y},{\cal A},{\cal B})$.
\end{thm}

\noindent\bf Proof.\rm ~We argue by contradiction. By Lemma \ref{E00}, we get 
$(\mathbb{X},\mathbb{Y},\mathbb{A},\mathbb{B})\sqsubseteq 
(2^\omega ,2^\omega ,\lceil T\rceil\cap\mathbb{E}^0_0,\lceil T\rceil\cap\mathbb{E}^1_0)$. This shows that $\mathbb{A},\mathbb{B}$ are locally countable. As (a) and (b) cannot hold simultaneously, $\mathbb{A}$ is not separable from $\mathbb{B}$ by a $\mbox{pot}(\bortwo )$ 
set. By Corollary \ref{cor2delta02} we get 
$$(2^\omega ,2^\omega ,\lceil T\rceil\cap\mathbb{E}^0_0,\lceil T\rceil\cap\mathbb{E}^1_0)
\sqsubseteq (\mathbb{X},\mathbb{Y},\mathbb{A},\mathbb{B})\mbox{,}$$ 
so that we may assume that $(\mathbb{X},\mathbb{Y},\mathbb{A},\mathbb{B})\! =\! 
(2^\omega ,2^\omega ,\lceil T\rceil\cap\mathbb{E}^0_0,\lceil T\rceil\cap\mathbb{E}^1_0)$.\bigskip

\noindent $\bullet$ In the proof of Theorem 15 in [L4], the author considers a set 
$A\! =\!\bigcup_{s\in (\omega\setminus\{ 0\})^{<\omega}}~\mbox{Gr}({l_s}_{\vert G})$, where the 
$l_s$'s are partial continuous open maps from $2^\omega$ into itself with dense open domain, and $G$ is the intersection of their domain. Moreover, the $l_s$'s have the properties that $l_s(x)\!\not=\! l_t(x)$ if $t\!\not=\! s$, and $l_s(x)$ is the limit of $\big( l_{sk}(x)\big)_{k\in\omega}$, for each $x\!\in\! G$. We set, for $\varepsilon\!\in\! 2$, $A_\varepsilon\! :=\!\bigcup_{s\in (\omega\setminus\{ 0\})^{<\omega},
\vert s\vert\equiv\varepsilon\ (\mbox{mod }2)}~\mbox{Gr}({l_s}_{\vert G})$, so that $A_0$ and 
$A_1$ are disjoint Borel sets.\bigskip

 Let us check that $A_0$ is not separable from $A_1$ by a $\mbox{pot}(\bortwo )$ set. We argue by contradiction, which gives $D\!\in\!\mbox{pot}(\bortwo )$ and a dense $G_\delta$ subset $H$ of $2^\omega$ such that $D\cap H^2\!\in\!\bortwo (H^2)$. We may assume that $H\!\subseteq\! G$. Note that $H\cap\bigcap_{s\in (\omega\setminus\{ 0\})^{<\omega}}~l_s^{-1}(H)$ is a dense 
$G_\delta$ subset of $2^\omega$, and thus contains a point $x$. The vertical section 
$A_x$ is contained in $H$. In particular, the disjoint sections $(A_0)_x$ and $(A_1)_x$ are separable by a $\bortwo$ subset $\cal D$ of the Polish space $H$. It remains to note that 
${\cal D}\cap\overline{A_x}^H$ is a dense and co-dense $\bortwo$ subset of $\overline{A_x}^H$, which contradicts Baire's theorem.\bigskip

 This gives $u\! :\! N_0\!\rightarrow\! 2^\omega$ and $v\! :\! N_1\!\rightarrow\! 2^\omega$ with 
$\lceil T\rceil\cap\mathbb{E}^\varepsilon_0\!\subseteq\! (u\!\times\! v)^{-1}(A_\varepsilon )$.\bigskip

\noindent $\bullet$ We set $B_1\! :=\!\lceil T\rceil\cap (\mathbb{E}^0_0\cup\mathbb{E}^1_0)$. Note that $B_1\!\notin\!\mbox{pot}(G_\delta )$, since otherwise $\lceil T\rceil\cap\mathbb{E}^0_0$ and 
${\lceil T\rceil\cap\mathbb{E}^1_0}$ are two disjoint $\mbox{pot}(G_\delta )$ sets, and thus 
$\mbox{pot}(\bortwo )$-separable. Then we can follow the proof of Theorem 15 in [L4]. This proof gives 
$U\! :\! F\!\rightarrow\! G$ and $V\! :\! F\!\rightarrow\! 2^\omega$ injective continuous satisfying the inclusion ${\bigcup_{n\in\omega}~\mbox{Gr}(f_n)\!\subseteq\! (U\!\times\! V)^{-1}(A)}$.\bigskip

The only thing to check is that there is $(c,d)$ in 
$\bigcup_{n\in\omega}~\omega^n\!\times\!\omega^{n+1}$ and a nonempty open subset $R$ of 
$D_{f_{c,d}}$ such that $\Big( U(x),V\big( f_{c,d}(x)\big)\Big)\!\notin\!\mbox{Gr}(l_\emptyset )$ for each $x\!\in\! R$. We argue by contradiction, which gives a dense $G_\delta$ subset $K$ of $F$ such that $\bigcup_{n\in\omega}~\mbox{Gr}({f_n}_{\vert K})\!\subseteq\! 
(U_{\vert K}\!\times\! V)^{-1}\big(\mbox{Gr}({l_\emptyset}_{\vert G} )\big)$. As 
$(U_{\vert K}\!\times\! V)^{-1}\big(\mbox{Gr}({l_\emptyset}_{\vert G} )\big)$ is the graph of a partial Borel map, $\bigcup_{n\in\omega}~\mbox{Gr}({f_n}_{\vert K})$ too. Therefore 
$\bigcup_{n\in\omega}~\mbox{Gr}({f_n}_{\vert K})\!\in\!\mbox{pot}(\bormone )\!\setminus\!
\mbox{pot}(G_\delta )$, which is absurd.\hfill{$\square$}

\vfill\eject

 This shows that we cannot completely remove the assumption that $A\cup B$ is s-acyclic or locally countable in Corollary \ref{cor1delta02}. This also shows that we cannot simply assume the disjointness of the analytic sets $A,B$ in Theorem \ref{delt02} and Corollary \ref{cor2delta02}.\bigskip

\noindent - By Theorem 2.16 in [L3], we cannot completely remove the assumption that 
$A\cup B$ is s-acyclic or locally countable in Corollary \ref{cor1pi01}. This also shows that we cannot simply assume disjointness in Theorem \ref{checkD2Sigma01} and Corollary \ref{cor2pi01}.\bigskip

 We saw that there is a version of Corollary \ref{partialsigma02} for ${\bf\Gamma}\! =\!\boraone$, where we replace the class $F_\sigma$ with the class of open sets. We cannot replace the class $F_\sigma$ with the class of closed sets.

\begin{prop} \label{negpartialSigma01} There is no triple 
$(\mathbb{X},\mathbb{A},\mathbb{B})$, where $\mathbb{X}$ is Polish and 
$\mathbb{A},\mathbb{B}$ are disjoint analytic relations on $\mathbb{X}$ such that $\mathbb{A}$ is contained in a potentially closed s-acyclic or locally countable  relation such that, for each triple 
$({\cal X},{\cal A},{\cal B})$ of the same type, exactly one of the following holds:\smallskip

(a) the set $\cal A$ is separable from $\cal B$ by a $\mbox{pot}(\boraone )$ set,\smallskip

(b) $(\mathbb{X},\mathbb{X},\mathbb{A},\mathbb{B})\sqsubseteq ({\cal X},{\cal X},{\cal A},{\cal B})$.
\end{prop} 

\noindent\bf Proof.\rm ~We argue by contradiction, which gives a triple. Note that 
$\mathbb{A}$ is not separable from $\mathbb{B}$ by a $\mbox{pot}(\boraone )$ set. Theorem 9 in [L5] gives $F,G\! :\! 2^\omega\!\rightarrow\!\mathbb{X}$ continuous such that 
$\Delta (2^\omega )\!\subseteq\! (F\!\times\! G)^{-1}(\mathbb{A})$ and 
$\mathbb{G}_0\!\subseteq\! (F\!\times\! G)^{-1}(\mathbb{B})$. We set 
$\mathbb{A}'\! :=\! (F\!\times\! G)[\Delta (2^\omega )]$, 
$\mathbb{B}'\! :=\! (F\!\times\! G)[\mathbb{G}_0]$ and 
$\mathbb{C}'\! :=\! (F\!\times\! G)[\overline{\mathbb{G}_0}]$. Note that 
$\mathbb{A}'$, $\mathbb{C}'$ are compact and $\mathbb{C}'$ is the locally countable disjoint union of $\mathbb{A}'$ and $\mathbb{B}'$. In particular, $\mathbb{B}'$ is $D_2(\boraone )$, 
$\mathbb{A}'\!\subseteq\!\mathbb{A}$, $\mathbb{B}'\!\subseteq\!\mathbb{B}$, and 
$\mathbb{A}'$ is not separable from $\mathbb{B}'$ by a pot$(\boraone )$ set. So we may assume that $\mathbb{A},\mathbb{B}$ are Borel with locally countable union which is the closure of 
$\mathbb{B}$. Corollary \ref{cor1pi01} gives $f',g'\! :\! 2^\omega\!\rightarrow\!\mathbb{X}$ injective continuous such that 
$\mathbb{G}_0\! =\!\overline{\mathbb{G}_0}\cap (f'\!\times\! g')^{-1}(\mathbb{B})$. In particular, 
$$\Delta (2^\omega )\!\subseteq\! (f'\!\times\! g')^{-1}(\overline{\mathbb{B}}\!\setminus\!\mathbb{B})\! =\! (f'\!\times\! g')^{-1}(\mathbb{A}).$$ 
This means that we may assume that 
$\mathbb{X}\! =\! 2^\omega$, $\mathbb{A}\! =\!\Delta (2^\omega )$ and 
$\mathbb{B}\! =\!\mathbb{G}_0$.\bigskip

 The proof of Theorem 10 in [L5] provides a Borel graph $\cal B$ on $X\! :=\! 2^\omega$ with no Borel countable coloring such that any locally countable Borel digraph contained in $\cal B$ has a Borel countable coloring. Consider the closed symmetric acyclic locally countable relation 
${\cal A}\! :=\!\Delta (2^\omega )$. As there is no Borel countable coloring of $\cal B$, $\cal A$ is not separable from $\cal B$ by a pot$(\boraone )$ set. If $f,g$ exist, then $f\! =\! g$ since $\mathbb{A}$ is contained in $(f\!\times\! g)^{-1}({\cal A})$. This implies that $f$ is a homomorphism from $\mathbb{G}_0$ into $\cal B$. The digraph $(f\!\times\! f)[\mathbb{G}_0]$ is locally countable and Borel since $f$ is injective. Thus it has a Borel countable coloring, and $\mathbb{G}_0$ too, which is absurd.
\hfill{$\square$}\bigskip

 For oriented graphs, we cannot completely remove the assumption that $G$ is s-acyclic or locally countable in Theorem \ref{Zog2}. Let us check it for ${\bf\Gamma}\! =\!\bortwo$.
 
\begin{prop} There is no tuple $(\mathbb{X},\mathbb{G})$, where $\mathbb{X}$ is Polish and 
$\mathbb{G}$ is an analytic oriented graph on $\mathbb{X}$, such that for any tuple $({\cal X},{\cal G})$ of this type, exactly one of the following holds:\smallskip  

(a) the set $\cal G$ is separable from ${\cal G}^{-1}$ by a $\mbox{pot}(\bortwo )$ set,\smallskip

(b) there is $f\! :\! 2^\omega\!\rightarrow\! X$ injective continuous such that 
$\mathbb{G}\!\subseteq\! (f\!\times\! f)^{-1}({\cal G})$.\end{prop} 
 
\noindent\bf Proof.\rm ~We use the notation of the proof of Theorem \ref{negdelta02}, and argue by contradiction. Recall the analytic s-acyclic oriented graph 
${\cal G}_{\bortwo}\! =\! (\lceil T\rceil\cap\mathbb{E}^0_0)\cup (\lceil T\rceil\cap\mathbb{E}^1_0)^{-1}$ considered in the proof of Theorem \ref{Zog2}. Note that there is 
$f_0\! :\!\mathbb{X}\!\rightarrow\! 2^\omega$ injective continuous such that 
$\mathbb{G}\!\subseteq\! (f_0\!\times\! f_0)^{-1}({\cal G}_{\bortwo})$. In particular, $\mathbb{G}$ is s-acyclic and Theorem \ref{Zog2} applies. This shows that we may assume that 
$(\mathbb{X},\mathbb{G})\! =\! (2^\omega ,{\cal G}_{\bortwo})$.\bigskip
 
 If $R$ is a relation on $2^\omega$, then we set 
$G_R\! :=\!\{ (0\alpha ,1\beta )\mid (\alpha ,\beta )\!\in\! R\}$. As $A_0$ is not separable from $A_1$ by a $\mbox{pot}(\bortwo )$ set, $G_{A_0}$ is not separable from $G_{A_1}$ by a $\mbox{pot}(\bortwo )$ set. As $G_{A_0}\cup G_{A_1}\!\subseteq\! N_0\!\times\! N_1$ and $G_{A_0},G_{A_1}$ are disjoint, 
$\mathbb{H}\! :=\! G_{A_0}\cup (G_{A_1})^{-1}$ is a Borel oriented graph, and $\mathbb{H}$ is not separable from $\mathbb{H}^{-1}$ by a $\mbox{pot}(\bortwo )$ set, as in the proof of Theorem 
\ref{Zog}. If $f\! :\! 2^\omega\!\rightarrow\! 2^\omega$ is injective continuous and 
$(\lceil T\rceil\cap\mathbb{E}^0_0)\cup (\lceil T\rceil\cap\mathbb{E}^1_0)^{-1}\!\subseteq\!
\mathbb{H}$, then on a nonempty clopen set $S\! :=\! N_{s_q}\!\times\! N_{t_q}$, the first coordinate is either preserved, or changed.\bigskip

 As in the proof of Lemma \ref{E00}, we see that $\lceil T\rceil\cap\mathbb{E}^0_0\cap S$ is not separable from $\lceil T\rceil\cap\mathbb{E}^1_0\cap S$ by a $\mbox{pot}(\bortwo )$ set. By Corollary \ref{cor1delta02}, there is $f\! :\! 2^\omega\!\rightarrow\! 2^\omega$ injective continuous such that 
$$\lceil T\rceil\cap\mathbb{E}^\varepsilon_0\!\subseteq\! 
(f\!\times\! f)^{-1}(\lceil T\rceil\cap\mathbb{E}^\varepsilon_0\cap S)$$ 
for each $\varepsilon\!\in\! 2$. This proves the existence of 
$g\! :\! 2^\omega\!\rightarrow\! 2^\omega$ injective continuous such that 
$$\lceil T\rceil\cap (\mathbb{E}^0_0\cup\mathbb{E}^1_0)\!\subseteq\! (g\!\times\! g)^{-1}(G_A).$$ This gives $u\! :\! N_0\!\rightarrow\! 2^\omega$ and $v\! :\! N_1\!\rightarrow\! 2^\omega$ injective continuous such that 
$\lceil T\rceil\cap (\mathbb{E}^0_0\cup\mathbb{E}^1_0)\!\subseteq\! (u\!\times\! v)^{-1}(A)$ since the maps $\varepsilon\alpha\!\mapsto\!\alpha$ are injective. But we saw that this is not possible in the proof of Theorem \ref{negdelta02}.\hfill{$\square$}\bigskip

\noindent\bf Question.\rm\ Are there versions of our results for the classes $D_\eta (\boratwo ),\check D_\eta (\boratwo )$  (when $\omega\!\leq\!\eta\! <\!\omega_1$) and 
$\Delta\big( D_\eta (\boratwo )\big)$ (when $2\!\leq\!\eta\! <\!\omega_1$)?

\vfill\eject

\section{$\!\!\!\!\!\!$ References}

\noindent [D-SR]\ \ G. Debs and J. Saint Raymond, Borel liftings of Borel sets: 
some decidable and undecidable statements,~\it Mem. Amer. Math. Soc.\rm ~187, 876 (2007)

\noindent [K]\ \ A. S. Kechris,~\it Classical descriptive set theory,~\rm Springer-Verlag, 1995

\noindent [K-S-T]\ \ A. S. Kechris, S. Solecki and S. Todor\v cevi\'c, Borel chromatic numbers,\ \it 
Adv. Math.\rm\ 141 (1999), 1-44

\noindent [L1]\ \ D. Lecomte, Classes de Wadge potentielles et 
th\'eor\`emes d'uniformisation partielle,\it ~Fund. Math.~\rm 143 (1993), 231-258

\noindent [L2]\ \ D. Lecomte, Uniformisations partielles et crit\`eres \`a la 
Hurewicz dans le plan,~\it Trans. Amer. Math. Soc.\rm ~347, 11 (1995), 4433-4460

\noindent [L3]\ \ D. Lecomte, Tests \`a la Hurewicz dans le plan,\it ~Fund. Math.~
\rm 156 (1998), 131-165

\noindent [L4]\ \ D. Lecomte, Complexit\'e des bor\'eliens~\`a coupes 
d\'enombrables,\ \it Fund. Math.~\rm 165 (2000), 139-174

\noindent [L5]\ \ D. Lecomte, On minimal non potentially closed subsets of the plane,\ \it Topology Appl.\rm\ 154, 1 (2007), 241-262

\noindent [L6]\ \ D. Lecomte, How can we recognize potentially $\bormxi$ subsets of the plane?,~\it J. Math. Log.\rm\ 9, 1 (2009), 39-62

\noindent [L7]\ \ D. Lecomte, A dichotomy characterizing analytic graphs of uncountable Borel chromatic number in any dimension,~\it Trans. Amer. Math. Soc.\rm~361 (2009), 4181-4193

\noindent [L8]\ \ D. Lecomte, Potential Wadge classes,~\it\ Mem. Amer. Math. Soc.,\rm ~221, 1038 (2013)

\noindent [Lo1]\ \ A. Louveau, Some results in the Wadge hierarchy of Borel sets,\ \it Cabal Sem. 
79-81, Lect. Notes in Math.\ \rm 1019 (1983), 28-55

\noindent [Lo2]\ \ A. Louveau, A separation theorem for $\Ana$ sets,\ \it Trans. A. M. S.\rm\  260 (1980), 363-378 

\noindent [Lo3]\ \ A. Louveau, Ensembles analytiques et bor\'eliens dans les 
espaces produit,~\it Ast\'erisque (S. M. F.)\ \rm 78 (1980)

\noindent [Lo4]\ \ A. Louveau, Some dichotomy results for analytic graphs,~\it manuscript\ \rm 

\noindent [Lo-SR]\ \ A. Louveau and J. Saint Raymond, The strength of Borel Wadge determinacy,
\ \it Cabal Seminar 81-85, Lecture Notes in Math.\ \rm 1333 (1988), 1-30

\noindent [M]\ \ Y. N. Moschovakis,~\it Descriptive set theory,~\rm North-Holland, 1980

\end{document}